\documentclass[12pt,reqno]{amsart}
\usepackage{amsthm}
\usepackage{amssymb}
\usepackage{amsmath}
\usepackage{mathrsfs}
\usepackage{amsfonts}
\usepackage{amssymb,amsmath}
\usepackage{graphicx}
\usepackage[colorlinks, linkcolor=black, citecolor=blue,   hypertexnames=false]{hyperref}
\usepackage[numbers,sort&compress]{natbib}
\usepackage{geometry}
\usepackage{abstract}

\numberwithin{equation}{section}
\newtheorem{thm}{Theorem}[section]
\newtheorem{lem}{Lemma}[section]
\newtheorem{prop}{Proposition}[section]
\newtheorem*{pf}{Proof}
\newtheorem{rk}{Remark}[section]
\newtheorem{definition}{Definition}[section]

\theoremstyle{definition}

\allowdisplaybreaks

\begin{document}
	\title[Stable blow-up for Yang-Mills heat flow]{\textbf{Stable blow-up solutions for the SO(\MakeLowercase{d})-equivariant supercritical Yang-Mills heat flow}}
	\author{Yezhou Yi}
	
	\address{ \centerline{Yezhou Yi}
		\newline\indent
		\centerline{School of Mathematics and Statistics}
		\newline\indent
		\centerline{Ningbo University, Ningbo, 315211, Zhejiang, P.R. China}
	}
	\email{yiyezhou@nbu.edu.cn}
	\date{}
    \maketitle
	
\begin{abstract}
\noindent\textbf{ABSTRACT.} We consider the $SO(d)$-equivariant Yang-Mills heat flow 
    \begin{equation*}
		\partial_t u-\partial_r^2 u-\frac{(d-3)}{r}\partial_r u+\frac{(d-2)}{r^2}u(1-u)(2-u)=0
	\end{equation*} in dimensions $d>10.$ We construct a family of $\mathcal{C}^{\infty}$ solutions which blow up in finite time via concentration of a universal profile \begin{equation*}
	u(t,r)\sim Q\left(\frac{r}{\lambda(t)}\right),
\end{equation*}where $Q$ is a stationary state of the equation and the blow-up rates are quantized by \begin{equation*}
\lambda(t)\sim c_{u}(T-t)^{\frac{l}{\gamma}},\,\,\,l\,\,\,\text{is any positive integer},\,\,\,\gamma=\gamma(d)=\frac{d-4-\sqrt{(d-6)^2-12}}{2}.
\end{equation*}
Moreover, such solutions are in fact $(l-1)$-codimension stable under pertubation of the initial data.
\end{abstract}
~\\
\noindent\textbf{KEYWORDS:} {Yang-Mills heat flow; Asymptotic behavior; $SO(d)$-equivariant; Stability.}\\

\noindent\textbf{MSC:} 35B40, 35K15, 35K55

\section{Introduction}

We denote $E$ a principle fibre bundle over a $d$ dimensional Riemannian manifold $M,$ with a semi-simple Lie group as structure group. Denoting $AdE$ the adjoint bundle to $E,$ a smooth connection on $E$ is a smooth map from $M$ to $AdE\otimes T^*M.$ Denoting $\mathcal{G}$ as the Lie algebra of $G,$ locally a connection $A$ is a $\mathcal{G}$-valued 1-form on the coordinate patches $U_{\alpha}$ of $M,$ as $A=A_{j}(x)\,\mathrm{d}x^j$ with $A_j\,\colon\,U_{\alpha}\rightarrow \mathcal{G}.$ Denoting $D_A$ the covariant derivative with respect to $A,$ the curvature $F_A$ of a connection $A$ is defined by $F_A=D_A \circ D_A,$ locally it is $\mathcal{G}$-valued 2-form $F_{j,k}\, \mathrm{d}x^j\mathrm{d}x^k,$ where
\begin{equation*}
	F_{j,k}:=\partial_j A_k-\partial_k A_j+[A_j,A_k].
\end{equation*} 
The Yang-Mills functional $\mathcal{F}$ is defined by 
\begin{equation*}
	\mathcal{F}(A)=\int_{M} F_{j,k}F^{j,k}\,\mathrm{d}vol_{M},
\end{equation*}
which is invariant under gauge transformations. The associated Euler-Lagrange equations read 
\begin{equation}\label{Yang-Mills connections equation}
	D^jF_{j,k}=0,
\end{equation}
where $D_j=:\partial_j+[A_j,\cdot].$ Solutions of (\ref{Yang-Mills connections equation}) are referred to as Yang-Mills connections. One way to find Yang-Mills connections is to study the $L^2$-gradient flow associated with $\mathcal{F},$ i.e. the initial value problem
\begin{equation}\label{eq of gradient flow}
	\left\{\begin{aligned}
		&\partial_tA_j(t,x)=-D^kF_{j,k}(t,x)\\ &A_j(0,x)=A_{0j}(x)
	\end{aligned}\right.
\end{equation}
for some initial connection $A_0.$ Usually (\ref{eq of gradient flow}) is referred to as the Yang-Mills heat flow.\par 
In this paper, we consider (\ref{eq of gradient flow}) under the situation $M=\mathbb{R}^d,$ $G=SO(d),$ $E$ is the trivial bundle $\mathbb{R}^d\times SO(d),$ and we investigate connections given by 
\begin{equation}\label{def of sod equivariant connnection}
	A_j(x)=\frac{u(r)}{r^2}\sigma_j(x),
\end{equation}
where $r=|x|,$ $u$ is a real-valued function on $[0,\infty),$ and $\{\sigma_j\}_{j=1}^d$ are a basis for the Lie algebra $so(d),$ given by 
\begin{equation*}
	(\sigma_j)_{\beta}^{\alpha}=\delta_j^{\alpha}x^{\beta}-\delta_j^{\beta}x^{\alpha},\,\,\,\text{for}\,\,\,1\le \alpha, \beta\le d.
\end{equation*}
Note that connections satisfying (\ref{def of sod equivariant connnection}) are equivariant with respect to $SO(d)$-action, and are referred to as $SO(d)$-equivariant connections. 
In this case, (\ref{eq of gradient flow}) becomes     
\begin{equation}\label{eq:Y-M heat}
	\left\{\begin{aligned}
		&\partial_t u-\partial_r^2 u-\frac{(d-3)}{r} \partial_r u+\frac{(d-2)}{r^2}u(1-u)(2-u)=0\\ &u(0,\cdot)=u_0(\cdot)
	\end{aligned}\right. .
\end{equation}We credit Dumitrascu \cite{dumitrascu1982equivariant} for the first derivation of equivariant supercritical Yang-Mills equation, the readers can also refer to Weinkove \cite{weinkove2004singularity} for more details.\par
 Let us briefly explain the meaning of energy supercritical. For any $\lambda>0,$ if $u(t,r)$ is a solution of (\ref{eq:Y-M heat}), then $u(\frac{t}{\lambda^2},\frac{r}{\lambda})$ is also a solution. If we denote the energy functional of (\ref{eq:Y-M heat}) as 
\[E(u(t)):=\frac12 \int_{0}^{+\infty} \left( |\partial_r u|^2+\frac{(d-2)u^2(2-u)^2}{2r^2}\right)r^{d-3}\,\mathrm{d}r,\]
then we have for any radial function $u_0\,\colon\, \mathbb{R}^d\rightarrow \mathbb{R},$ 
\[E\left(u_0\left(\frac{r}{\lambda}\right)\right)=\lambda^{d-4}E(u_0(r)).\]
Therefore, $d<4$ corresponds to the energy subcritical cases, $d=4$ corresponds to energy-critical case, and $d\ge 5$ corresponds to the energy supercritical cases. \par   

Historically, there has been a lot of work devoted to the study of Yang-Mills heat flow. In the case $d=2$ or $3,$ R\aa de \cite{rade1992yang} proved the flow of (\ref{eq of gradient flow}) exists for all time and converges to a Yang-Mills connection. In the case $d=4,$ global existence of solutions of (\ref{eq:Y-M heat}) was established by Schlatter, Struwe and Tahvildar-Zadeh \cite{schlatter1998global}, and Waldron \cite{waldron2019long} for more general geometric situations of (\ref{eq of gradient flow}). In the case $d\ge 5,$ solutions of (\ref{eq:Y-M heat}) may blow up in finite time, see the works of Naito \cite{naito1994finite}, Grotowski \cite{grotowski2001finite} and Gastel \cite{gastel2002singularities}. However, they did not give any structure of blow-up solutions while this paper manages to describe it. Weinkove \cite{weinkove2004singularity} investigated the nature of singularities of Yang-Mills heat flow over a compact manifold and showed that under some assumptions of the blow-up rate, homothetically shrinking solitons appear as blow-up limits at singular points. Such objects correspond to self-similar solutions of Yang-Mills heat flow on the trivial bundle over $\mathbb{R}^d,$ which were also described explicitly in Section 4 in \cite{weinkove2004singularity} for $5\le d\le 9.$ On
Weinkove's self-similar blow-up solutions of (\ref{eq:Y-M heat}), Donninger and Sch\"orkhuber \cite{donninger2019stable}, Glogi{\'c} and Sch{\"o}rkhuber \cite{glogic2020nonlinear} proved that these blow-up are stable when $5\le d\le 9.$ For higher dimension cases $d\ge 10,$ it is proved in Bizo{\'n} and Wasserman \cite{bizon2015nonexistence} Theorem 2 that for $d\ge 10,$ there exists no self-similar blow-up solution of (\ref{eq:Y-M heat}). It is then expected that type II blow-up solutions exist.\par     
In this paper, we will construct blow up solutions of (\ref{eq:Y-M heat}) for $d>10.$ In fact, high dimensional blow-up phenomenon has been widely studied for various types of partial differential equations. For the semilinear heat equation
\begin{equation*}
	\partial_t u=\Delta u+|u|^{p-1}u 
\end{equation*}
with $d\ge 11$ and $p>1+\frac{4}{d-4-2\sqrt{d-1}}.$ Herrero and Vel{\'a}zquez \cite{herrero1994explosion} formally showed the existence of type II blow up with 
\begin{equation*}
	\|u(t)\|_{L^{\infty}}\sim \frac{1}{(T-t)^{\frac{2\alpha l}{p-1}}},\,\,\,l\in \mathbb{N}_{+},\,\,\,2\alpha l>1.
\end{equation*}
The formal result was rigorously clarified by the works of Mizoguchi \cite{mizoguchi2007rate}, Matano and Merle \cite{matano2009classification}, Collot \cite{collot2017nonradial}. For energy supercritical nonlinear Schr\"odinger equation in dimensions $d\ge 11,$ Merle, Rapha{\"e}l and Rodnianski \cite{merle2015type} constructed smooth blow-up solutions via a robust ``modulation" method. For the energy critical focusing nonlinear Schr\"odinger equation in space dimensions $d\ge 7,$ Jendrej \cite{jendrej2017construction} proved the existence of pure two bubbles, one of the bubbles develops at scale $1,$ whereas the length scale of the other converges to $0.$ Jendrej used energy-virial functional along with a well-designed approximation to operators which eliminates the unbounded part, such method has been adapted successfully in other dispersive equations, for example in energy critical wave equation when $d=6,$ see Jendrej \cite{jendrej2019construction}. For energy supercritical wave equation in $d\ge 11,$ Collot \cite{collot2018type} constructed blow-up with one bubble in the radial case, which also studied corresponding blow-up manifold and introduced an interesting Morawetz term to control local energy. For high dimensional harmonic heat flow, Ghoul, Ibrahim and Nguyen \cite{ghoul2018stability} showed one bubble blow-up exists for 1-corotational supercritical harmonic heat flow in $d\ge 7.$ Also when $d=7,$ Ghoul \cite{ghoul2017stable} proved there exists the same blow-up structure with the blow-up rate as
\begin{equation}\label{critical blow up rate}
	\lambda(t)\simeq \frac{\sqrt{(T-t)}}{|\log (T-t)|}.
\end{equation}
For energy supercritical wave maps with $d\ge 7,$ one may refer to Ghoul, Ibrahim and Nguyen \cite{ghoul2018construction}. However, whether similar blow-up phenomenon happens remains an open problem for Yang-Mills heat flow.\par  
Next we introduce the main result of this paper. Denoting $Q(r)$ as the ground state solution of (\ref{eq:Y-M heat}), 
i.e. it satisfies the equation
\begin{equation}\label{eq:Y-M heat ground state}
	\left\{ \begin{aligned}
		&-\partial_r^2 Q-\frac{(d-3)}{r}\partial_r Q+(d-2)\frac{Q(1-Q)(2-Q)}{r^2}=0\\
		&Q(0)=\partial_r Q(0)=0
	\end{aligned}\right. .
\end{equation}
In the author's previous work \cite{yi2023asymptotic}, it shows that (\ref{eq:Y-M heat ground state}) admits a solution which satisfies the asymptotics (when $d>10$) 
\begin{equation}\label{asymp: ground state}
	Q(r)=\left\{\begin{aligned}
		&\frac{1}{2}r^2+O(r^4)\,\,\,\text{as}\,\,\,r\rightarrow 0\\
		&1-\alpha r^{-\gamma}(1+O(r^{-2\gamma}))\,\,\,\text{as}\,\,\,r\rightarrow \infty
	\end{aligned}\right. ,	
\end{equation}
where
\begin{equation}\label{def of gamma}
	\alpha>0,\,\,\,\gamma=\gamma(d)=\frac{d-4-\sqrt{(d-6)^2-12}}{2}.
\end{equation}
Note that when $d>10,$ $\gamma\in (1,2).$\par 
Our goal is to study potential blow-up phenomenon of (\ref{eq:Y-M heat}), and our main result is the following.\par 
\begin{thm}\label{main thm}
	Let $d>10,$ $\gamma$ as in (\ref{def of gamma}), let $l$ be any positive integer, denoting
	\begin{equation}\label{def: def of h}
	\hbar:=\left[\frac{1}{2}\Big(\frac{d-2}{2}-\gamma\Big)\right],
	\end{equation} given $L\gg l$ a large integer and defining $\Bbbk:=L+\hbar+1.$ Then there exists a smooth radial initial data $u_0$ such that the corresponding solution to (\ref{eq:Y-M heat}) has the decomposition
	\begin{equation}\label{eq: form of the solution}
		u(t,r)=Q\left(\frac{r}{\lambda(t)}\right)+q\left(t,\frac{r}{\lambda(t)}\right),
	\end{equation}
where \begin{equation}\label{blow up rates}
	\lambda(t)=c(u_0)(T-t)^{\frac{l}{\gamma}}(1+o_{t\rightarrow T}(1))\,\,\,\text{with}\,\,\, c(u_0)>0,
\end{equation}
and \begin{equation}
	\lim\limits_{t\rightarrow T} \|\nabla^{\sigma} q(t)\|_{L^2(r^{d-3}\mathrm{d}r)}=0\,\,\,\text{for all}\,\,\,\sigma\in \left[2\hbar+4,2\Bbbk\right].
\end{equation} 
Moreover, the blow-up solution is $(l-1)$-codimension stable.
\end{thm}

\begin{rk}
	\normalfont
	Let us briefly explain the sense of $(l-1)$-codimension stable. Our initial data is of the form
	\begin{equation}\label{eq: form of initial data}
		u_0=Q_{b(0)}+q_0,
	\end{equation}
where $Q_b$ is a deformation of $Q$ and $b=(b_1,\cdots,b_L)$ corresponds to possible unstable directions in a suitable neighborhood of $Q.$ We will prove that for all $q_0\in \dot H^{\sigma}\cap \dot H^{2\Bbbk}$ small enough, for all $(b_1(0),b_{l+1}(0),\cdots,b_L(0))$ small enough, there exists a choice of unstable directions $(b_2(0),\cdots,b_{l}(0))$ such that the solution of (\ref{eq:Y-M heat}) with initial data (\ref{eq: form of initial data}) satisfies the conclusion of Theorem \ref{main thm}. This implies the constructed solution is $(l-1)$-codimension stable.   
\end{rk}\par 

\begin{rk}
	\normalfont
Again in view of Theorem 2 in Bizo{\'n} and Wasserman \cite{bizon2015nonexistence}, one is supposed to find type II blow-up for all $d\ge 10.$ However, we consider only $d>10$ and exclude the limit case $d=10$ due to technical reasons, in $d>10$ cases we have $1<\gamma<2,$ $0<\delta:=\frac{1}{2}\Big(\frac{d-2}{2}-\gamma\Big)-\hbar<1,$ and $d-2\gamma>6$ which ensure involved estimates good enough. We conjecture that it is possible to construct type II blow-up solutions for the case $d=10$, with explicitly different blow-up rates than those in (\ref{blow up rates}).       
\end{rk}\par 

It appears that Yang-Mills heat flow share similar properties to harmonic heat flow in high dimensions, thus this paper has borrowed techniques from Ghoul et al. \cite{ghoul2018stability}. Also the preprint of Bensouilah, Duong and Ghoul \cite{bensouilah2022non} which appears after the preprint of this paper, has established similar blow-up results compared to Theorem \ref{main thm} by different methods. Note that the main idea of this paper, in the author's points of view, originates from Merle, Rapha{\"e}l and Rodnianski \cite{merle2015type} which seems astonishing since they studied blow-up for supercritical nonlinear Schr\"odinger equation. However, there are technical difficulties in order to get appropriate energy estimates which are supposed to decay well enough to close bootstrap after integration, especially when dealing with nonlinear term for energy estimates, those may be considered as the original part of this paper.\par 
A closely related topic is the blow-up behavior for hyperbolic version of (\ref{eq:Y-M heat}) as
\begin{equation}\label{eq: Yang-Mills hyperbolic}
	\partial_t^2 u-\partial_r^2 u-\frac{(d-3)}{r} \partial_r u+\frac{(d-2)}{r^2}u(1-u)(2-u)=0,
\end{equation}
where we omit more general geometric backgound. Historically, the existence of blow-up for (\ref{eq: Yang-Mills hyperbolic}) was first proved by Cazenave, Shatah and Tahvildar-Zadeh in \cite{cazenave1998harmonic}, they constructed singular traveling waves by using
self-similar solutions. The self-similar blow-up for (\ref{eq: Yang-Mills hyperbolic}) is proved to be stable for all odd dimensions $d\ge 5,$ see Donninger \cite{donninger2014stable} and Glogi{\'c} \cite{glogic2021stable}. Rapha{\"e}l and Rodnianski \cite{raphael2012stable} constructed stable one bubble blow-up when $d=4,$ i.e. the energy-critical case. Also in dimension four, Jendrej \cite{jendrej2019construction} constructed two bubbles, and Krieger, Schlag, Tataru \cite{krieger2009renormalization} showed the existence of a family of one bubble where the blow-up rates are a modification of the self-similar rate by a power of logarithm. We conjecture that one may establish similar results to Theorem \ref{main thm} for (\ref{eq: Yang-Mills hyperbolic}) with exactly the same blow-up rates $\lambda.$\par 

This paper is organized as follows. In section \ref{On the linearized operator L}, we show fundamental calculations about the linearized operator of (\ref{eq:Y-M heat}) and establish coercivity properties which are crucial for both modulation estimates and energy estimates later. In section \ref{Construction of an approximate solution}, we construct an approximate solution and estimate error terms. In section \ref{Linearization bk for k from 1 to l}, one perturbs the modulation parameter equation near a set of explicit solutions (The author doubts that whether one can find explicitly non-polynomial exact solutions or first approximate solutions to further perturb with and succeed in closing the whole bootstrap argument). In later sections, it is a somewhat standard procedure in modulation methods, we decompose the solution, describe initial data, provide bootstrap assumptions, then under such assumptions, we estimate modulation parameters by taking $L^2$ inner product of the ``orthogonality" elements with the equation of the residue which also arises in previous decomposition, then derive the most crucial monotone estimates called energy estimates so that integrating it back helps us close bootstrap (since almost all bootstrap assumptions can be estimated through proper energy) and finally a contradictory argument relying on basic topological facts concludes the main Proposition, that is, Proposition \ref{existence od sol trapped for large rescaled time}.      

\subsection*{Notations}\label{Structure of the ground state and notation}
Now we introduce some notations. Denoting $\Lambda Q(y):=y\partial_y Q(y).$ By (\ref{asymp: ground state}), we have
\begin{equation}\label{asymp: Lambda Q}
	\Lambda Q(y)=\left\{ \begin{aligned}
		&y^2+O(y^4)\,\,\,\text{as}\,\,\, y\rightarrow 0\\
		&\frac{\alpha \gamma}{y^{\gamma}}\left(1+O\left(\frac{1}{y^{2\gamma}}\right)\right)\,\,\,\text{as}\,\,\, y\rightarrow \infty
	\end{aligned}\right. .
\end{equation}
Denoting the linearized operator $\mathscr{L}:=-\partial_y^2-\frac{(d-3)}{y}\partial_y +\frac{Z(y)}{y^2},$ where $Z(y):=(d-2)f'(Q(y))$ and $f(u):=u(1-u)(2-u).$ Note that by substituting $Q_\lambda(r):=Q(\frac{r}{\lambda})$ into (\ref{eq:Y-M heat ground state}) in the variable $y=\frac{r}{\lambda}$ and acting on $\partial_{\lambda}|_{\lambda=1}$, one gets $\mathscr{L}(\Lambda Q)=0.$ Furthermore, denoting $T_k:=(-1)^k (\mathscr{L}^{-1})^k (\Lambda Q),$ for $0\le k\le L.$\par 
For any two radial functions $f_1$ and $f_2$, denoting their inner product as \begin{equation*}
	\langle f_1,f_2 \rangle:=\int_0^{\infty} f_1(y)f_2(y)y^{d-3}\,\mathrm{d}y.
\end{equation*}
For convenience, let $\int f_1:=\int_{0}^{\infty} f_1(y) y^{d-3}\,\mathrm{d}y.$\par 
Defining $\chi$ as a smooth radial cut-off function such that $\chi(y)=1$ for $0\le y\le 1,$ $\chi(y)=0$ for $y\ge 2$ and  $0<\chi(y)<1$ for $1<y<2.$ Then we denote $\chi_M(y):=\chi(\frac{y}{M}).$\par 
For any smooth radial function $g,$ denoting $g_{2k}:=\mathscr{L}^k g,$ $g_{2k+1}:=\mathscr{A}\mathscr{L}^k g,$ for any $k\in \mathbb{N}.$ Denoting $\mathscr{L}_{\lambda}:=-\partial_r^2-\frac{(d-3)}{r}\partial_r+\frac{Z_{\lambda}(r)}{r^2}$ and etc, then we write $g_{2k}^*:=\mathscr{L}_{\lambda}^k g,$ $g_{2k+1}^*:=\mathscr{A}_{\lambda}\mathscr{L}_{\lambda}^k g.$\par 
For any $b_1>0,$ we define $B_0:=b_1^{-\frac{1}{2}},$ $B_1:=B_0^{1+\eta},$ where $0<\eta\ll 1$ is to be choosen later.\par 
We emphasize again that for $\gamma$ as (\ref{def of gamma}) and $\hbar$ as (\ref{def: def of h}), we define $\delta:=\frac{1}{2}\Big(\frac{d-2}{2}-\gamma\Big)-\hbar$ throughout this paper.\par 
The notation $A\lesssim B$ means that there exists a positive constant $C$ such that $A\le CB.$ And we denote $A\simeq B$ if $A\lesssim B$ and $B\lesssim A$ are both valid.

\vspace{\baselineskip}
\section{On the linearized operator $\mathscr{L}$}\label{On the linearized operator L}
In this section, we make preparations related to $\mathscr{L}.$ We shall omit some detailed proofs if they are obtained by direct computations.

\subsection{Decomposition, kernel and computation of the inverse of $\mathscr{L}$}
\begin{lem}\label{factorization of L}
	The operator $\mathscr{L}$ factorizes as $\mathscr{L}$=$\mathscr{A^*}$$\mathscr{A}$ with
	\begin{align}
		&\mathscr{A}\omega:=\left(-\partial_y+\frac{V(y)}{y}\right)\omega=-\Lambda Q\partial_y \left(\frac{\omega}{\Lambda Q}\right),\label{def: def of A}\\
		&\mathscr{A^*}\omega:=\left(\partial_y +\frac{d-3+V(y)}{y}\right)\omega=\frac{1}{y^{(d-3)} \Lambda Q}\partial_y (y^{(d-3)} \Lambda Q \omega),\label{def: def of A*}
	\end{align}
where  
\begin{equation}\label{def: def of V}
	V(y):=\Lambda \ln(\Lambda Q)=\left\{\begin{aligned}
		&2+O(y^2)\,\,\,\text{as}\,\,\,y\rightarrow 0\\
		&-\gamma+O\left(\frac{1}{y^{2\gamma}}\right)\,\,\,\text{as}\,\,\,y\rightarrow \infty
	\end{aligned}\right. .
\end{equation}	
\end{lem}
\begin{rk}
	\normalfont
	Note that
	\begin{equation}\label{commutator between L and Lambda}
		\left[\mathscr{L},\Lambda\right]=2\mathscr{L}-\frac{\Lambda Z(y)}{y^2}.
	\end{equation}
	Denoting $\widetilde{\mathscr{L}}:=\mathscr{A}\mathscr{A^*}.$ By Lemma \ref{factorization of L},
	\begin{equation}\label{exp of tilde L}
		\widetilde{\mathscr{L}}=-\partial_y^2-\frac{(d-3)}{y}\partial_y+\frac{\widetilde{Z}(y)}{y^2},
	\end{equation}
	where 
	\begin{equation}\label{def: def of tilde Z}
		\widetilde{Z}(y):=(V+1)^2+(d-4)(V+1)-\Lambda V.
	\end{equation}
\end{rk}
Next we find the other kernel of $\mathscr{L}.$ If $\mathscr{L}\Gamma=0,$ then $\mathscr{A}\Gamma$ lies in the kernel of $\mathscr{A^*}.$ By (\ref{def: def of A*}), $\mathscr{A}\Gamma\in \textnormal{Span}\{\frac{1}{y^{d-3}\Lambda Q}\}.$ By definition (\ref{def: def of A}), we can impose that 
\[-\partial_y \Gamma+\frac{V(y)}{y}\Gamma=\frac{-1}{y^{d-3}\Lambda Q}.\]
It has a solution of the form
\[\Gamma(y)=\Lambda Q(y)\int_{1}^{y}\frac{\mathrm{d}\xi}{\xi^{d-3}(\Lambda Q(\xi))^2},\]
and the asymptotics
\begin{equation}\label{asymp: Gamma}
	\Gamma(y)\simeq\left\{\begin{aligned}
		&\frac{c}{y^{d-2}}\,\,\,\text{as}\,\,\,y\rightarrow 0\\
		&\frac{c}{y^{d-4-\gamma}}\,\,\,\text{as}\,\,\,y\rightarrow \infty
	\end{aligned}\right. .
\end{equation}\par 
Then we introduce calculations on the inverse of $\mathscr{L}.$ By standard ODE theory, for any radial function $g$, there exists a solution to $\mathscr{L}\omega=g$ written as  
\begin{equation}\label{standard inverse of L}
	\omega=\mathscr{L}^{-1}g=-\Gamma(y)\int_{0}^{y} g(x)\Lambda Q(x)x^{d-3}\,\mathrm{d}x+\Lambda Q(y)\int_{0}^{y} g(x)\Gamma(x)x^{d-3}\,\mathrm{d}x.
\end{equation}
For the convenience of calculation, there is a two-step method to compute the above $\mathscr{L}^{-1}$ as follows.
\begin{lem}\label{two step method to calculate the inverse of L}
	Let $g\in \mathcal{C}_{rad}^{\infty},$ then $\mathscr{L}\omega=g$ exists a solution solved by 
	\begin{align*}
		\mathscr{A}\omega&=\frac{1}{y^{d-3}\Lambda Q}\int_{0}^{y} g(x)\Lambda Q(x)x^{d-3}\,\mathrm{d}x,\\
		\omega&=-\Lambda Q\int_{0}^{y} \frac{\mathscr{A}\omega(x)}{\Lambda Q(x)}\,\mathrm{d}x.
	\end{align*}
\end{lem}
\begin{pf}
	\normalfont
	Acting $\mathscr{A}$ on (\ref{standard inverse of L}) and making use of (\ref{def: def of A}).
\end{pf}
\vspace{\baselineskip}
\subsection{Coercivity of $\mathscr{L}$}
The proof of the following Hardy inequality is similar to Lemma B.1 in \cite{merle2015type}.
\begin{lem}\label{hardy ineq}
	Let $\alpha>0,$ $\alpha\neq \frac{d-4}{2}$ and $u\in \mathcal{D}_{rad}:=\{u\in \mathcal{C}_{c}^{\infty}\,\text{with radial symmetry}\},$ then
	\begin{equation*}
		\int_{1}^{\infty} \frac{|\partial_y u|^2}{y^{2\alpha}}\ge \Big(\frac{d-(2\alpha+4)}{2}\Big)^2 \int_{1}^{\infty} \frac{u^2}{y^{2+2\alpha}}-C_{\alpha,d}u^2(1).
	\end{equation*}
\end{lem}
Then one can follow the road map in Appendix A in \cite{ghoul2018stability}. Firstly, we establish coercivity of $\mathscr{A^*}$ as follows.
\begin{lem}\label{coercivity of A*}
	Let $\alpha\ge 0.$ There exists $C_{\alpha}>0$ such that for all $u\in \mathcal{D}_{rad},$ $i=0,$ $1,$ $2,$
	\begin{equation*}
		\int \frac{|\mathscr{A^*}u|^2}{y^{2i}(1+y^{2\alpha})}\ge C_{\alpha} \left(\int \frac{|\partial_y u|^2}{y^{2i}(1+y^{2\alpha})}+\int \frac{u^2}{y^{2i+2}(1+y^{2\alpha})}\right).
	\end{equation*} 
\end{lem} 

Denoting \begin{equation}\label{def: def of PhiM}
	\Phi_M:=\sum\limits_{k=0}^L c_{k,M}\mathscr{L}^k(\chi_M \Lambda Q),
\end{equation} 
where 
\begin{equation}\label{def: def of ckM}
	c_{0,M}:=1,\,\,c_{k,M}:=(-1)^{k+1}\cdot \frac{\sum\limits_{j=0}^{k-1}c_{j,M}\langle \mathscr{L}^j(\chi_M \Lambda Q),T_k\rangle}{\langle \chi_M \Lambda Q,\Lambda Q\rangle},\,\,1\le k\le L.
\end{equation}
Note that the choices of $c_{k,M}$ are equivalent to 
\begin{equation}\label{property of Phi}
	\left\{\begin{aligned}
		\langle \Phi_M,\Lambda Q\rangle&=\langle \chi_M \Lambda Q,\Lambda Q\rangle\\
		\langle \Phi_M,T_k\rangle&=0\,\,\,\text{for}\,\,\,1\le k\le L
	\end{aligned}\right. .
\end{equation} 
In particular,
\begin{equation}\label{more property of Phi}
	\langle \mathscr{L}^i T_k,\Phi_M\rangle=(-1)^k \langle \chi_M \Lambda Q,\Lambda Q\rangle \delta_{i,k},\,\,\,\text{for}\,\,\,0\le i,k\le L,
\end{equation}
where \begin{equation*}
	\delta_{i,k}:=\left\{\begin{aligned}
		&1\,\,\, \text{if}\,\,\, i=k\\
		&0\,\,\,\text{if}\,\,\, i\neq k
	\end{aligned}\right. .
\end{equation*}

Then we establish coercivity of $\mathscr{A}$ as follows.
\begin{lem}\label{coercivity of A}
	Let $p\ge 0,$ $i=0,$ $1,$ $2$ and $2i+2p-(d-2\gamma-4)\neq 0.$ Assuming in addition $\langle u,\Phi_M\rangle=0,$ if $2i+2p>d-2\gamma-4.$ Then we have
	\begin{equation*}
		\int \frac{|\mathscr{A}u|^2}{y^{2i}(1+y^{2p})}\gtrsim \int \frac{|\partial_y u|^2}{y^{2i}(1+y^{2p})}+\int \frac{u^2}{y^{2i+2}(1+y^{2p})}.
	\end{equation*}
\end{lem}
Next we are in place to establish coercivity of $\mathscr{L}$ as follows.
\begin{lem}\label{coercivity of L}
	Let $k\in \mathbb{N},$ $i=0,$ $1,$ $2$ and $M=M(k)$ large enough. Then there exists $c_{M,k}>0$ such that for all $u\in \mathcal{D}_{rad}$ with $\langle u,\Phi_M\rangle=0$ if $2i+2k>d-2\gamma-6,$ we have
	\begin{equation*}
		\int \frac{|\mathscr{L}u|^2}{y^{2i}(1+y^{2k})}\ge c_{M,k}\int \left(\frac{|\partial_y^2 u|^2}{y^{2i}(1+y^{2k})}+\frac{|\partial_y u|^2}{y^{2i}(1+y^{2k+2})}+\frac{u^2}{y^{2i+2}(1+y^{2k+2})}\right),
	\end{equation*}
and
\begin{equation*}
	\int \frac{|\mathscr{L}u|^2}{y^{2i}(1+y^{2k})}\ge c_{M,k}\int \left(\frac{|\mathscr{A}u|^2}{y^{2i+2}(1+y^{2k})}+\frac{u^2}{y^{2i}(1+y^{2k+4})}\right).
\end{equation*}
\end{lem} 
Finally we give the coercivity property of iterate of $\mathscr{L}$ as follows.
\begin{lem}\label{coercivity of iterate of L}
	Let $k\in \mathbb{N},$ $M=M(k)$ large enough. Then there exists $c_{M,k}>0$ such that for any $u\in \mathcal{D}_{rad}$ with $\langle u,\mathscr{L}^m \Phi_M\rangle=0,$ $0\le m\le k-\hbar,$ we have
	\begin{align*}
	\mathscr{E}_{2k+2}(u):&=\int |\mathscr{L}^{k+1} u|^2\\ &\ge c_{M,k}\left(\sum\limits_{j=0}^k\int \frac{|\mathscr{L}^j u|^2}{y^4(1+y^{4(k-j)})}+\int \frac{|\mathscr{A}(\mathscr{L}^k u)|^2}{y^2}+\sum\limits_{j=0}^{k-1}\int \frac{|\mathscr{A}(\mathscr{L}^j u)|^2}{y^6(1+y^{4(k-j-1)})}\right).
	\end{align*}
\end{lem}
\begin{rk}
	\normalfont
	We point out that in Lemma \ref{coercivity of iterate of L}, when verifying the case for $k=0,$ there is no need for orthogonal conditions, one just apply Lemma \ref{coercivity of A*} and Lemma \ref{coercivity of A} to get 
	\begin{equation*}
		\int |\mathscr{L} u|^2\gtrsim \int \frac{|\mathscr{A}u|^2}{y^2}\gtrsim \int \frac{u^2}{y^4}+\int \frac{|\mathscr{A}u|^2}{y^2}.
	\end{equation*}
Also note that when $d>10,$ $d-2\gamma-4>2$ holds, thus the assumption in Lemma \ref{coercivity of A} meets.
\end{rk}

\vspace{\baselineskip}
\subsection{Leibniz rule for the iteration of $\mathscr{L}$}
We introduce Leibniz rules for $\mathscr{L}^k$ and $\mathscr{A}\mathscr{L}^k$ as follows. One can prove it by induction on $k,$ a similar detailed proof is given in Lemma C.1 in \cite{ghoul2018stability}.
\begin{lem}\label{leibniz rule}
	For any smooth radial function $\phi,$ $g$ and any $k\in \mathbb{N},$ we have
	\begin{align}
		\mathscr{L}^{k+1}(\phi g)&=\sum\limits_{m-=0}^{k+1} g_{2m}\phi_{2k+2,2m}+\sum\limits_{m=0}^k g_{2m+1}\phi_{2k+2,2m+1},\label{leibniz for iterate of L}\\
		\mathscr{A}\mathscr{L}^k(\phi g)&=\sum\limits_{m=0}^k g_{2m+1}\phi_{2k+1,2m+1}+\sum\limits_{m=0}^k g_{2m}\phi_{2k+1,2m},\label{leibniz for A composite iterate of L}
	\end{align}
where for $k=0,$ 
\begin{align*}
	\phi_{1,0}:&=-\partial_y \phi,\,\,\,\phi_{1,1}:=\phi,\\
	\phi_{2,0}:&=-\partial_y^2 \phi-\frac{(d-3+2V)}{y}\partial_y \phi,\,\,\,\phi_{2,1}:=2\partial_y \phi,\,\,\,\phi_{2,2}:=\phi,
\end{align*}
for $k\ge 1,$
\begin{align*}
	\phi_{2k+1,0}:&=-\partial_y \phi_{2k,0},\\
	\phi_{2k+1,2i}:&=-\partial_y \phi_{2k,2i}-\phi_{2k,2i-1},\,\,\,1\le i\le k,\\
	\phi_{2k+1,2i+1}:&=\phi_{2k,2i}+\frac{(d-3+2V)}{y}\phi_{2k,2i+1}-\partial_y \phi_{2k,2i+1},\,\,\,0\le i\le k-1,\\
	\phi_{2k+1,2k+1}:&=\phi_{2k,2k}=\phi,\\
	\phi_{2k+2,0}:&=\partial_y \phi_{2k+1,0}+\frac{(d-3+2V)}{y}\phi_{2k+1,0},\\
	\phi_{2k+2,2i}:&=\phi_{2k+1,2i-1}+\partial_y \phi_{2k+1,2i}+\frac{(d-3+2V)}{y}\phi_{2k+1,2i},\,\,\,1\le i\le k,\\
	\phi_{2k+2,2i+1}:&=\partial_y\phi_{2k+1,2i+1}-\phi_{2k+1,2i},\,\,\,0\le i\le k,\\
	\phi_{2k+2,2k+2}:&=\phi_{2k+1,2k+1}=\phi.
\end{align*}
\end{lem}

\vspace{\baselineskip}
\section{Construction of an approximate solution}\label{Construction of an approximate solution}
Denoting the self-similar change of variables as
\begin{equation}\label{change of variable}
	\omega(s,y):=u(t,r),\,\,\,y:=\frac{r}{\lambda(t)},\,\,\,s:=s_0+\int_{0}^{t} \frac{\mathrm{d}\tau}{\lambda^2(\tau)}.
\end{equation}
We shall derive more information on the parameter $\lambda(t)$ later. Substituting (\ref{change of variable}) into (\ref{eq:Y-M heat}), we get the renormalized flow
\begin{equation}\label{eq: renormalized flow}
	\partial_s \omega-\partial_y^2 \omega-\frac{(d-3)}{y}\partial_y \omega-\frac{\lambda_s}{\lambda}\Lambda \omega+\frac{(d-2)}{y^2}\omega(1-\omega)(2-\omega)=0.
\end{equation}
In this section, we construct approximate solutions with respect to (\ref{eq: renormalized flow}).
\subsection{Definition and properties of degree and homogeneous admissable functions}
Firstly, we shall sum up function properties which we will encounter often and express them in a systematic and unified way.
\begin{definition}
	Admissible function:
	we say $g\in\mathcal{C}_{rad}^{\infty}$ is admissible of degree $(p_1,p_2)\in \mathbb{N}\times \mathbb{Z}$ if\\
	\textnormal{(\romannumeral1)} For $y$ close to $0$, $g(y)=\sum\limits_{k=p_1}^p c_k y^{2k+2}+O(y^{2p+4}).$\\
	 \textnormal{(\romannumeral2)} For $y\ge 1,$ for all $k\in \mathbb{N},$ $|\partial_y^k g(y)|\lesssim y^{2p_2-\gamma-k}.$\\
	 We abbreviate it as $g\sim (p_1,p_2).$ 
\end{definition}
Under certain operations, the degree has the following properties. One just apply Lemma \ref{two step method to calculate the inverse of L} and use induction method, we shall omit the details.
\begin{lem}\label{property of degree}
	Let $g$ be an admissible function of degree $(p_1,p_2)
	\in \mathbb{N}\times \mathbb{Z},$ then\\
	\textnormal{(\romannumeral1)} $\Lambda g\sim (p_1,p_2).$\\
	\textnormal{(\romannumeral2)} $\mathscr{L}g\sim (p_1-1,p_2-1),$ for $p_1\ge 1.$\\
	\textnormal{(\romannumeral3)} $\mathscr{L}^{-1}g\sim (p_1+1,p_2+1).$\\
    \textnormal{(\romannumeral4)} $T_k\sim (k,k),$ for all $k\in \mathbb{N}.$\\
    \textnormal{(\romannumeral5)} $\Lambda T_k-(2k-\gamma)T_k\sim (k,k-1),$ for all $k\in \mathbb{N}_{+}.$
\end{lem}

\begin{definition}
	Homogeneous admissible function: Let $L\gg 1$ be an integer and $m:=(m_1,\cdots,m_L)\in \mathbb{N}^L,$ $b:=(b_1,\cdots,b_L).$ We say that a radial function $g(b,y)$ is homogeneous of degree $(p_1,p_2,p_3)\in \mathbb{N}\times\mathbb{Z}\times\mathbb{N},$ if it is a finite linear combination of monomials $\widetilde{g}(y)\prod\limits_{k=1}^L b_k^{m_k}$ with $\widetilde{g}(y)\sim (p_1,p_2)$ and $\sum\limits_{k=1}^L km_k=p_3.$ We abbreviate it as $g\sim (p_1,p_2,p_3).$ 
\end{definition}

\subsection{Estimate of the approximate profile and error terms}
\begin{prop}\label{first approximation}
	Let $d>10$ and $L\gg 1$ be an integer. Then there exists a small enough universal constant $b^{*}>0$ such that the following holds true. Let $b=(b_1,\cdots,b_L)\,\colon\,[s_0,s_1]\rightarrow (-b^*,b^*)^L$ be a $\mathcal{C}^1$ map with a priori bound on $[s_0,s_1]:$
	\begin{equation}\label{assump on b_k}
		0<b_1<b^*,\,\,\,|b_k|\lesssim b_1^k,\,\,\,\text{for}\,\,\,2\le k\le L.
	\end{equation} 
Then there exist profiles $S_1=0,$ $S_k=S_k(b,y),$ $2\le k\le L+2$ such that 
\begin{equation}\label{form of first approxiamtion}
	Q_{b(s)}(y):=Q(y)+\sum\limits_{k=1}^L b_k(s)T_k(y)+\sum\limits_{k=2}^{L+2} S_k(b,y)=:Q(y)+\Theta_{b(s)}(y)
\end{equation} 
as an approximation to (\ref{eq: renormalized flow}) satisfies
\begin{equation}\label{eq: eq of first approximation}
	\partial_s Q_b-\partial_y^2 Q_b-\frac{(d-3)}{y}\partial_y Q_b+b_1 \Lambda Q_b+\frac{(d-2)}{y^2}f(Q_b)=Mod(t)+\Psi_b
\end{equation}
with the following properties.\\
\textnormal{(\romannumeral1)} $Mod(t)=\sum\limits_{k=1}^L [(b_k)_s+(2k-\gamma)b_1b_k-b_{k+1}][T_k+\sum\limits_{j=k+1}^{L+2} \frac{\partial S_j}{\partial b_k}].$\\
\textnormal{(\romannumeral2)} The tails $S_k$ satisfy \begin{align}
	&S_k\sim (k,k-1,k),\,\,\, \text{for}\,\,\,2\le k\le L+2.\label{degree of S_k}\\ &\frac{\partial S_k}{\partial b_m}=0,\,\,\,\text{for}\,\,\,2\le k\le m\le L+2.\label{dependence for S_k with respect to b_m}
\end{align}
\textnormal{(\romannumeral3)} For all $0\le m\le L,$
\begin{equation}\label{esti of Psib in the scale B1}
	\int_{y\le 2B_1} |\mathscr{L}^{\hbar+m+1}\Psi_b|^2+\int_{y\le 2B_1} \frac{|\Psi_b|^2}{1+y^{4(\hbar+m+1)}}\lesssim b_1^{2m+4+2(1-\delta)-C_L \eta}.
\end{equation}
For all $M\ge 1,$
\begin{equation}\label{esti of Psib in the scale M}
	\int_{y\le 2M} |\mathscr{L}^{\hbar+m+1} \Psi_b|^2\lesssim M^{C} b_1^{2L+6}.
\end{equation}
\end{prop}

\begin{pf}
	\normalfont
	Defining the approximate solution to (\ref{eq: renormalized flow}) as (\ref{form of first approxiamtion}). In addition, assuming (\ref{assump on b_k}), $S_1=0$ and (\ref{dependence for S_k with respect to b_m}). Then we shall construct $S_k$ and verify (\ref{dependence for S_k with respect to b_m}) holds indeed. Applying (\ref{eq:Y-M heat ground state}) we get
	\begin{align*}
		&\partial_s Q_b-\partial_y^2 Q_b-\frac{(d-3)}{y}\partial_y Q_b+b_1 \Lambda Q_b+\frac{(d-2)}{y^2}f(Q_b)\\=&\partial_s \Theta_{b}+\mathscr{L}\Theta_{b}+b_1\Lambda Q+b_1\Lambda \Theta_{b}\\&+\frac{(d-2)}{y^2}[f(Q+\Theta_{b})-f(Q)-f'(Q)\Theta_{b}]\\ =:&A_1+A_2.
	\end{align*}
Direct computation gives
\begin{align*}
	A_1&=\sum\limits_{k=1}^L [(b_k)_s+(2k-\gamma)b_1b_k-b_{k+1}]T_k\\ &+\sum\limits_{k=1}^L [b_1b_k(\Lambda T_k-(2k-\gamma)T_k)+b_1\Lambda S_k]+b_1\Lambda S_{L+1}+b_1\Lambda S_{L+2}\\ &+\sum\limits_{k=1}^{L+1} \mathscr{L}S_{k+1}+\sum\limits_{k=2}^{L+2}\partial_s S_k.
\end{align*}
Note that 
\begin{equation*}
	\partial_s S_k=\sum\limits_{j=1}^L [(b_j)_s+(2j-\gamma)b_1b_j-b_{j+1}]\frac{\partial S_k}{\partial b_j}-\sum\limits_{j=1}^L [(b_j)_s+(2j-\gamma)b_1b_j-b_{j+1}]\frac{\partial S_k}{\partial b_j},
\end{equation*}
hence 
\begin{equation*}
	A_1=Mod(t)+\sum\limits_{k=1}^{L+1}[\mathscr{L}S_{k+1}+E_k]+E_{L+2},
\end{equation*}
where for $k=1,\cdots,L,$ \begin{equation*}
	E_k:=b_1b_k[\Lambda T_k-(2k-\gamma)T_k]+b_1\Lambda S_k-\sum\limits_{j=1}^{k-1}[(2j-\gamma)b_1b_j-b_{j+1}]\frac{\partial S_k}{\partial b_j},
\end{equation*}
for $k=L+1,$ $L+2,$ 
\begin{equation*}
	E_k:=b_1\Lambda S_k-\sum\limits_{j=1}^L[(2j-\gamma)b_1b_j-b_{j+1}]\frac{\partial S_k}{\partial b_j}.
\end{equation*}
For the expansion of $A_2,$ by Taylor expansion with integral remainder, one gets
\begin{equation*}
	A_2=\frac{(d-2)}{y^2}\left[\sum\limits_{i=2}^{L+2}P_i+R_1+R_2\right],
\end{equation*}
where 
\begin{align*}
	P_i&:=\sum\limits_{j=2}^{L+2}\frac{f^{(j)}(Q)}{j!}\sum_{\substack{|J|_1=j\\|J|_2=i}} c_J \prod_{k=1}^{L}b_k^{i_k}T_k^{i_k}\prod_{k=2}^{L+2}S_k^{j_k},\\
	R_1&:=\sum_{j=2}^{L+2}\frac{f^{(j)}(Q)}{j!}\sum_{\substack{|J|_1=j\\|J|_2\ge L+3}}c_J \prod_{k=1}^{L}b_k^{i_k}T_k^{i_k}\prod_{k=2}^{L+2}S_k^{j_k},\\
	R_2&:=\frac{\Theta_{b}^{L+3}}{(L+2)!}\int_{0}^{1}(1-\tau)^{L+2}f^{(L+3)}(Q+\tau \Theta_{b})\,\mathrm{d}\tau,
\end{align*}
with $J=(i_1,\cdots,i_L,j_2,\cdots,j_{L+2})\in \mathbb{N}^{2L+1}$ and 
\begin{equation*}
	|J|_1=\sum_{k=1}^{L}i_k+\sum_{k=2}^{L+2}j_k,\,\,\,|J|_2=\sum_{k=1}^{L}ki_k+\sum_{k=2}^{L+2}kj_k.
\end{equation*}
Note that we take $L$ large enough so that $R_2=0.$
Thus 
\begin{equation}\label{zheng he bi jin jie yu renormalized flow}
	\partial_s Q_b-\partial_y^2 Q_b-\frac{(d-3)}{y}\partial_y Q_b+b_1\Lambda Q_b+\frac{(d-2)}{y^2}f(Q_b)=Mod(t)+\Psi_b,
\end{equation} 
with 
\begin{equation}\label{expression of Psi_b}
	\Psi_b:=\sum_{k=1}^{L+1}[\mathscr{L}S_{k+1}+E_k+\frac{(d-2)}{y^2}P_{k+1}]+E_{L+2}+\frac{(d-2)}{y^2}R_1.
\end{equation}\par 
Motivated by (\ref{expression of Psi_b}), we define $\{S_k\}_{k=1}^{L+2}$ as 
\begin{equation*}
	\left\{\begin{aligned}
		S_1&=0\\ S_k&=-\mathscr{L}^{-1}F_k
	\end{aligned}\right.
\end{equation*}
with 
\begin{equation*}
	F_k:=E_{k-1}+\frac{(d-2)}{y^2}P_k,\,\,\,2\le k\le L+2.
\end{equation*}
Then we aim at proving (\ref{degree of S_k}) and (\ref{dependence for S_k with respect to b_m}). Claim: 
\begin{equation}\label{degree of F_k}
	F_k\sim (k-1,k-2,k)\,\,\,\text{and}\,\,\,\frac{\partial F_k}{\partial b_m}=0,\,\,\,2\le k\le m\le L+2.
\end{equation}
One can prove it by induction.\par 
When $k=2,$ note that by (\ref{asymp: ground state}) and \textnormal{(\romannumeral4)} of Lemma \ref{property of degree}, \begin{equation*}
	\frac{f^{(2)}(Q)}{y^2}T_1^2\lesssim \left\{\begin{aligned}
		&y^6\ll y^4\,\,\,\text{as}\,\,\,y\rightarrow 0\\
		&y^{2-3\gamma}\ll y^{-\gamma}\,\,\,\text{as}\,\,\,y\rightarrow \infty
	\end{aligned}\right. .
\end{equation*}
Combined with \textnormal{(\romannumeral5)} of Lemma \ref{property of degree}, we have 
\begin{equation*}
	F_2=b_1^2\left(\Lambda T_1-(2-\gamma)T_1+c\frac{f^{(2)}(Q)}{y^2}T_1^2\right)\sim (1,0,2).
\end{equation*}\par 
Then we show $\le k\Longrightarrow k+1,$ specifically speaking we need to prove
\begin{equation}\label{need to prove on F_k}
	F_{k+1}\sim (k,k-1,k+1)\,\,\,\text{and}\,\,\,\frac{\partial F_{k+1}}{b_m}=0,\,\,\,k+1\le m.
\end{equation}
By induction hypothesis, $F_j\sim (j-1,j-2,j)$ and $\frac{\partial F_j}{\partial b_m}=0,$ $j\le m,$ then by \textnormal{(\romannumeral3)} of Lemma \ref{property of degree}, 
\begin{equation}\label{under induction hypo, esti of S_k}
	S_j\sim (j,j-1,j)\,\,\,\text{and}\,\,\,\frac{\partial S_j}{\partial b_m}=0,\,\,\,j\le m,\,\,\,\text{for any}\,\,\,2\le j\le k.
\end{equation}
Let us estimate $E_k$ and $\frac{P_{k+1}}{y^2}$ separately. On $E_k,$ by \textnormal{(\romannumeral5)} of Lemma \ref{property of degree}, (\ref{assump on b_k}), (\ref{under induction hypo, esti of S_k}) and \textnormal{(\romannumeral1)} of Lemma \ref{property of degree}, we have for the components of $E_k,$
\begin{align*}
	b_1b_k(\Lambda T_k-(2k-\gamma)T_k)&\sim (k,k-1,k+1),\\
	b_1\Lambda S_k&\sim (k,k-1,k+1),\\
	\left[(2j-\gamma)b_1-\frac{b_{j+1}}{b_j}\right]\left(b_j\frac{\partial S_k}{\partial b_j}\right)&\sim (k,k-1,k+1).
\end{align*}
Hence $E_k\sim (k,k-1,k+1).$ On $\frac{P_{k+1}}{y^2},$ note that it is the finite linear combinations of terms of the form
\begin{equation*}
	M_J:=\frac{f^{(j)}(Q)}{y^2}\prod_{m=1}^{L}b_m^{i_m}T_m^{i_m}\prod_{m=2}^{L+2}S_m^{j_m},
\end{equation*}
where $J=(i_1,\cdots,i_L,j_2,\cdots,j_{L+2}),$ $|J|_1=j,$ $|J|_2=k+1,$ $2\le j\le \min \{k+1,L+2\}.$
Note that by (\ref{asymp: ground state}), 
\begin{align*}
	&\text{when}\,\,\,y\rightarrow 0,\,\,\,f^{(j)}(Q)\lesssim 1.\\
	&\text{when}\,\,\,y\rightarrow \infty,\,\,\,f^{(j)}(Q)\lesssim \left\{\begin{aligned}
		&y^{-\gamma}\,\,\,\text{for}\,\,\,j\,\,\,\text{even}\\ &1\,\,\,\text{for}\,\,\,j\,\,\,\text{odd}
	\end{aligned}\right. .
\end{align*}
Combined with (\ref{assump on b_k}), \textnormal{(\romannumeral4)} of Lemma \ref{property of degree} and (\ref{under induction hypo, esti of S_k}), we get
\begin{align*}
	&\text{when}\,\,\,y\rightarrow 0,\,\,\,M_J\lesssim b_1^{k+1}y^{\sum_{m=1}^{L}(2m+2)i_m+\sum_{m=2}^{L+2}(2m+2)j_m-2}\ll b_1^{k+1}y^{2k+2}.\\
	&\text{when}\,\,\,y\rightarrow \infty,\,\,\,M_J\lesssim\left\{\begin{aligned}
		&b_1^{k+1}y^{2k-\gamma(j+1)-2\sum_{m=2}^{L+2}j_m}\,\,\,\text{for}\,\,\,\text{even}\,\,\,j\\ &b_1^{k+1}y^{2k-\gamma j-2\sum_{m=2}^{L+2}j_m}\,\,\,\text{for}\,\,\,\text{odd}\,\,\,j
	\end{aligned}\right. \ll b_1^{k+1}y^{2(k-1)-\gamma}.
\end{align*}
Hence $M_J\sim (k,k-1,k+1)$ and same holds for $F_{k+1}.$ By the definition of $F_{k+1}$ and induction hypothesis, $\frac{\partial F_k}{\partial b_m}=0$ for $2\le k\le m\le L+2$ is easily verified. This completes the proof of (\ref{degree of F_k}). Then again by \textnormal{(\romannumeral3)} of Lemma \ref{property of degree}, (\ref{degree of S_k}) and (\ref{dependence for S_k with respect to b_m}) hold true.\par 
Here we omit the details for the proof of (\ref{esti of Psib in the scale B1}) and (\ref{esti of Psib in the scale M}) since one needs only to use the degrees of the components of $\Psi_b$ which are already known. This concludes the whole proof. 
\end{pf}
\vspace{\baselineskip}

\subsection{Localized approximation}
\begin{prop}\label{localized approximation}
Under the assumptions in Proposition \ref{first approximation}, assuming in addition $|(b_1)_s|\lesssim b_1^2.$ Defining the localized approximation of (\ref{eq: renormalized flow}) as 
\begin{equation}\label{form of approximation}
	\widetilde{Q}_{b(s)}(y):=Q(y)+\sum\limits_{k=1}^L b_k \widetilde{T}_k+\sum\limits_{k=2}^{L+2}\widetilde{S}_k\,\,\,\text{with}\,\,\,\widetilde{T}_k:=\chi_{B_1}T_k,\,\,\,\widetilde{S}_k:=\chi_{B_1}S_k.
\end{equation}
Then $\widetilde{Q}_b$ satisfies the equation
\begin{equation}\label{eq: eq of approximation}
	\partial_s \widetilde{Q}_b-\partial_y^2 \widetilde{Q}_b-\frac{(d-3)}{y}\partial_y \widetilde{Q}_b+b_1 \Lambda \widetilde{Q}_b+\frac{(d-2)}{y^2}f(\widetilde{Q}_b)=\widetilde{\Psi}_b+\chi_{B_1}Mod(t)
\end{equation}
with $\widetilde{\Psi}_b$ satisfying the following properties.\\
\textnormal{(\romannumeral1)} For all $0\le m\le L-1,$
\begin{equation}\label{1 esti of tilde Psib}
	\int |\mathscr{L}^{\hbar+m+1}\widetilde{\Psi}_b|^2+\int \frac{|\mathscr{A}\mathscr{L}^{\hbar+m}\widetilde{\Psi}_b|^2}{1+y^2}+\int \frac{|\mathscr{L}^{\hbar+m}\widetilde{\Psi}_b|^2}{1+y^4}+\int \frac{|\widetilde{\Psi}_b|^2}{1+y^{4(\hbar+m+1)}}\lesssim b_1^{2m+2+2(1-\delta)-C_L \eta},
\end{equation}
and
\begin{equation}\label{2 esti of tilde Psib}
	\int |\mathscr{L}^{\hbar+L+1}\widetilde{\Psi}_b|^2+\int \frac{|\mathscr{A}\mathscr{L}^{\hbar+L}\widetilde{\Psi}_b|^2}{1+y^2}+\int \frac{|\mathscr{L}^{\hbar+L}\widetilde{\Psi}_b|^2}{1+y^4}+\int \frac{|\widetilde{\Psi}_b|^2}{1+y^{4(\hbar+L+1)}}\lesssim b_1^{2L+2+2(1-\delta)(1+\eta)}.
\end{equation}
\textnormal{(\romannumeral2)} For all $M\le \frac{B_1}{2}$ and $0\le m\le L,$
\begin{equation}\label{3 esti of tilde Psib}
	\int_{y\le 2M} |\mathscr{L}^{\hbar+m+1}\widetilde{\Psi}_b|^2\lesssim M^{C}b_1^{2L+6}.
\end{equation}
\textnormal{(\romannumeral3)} For all $0\le m\le L,$
\begin{equation}\label{4 esti of tilde Psib}
	\int_{y\le 2 B_0}|\mathscr{L}^{\hbar+m+1}\widetilde{\Psi}_b|^2+\int_{y\le 2 B_0} \frac{|\widetilde{\Psi}_b|^2}{1+y^{4(\hbar+m+1)}}\lesssim b_1^{2m+4+2(1-\delta)-C_L \eta}.
\end{equation}
\end{prop}
\begin{pf}
	\normalfont
	Direct computation gives 
	\begin{align*}
		&\partial_s \widetilde{Q}_b-\partial_y^2 \widetilde{Q}_b-\frac{(d-3)}{y}\partial_y \widetilde{Q}_b+b_1 \Lambda \widetilde{Q}_b+\frac{(d-2)}{y^2}f(\widetilde{Q}_b)\\=&\chi_{B_1}\left[\partial_s {Q}_b-\partial_y^2 {Q}_b-\frac{(d-3)}{y}\partial_y {Q}_b+b_1 \Lambda {Q}_b+\frac{(d-2)}{y^2}f({Q}_b)\right]\\&+\Theta_{b}\left[\partial_s \chi_{B_1}-\left(\partial_y^2 \chi_{B_1}+\frac{(d-3)}{y}\partial_y \chi_{B_1}\right)+b_1\Lambda \chi_{B_1}\right]-2\partial_y\chi_{B_1}\partial_y \Theta_{b}+b_1(1-\chi_{B_1})\Lambda Q\\ &+\frac{(d-2)}{y^2}\left[f(\widetilde{Q}_b)-f(Q)-\chi_{B_1}\left(f(Q_b)-f(Q)\right)\right].
	\end{align*}
Then by (\ref{eq: eq of first approximation}),
\begin{equation*}
	\partial_s \widetilde{Q}_b-\partial_y^2 \widetilde{Q}_b-\frac{(d-3)}{y}\partial_y \widetilde{Q}_b+b_1 \Lambda \widetilde{Q}_b+\frac{(d-2)}{y^2}f(\widetilde{Q}_b)\\=:\chi_{B_1}Mod(t)+\widetilde{\Psi}_b,
\end{equation*}
where 
\begin{align*}
	\widetilde{\Psi}_b&:=\chi_{B_1}\Psi_b+\widetilde{\Psi}_b^{(1)}+\widetilde{\Psi}_b^{(2)}+\widetilde{\Psi}_b^{(3)},\\ \widetilde{\Psi}_b^{(1)}&:=b_1(1-\chi_{B_1})\Lambda Q,\\ \widetilde{\Psi}_b^{(2)}&:=\frac{(d-2)}{y^2}\left[f(\widetilde{Q}_b)-f(Q)-\chi_{B_1}\left(f(Q_b)-f(Q)\right)\right],\\ \widetilde{\Psi}_b^{(3)}&:=\Theta_{b}\left[\partial_s \chi_{B_1}-\left(\partial_y^2 \chi_{B_1}+\frac{(d-3)}{y}\partial_y \chi_{B_1}\right)+b_1\Lambda \chi_{B_1}\right]-2\partial_y\chi_{B_1}\partial_y \Theta_{b}.
\end{align*}
We only estimate the contribution of $\widetilde{\Psi}_b^{(2)}$ in (\ref{1 esti of tilde Psib})-(\ref{4 esti of tilde Psib}). Note that in $\widetilde{\Psi}_b^{(2)},$ $y$ is supported in $B_1\le y\le 2B_1,$ thus its contribution to (\ref{3 esti of tilde Psib}) and (\ref{4 esti of tilde Psib}) hold trivially. Let us now estimate its contribution to (\ref{1 esti of tilde Psib}) and (\ref{2 esti of tilde Psib}). By Taylor expansion,
\begin{equation}\label{taylor exp of fQb-fQ}
	f(Q_b)-f(Q)=f'(Q)\Theta_{b}+\frac{f''(Q)}{2}\Theta_{b}^2+\Theta_{b}^3,\,\,\,\text{with}\,\,\,B_1\le y\le 2B_1.
\end{equation}
Note that for $2\le k\le L,$ by (\ref{assump on b_k}) and (\ref{degree of S_k}) we see that $|S_k|\lesssim b_1^ky^{2(k-1)-\gamma}.$ In comparison, $|b_kT_k|\lesssim b_1^ky^{2k-\gamma},$ which follows from (\ref{assump on b_k}) and \textnormal{(\romannumeral4)} of Lemma \ref{property of degree}. Similarly, $|S_{L+1}|\lesssim b_1^{L+1}y^{2L-\gamma}$ and $|S_{L+2}|\lesssim b_1^{L+2}y^{2(L+1)-\gamma},$ in comparison, $|b_LT_L|\lesssim b_1^Ly^{2L-\gamma}.$ Therefore, the main order term of $\Theta_{b}$ is $\sum\limits_{k=1}^L b_kT_k,$ or say  
\begin{equation}\label{esti of Thetab for y sim B1}
	|\Theta_b|\lesssim \sum_{k=1}^{L}b_1^ky^{2k-\gamma}1_{B_1\le y\le 2B_1}.
\end{equation} 
In particular, since $b_1^ky^{2k-\gamma}1_{B_1\le y\le 2B_1}\lesssim b_1^{\frac{\gamma}{2}+\eta(\frac{\gamma}{2}-k)},$ we have $|\Theta_{b}|\ll 1.$
Substituting (\ref{esti of Thetab for y sim B1}) into (\ref{taylor exp of fQb-fQ}) and making use of (\ref{asymp: ground state}), we get
\begin{equation*}
	\left\{\begin{aligned}
		|f(Q_b)-f(Q)|&\lesssim |\Theta_{b}|\lesssim \sum_{k=1}^{L}b_1^ky^{2k-\gamma}1_{B_1\le y\le 2B_1}\\ |f(\widetilde{Q}_b)-f(Q)|&\lesssim \chi_{B_1}|\Theta_{b}|\lesssim \sum_{k=1}^{L}b_1^ky^{2k-\gamma}1_{B_1\le y\le 2B_1}
	\end{aligned}\right. \Longrightarrow |\widetilde{\Psi}_b^{(2)}|\lesssim \sum_{k=1}^{L}b_1^ky^{2(k-1)-\gamma}1_{B_1\le y\le 2B_1}.
\end{equation*}
We further estimate that $|\widetilde{\Psi}_b^{(2)}|\lesssim \sum\limits_{k=1}^L b_1^kB_1^{2(k-1)}y^{-\gamma}1_{B_1\le y\le 2B_1}=b_1y^{\gamma}\sum\limits_{k=1}^L b_1^{-\eta (k-1)}1_{B_1\le y\le 2B_1},$ then 
\begin{align*}
	\int |\mathscr{L}^{\hbar+m+1}\widetilde{\Psi}_b^{(2)}|^2&\lesssim b_1^2\sum\limits_{k=1}^L b_1^{-2(k-1)\eta}\int_{B_1\le y\le 2B_1} |y^{-\gamma-2(\hbar+m+1)}|^2y^{d-3}\,\mathrm{d}y\\ &\lesssim b_1^{2m+2+2(1+\eta)(1-\delta)}\sum\limits_{k=1}^L b_1^{(2m-2k+2)\eta},\,\,\,\text{for all}\,\,\,0\le m\le L.
\end{align*}
This concludes the proof.
\end{pf}

\vspace{\baselineskip}
\section{Linearization of $\{b_k\}_{k=1}^L$}\label{Linearization bk for k from 1 to l}
Denoting $\{b_k^e\}_{k=1}^L$ as the solution of 
\begin{equation}\label{eq: eq for bke}
	\left\{\begin{aligned}
		&(b_k^e)_s+(2k-\gamma)b_1^e b_k^e-b_{k+1}^e=0\\
		&b_{l+1}^e=b_{l+2}^e=\cdots=b_L^e=0
	\end{aligned}\right. ,
\end{equation}
where $b_{L+1}^e:=0$ and $L$ satisfying  
\begin{equation}\label{def of l}
	\frac{\gamma}{2}<l\ll L
\end{equation} is an integer to be chosen later. One can find a set of solution explicitly in the form $b_k^e=c_k s^{-k},$ with 
\begin{equation}\label{coefficient of bke}
	\left\{\begin{aligned}
		c_1&=\frac{l}{2l-\gamma},\\
		c_{k+1}&=-\frac{\gamma(l-k)}{2l-\gamma}c_k,\,\,\,1\le k\le l-1,\\
		c_{l+1}&=c_{l+2}=\cdots=c_L=0.
	\end{aligned}\right.
\end{equation}
We abuse notation and still denote this set of solution as $\{b_k^e\}_{k=1}^L$ since we will not use any other solutions of (\ref{eq: eq for bke}). Then one can calculate pertubations near $\{b_k^e\}_{k=1}^L$ as the following Lemma, the details of the proof are the same as in Lemma 3.7 in \cite{merle2015type}. 
\begin{lem}\label{linearization of bk from 1 to l}
	Let $b_k(s)=b_k^e(s)+\frac{\mathcal{U}_k(s)}{s^k},$ $1\le k\le l,$ denoting $\mathcal{U}:=(\mathcal{U}_1,\cdots,\mathcal{U}_l).$ Then for $1\le k\le l-1,$
\begin{align}
		(b_k)_s+(2k-\gamma)b_1 b_k-b_{k+1}&=\frac{1}{s^{k+1}}[s(\mathcal{U}_k)_s-(A_l\mathcal{U})_k+O(|\mathcal{U}|^2)],\label{linearization for bk form 1 to l-1}\\
		(b_l)_s-(2l-\gamma)b_1 b_l&=\frac{1}{s^{l+1}}[s(\mathcal{U}_l)_s-(A_l\mathcal{U})_l+O(|\mathcal{U}|^2)],\label{linearization for bl}
\end{align}
where $A_l=(a_{i,j})_{l\times l}$ with 
\begin{equation*}
	\left\{\begin{aligned}
		a_{1,1}&=\frac{\gamma(l-1)}{2l-\gamma}-(2-\gamma)c_1,\\
		a_{i,i}&=\frac{\gamma(l-i)}{2l-\gamma},\,\,\,2\le i\le l,\\
		a_{i,i+1}&=1,\,\,\,1\le i\le l-1,\\
		a_{i,1}&=-(2i-\gamma)c_i,\,\,\,2\le i\le l,\\
		a_{i,j}&=0,\,\,\,\text{otherwise}.
	\end{aligned}\right.
\end{equation*}
Moreover, $A_l$ is diagonalizable with 
\begin{equation}\label{dig of Al}
	A_l=P_l^{-1}D_l P_l,\,\,\,D_l=diag\left\{-1,\frac{2\gamma}{2l-\gamma},\frac{3\gamma}{2l-\gamma},\cdots,\frac{l\gamma}{2l-\gamma}\right\}.
\end{equation}
\end{lem}  
\begin{rk}
	\normalfont
	Standard linearization of the $b$-system near $\{b_k^e\}_{k=1}^L$ yields a Jacobian similar to the structure of $A_l.$ However, by further shaping the remainder as $\frac{\mathcal{U}_k(s)}{s^k},$ one finds advantages in calculation. 
\end{rk}

\vspace{\baselineskip}
\section{Decomposition of the solution and coercivity-determined estimates on the remainder term}\label{Decomposition of the solution and coercivity-determined estimates on the remainder term}
\subsection{Decomposition of the solution}
In this section, we shall use implicit function theorem to show the existence of decomposition 
\begin{equation}\label{decomposition}
	u(t,r)=(\widetilde{Q}_b+q)\Big(t,\frac{r}{\lambda(t)}\Big)
\end{equation}
satisfying orthogonality conditions 
\begin{equation}\label{orthogonal condition}
	\langle q(s,y),\mathscr{L}^i \Phi_M\rangle=0,\,\,\,\text{for}\,\,\,0\le i\le L.
\end{equation}\par 
Indeed, let $u_0$ be close to $Q$ in some sense, then this closeness is propagated on a small time interval $[0,t_1).$ Defining the map
\begin{equation*}
	\mathcal{T}\,\colon\,(t,\lambda,b_1,\cdots,b_L)\longmapsto \Big(\langle u(t)-(\widetilde{Q}_b)_{\lambda},(\mathscr{L}^i \Phi_M)_{\lambda}\rangle\Big)_{i=0,\cdots,L}.
\end{equation*}
Then we choose $u_0$ such that $\mathcal{T}$ maps $(0,\lambda^*,b_1^*,\cdots,b_L^*)$ to the zero vector, for some $\lambda^*$ close to $1$ and $b_i^*$ close to $0$ for all $1\le i\le L.$ Note that by direct computation, the Jacobian of $\mathcal{T}$ at $t=0,$ $\lambda=1,$ $b=0$ is 
\begin{equation*}
	(-1)^{\frac{(1+L)L}{2}}\langle \chi_M \Lambda Q,\Lambda Q\rangle^{L+1}+\text{small correction},
\end{equation*}
which is nonzero. Then by implicit function theorem, there exists unique functions $\lambda=\lambda(t),$ $b=b(t)$ such that $\left(\big\langle u(t)-(\widetilde{Q}_{b(t)})_{\lambda(t)},(\mathscr{L}^i \Phi_M)_{\lambda(t)}\big\rangle \right)_{i=0,\cdots,L}\equiv 0$ on some time interval $[0,t_1^*).$

\vspace{\baselineskip}
\subsection{Coercivity-determined estimates on $q$}
Substituting (\ref{decomposition}) into (\ref{eq: renormalized flow}) and making use of (\ref{orthogonal condition}), we get the equation for the remainder term $q$ as
\begin{equation}\label{eq: eq of q}
	\partial_s q-\frac{\lambda_s}{\lambda}\Lambda q+\mathscr{L}q=-\widetilde{\Psi}_b-\widehat{Mod}+\mathcal{H}(q)-\mathcal{N}(q)=:\mathcal{F},
\end{equation}
where
\begin{align}
	\widehat{Mod}:&=-\Big(\frac{\lambda_s}{\lambda}+b_1\Big)\Lambda \widetilde{Q}_b+\chi_{B_1}Mod(t),\label{def of widehat Mod}\\
	\mathcal{H}(q):&=\frac{(d-2)}{y^2}[f'(Q)-f'(\widetilde{Q}_b)]q,\label{def of Hq}\\
	\mathcal{N}(q):&=\frac{(d-2)}{y^2}[f(\widetilde{Q}_b+q)-f(\widetilde{Q}_b)-f'(\widetilde{Q}_b)q].\label{def of Nq}
\end{align}\par
In original variable, denoting $v(t,r):=q(s,y),$ then (\ref{eq: eq of q}) becomes
\begin{equation}\label{eq: eq of v}
	\partial_t v+\mathscr{L}_{\lambda}v=\frac{1}{\lambda^2}\mathcal{F}_{\lambda}.
\end{equation}
Before we dig more into (\ref{eq: eq of q}) or (\ref{eq: eq of v}) for dynamic-determined estimates, which we shall call modulation estimates and energy estimates later. We make some preparations about estimates on $q$ which are coercivity-determined. For simplicity, denoting $\mathscr{E}_{2i}:=\mathscr{E}_{2i}(q)=\int |\mathscr{L}^i q|^2,$ for $i\in \mathbb{N}.$

\begin{lem}\label{coercivity-determined esti on q}
Recalling that $\Bbbk:=L+\hbar+1,$ $L\gg l$ is a large integer and $\hbar$ is defined as in (\ref{def: def of h}), then we have the following coercivity-determined estimates.\\
\textnormal{(\romannumeral1)} Near the origin $q$ has the Taylor expansion in the form
\begin{equation}\label{taylor expansion of q near zero}
	q=\sum\limits_{i=1}^{\Bbbk} c_i T_{\Bbbk-i}+r_q, 
\end{equation}
with bounds 
\begin{align}
	|c_i|&\lesssim \sqrt{\mathscr{E}_{2\Bbbk}},\label{esti for ci}\\
	|\partial_y^i r_q|&\lesssim y^{2\Bbbk+1-\frac{d}{2}-j}|\ln y|^{\Bbbk}\sqrt{\mathscr{E}_{2\Bbbk}},\,\,\,0\le j\le 2\Bbbk-1,\,\,\,y<1.\label{esti for rq}
\end{align}
\textnormal{(\romannumeral2)} Pointwise bound near the origin: for $y<1,$
\begin{align}
	|q_{2i}|+|\partial_y^{2i} q|&\lesssim y^{-\frac{d}{2}+3}|\ln y|^{\Bbbk-i}\sqrt{\mathscr{E}_{2\Bbbk}}\,\,\,\text{for}\,\,\,0\le i\le \Bbbk-1,\label{esti for even deri of q near 0}\\
	|q_{2i-1}|+|\partial_y^{2i-1} q|&\lesssim y^{-\frac{d}{2}+2}|\ln y|^{\Bbbk-i}\sqrt{\mathscr{E}_{2\Bbbk}}\,\,\,\text{for}\,\,\,1\le i\le \Bbbk.\label{esti for odd deri of q near 0}
\end{align}
\textnormal{(\romannumeral3)} Weighted bounds: for $1\le m\le \Bbbk,$
\begin{equation}\label{weighted bd sum version}
	\sum\limits_{i=0}^{2m} \int \frac{|\partial_y^i q|^2}{1+y^{4m-2i}}\lesssim \mathscr{E}_{2m}.
\end{equation}
Moreover, let $(i,j)\in \mathbb{N}\times \mathbb{N}_{+}$ with $2\le i+j\le 2\Bbbk,$ then
\begin{equation}\label{weighted bd}
	\int \frac{|\partial_y^i q|^2}{1+y^{2j}}\lesssim \left\{\begin{aligned}
		&\mathscr{E}_{2m},\,\,\,\text{for}\,\,\, i+j=2m,\,\,\, 1\le m\le \Bbbk.\\
		&\sqrt{\mathscr{E}_{2m}}\sqrt{\mathscr{E}_{2(m+1)}},\,\,\,\text{for}\,\,\, i+j=2m+1,\,\,\,1\le m\le \Bbbk-1.
	\end{aligned}\right.
\end{equation}
\textnormal{(\romannumeral4)} Pointwise bound far away: let $(i,j)\in \mathbb{N}\times\mathbb{N}$ with $1\le i+j\le 2\Bbbk-1,$ then for $y\ge 1,$
\begin{equation}\label{pointwise bd away from origin}
	\left|\frac{\partial_y^i q}{y^j}\right|^2\lesssim \frac{1}{y^{d-4}}\left\{\begin{aligned}
		&\mathscr{E}_{2m},\,\,\,\text{for}\,\,\, i+j+1=2m,\,\,\,1\le m\le \Bbbk.\\
		&\sqrt{\mathscr{E}_{2m}}\sqrt{\mathscr{E}_{2(m+1)}},\,\,\,\text{for}\,\,\, i+j=2m,\,\,\,1\le m\le \Bbbk-1.
	\end{aligned}\right.
\end{equation}
\end{lem} 
\begin{pf}
	\normalfont
	This can be proved in a similar way as in Appendix B in \cite{ghoul2018stability}, we shall still provide details since different asymptotics are involved.\par
	Proof of \textnormal{(\romannumeral1)}. We claim: for $1\le m\le \Bbbk,$ $q_{2\Bbbk-2m}$ admits the following expansion and estimates near the origin:
	\begin{equation}\label{expansion of q2k-2m at origin}
		\left\{\begin{aligned}
			q_{2\Bbbk-2m}&=\sum\limits_{i=1}^m c_{i,m}T_{m-i}+r_{2m},\\
			|c_{i,m}|&\lesssim \sqrt{\mathscr{E}_{2\Bbbk}},\\
			|\partial_y^i r_{2m}|&\lesssim y^{2m+1-\frac{d}{2}-j}|\ln y|^m \sqrt{\mathscr{E}_{2\Bbbk}},\,\,\,\text{for}\,\,\,0\le j\le 2m-1.
		\end{aligned}\right.
	\end{equation}
	Then \textnormal{(\romannumeral1)} immediately follows when one chooses $m=\Bbbk,$ and we focus on the proof of the claim. Defining the sequence $\{r_i\}_{i=0}^{2\Bbbk}$ recursively as $r_0:=q_{2\Bbbk},$ and
	\begin{equation*}
		\left\{\begin{aligned}
			r_{2i+1}(y)&:=\frac{1}{y^{d-3}\Lambda Q}\int_{0}^{y} r_{2i}\Lambda Q\cdot x^{d-3}\,\mathrm{d}x\\
			r_{2i+2}(y)&:=-\Lambda Q\int_{a}^{y} \frac{r_{2i+1}}{\Lambda Q}\,\mathrm{d}x
		\end{aligned}\right.,\,\,\,\text{for}\,\,\,0\le i\le \Bbbk-1,
	\end{equation*} 
	 where we choose $a\in \left(\frac{1}{2},1\right)$ such that 
	 \begin{equation*}
	 	|q_{2\Bbbk-i}(a)|^2\lesssim \int_{y\le 1} |q_{2\Bbbk-i}|^2\lesssim \mathscr{E}_{2\Bbbk},\,\,\,\text{for}\,\,\,1\le i \le 2\Bbbk.
	 \end{equation*}
	 Note that Lemma \ref{coercivity of iterate of L} justifies the last inequality, and the definition fulfills the relation $\mathscr{A^*}r_{2i+1}=r_{2i},$ $\mathscr{A}r_{2i+2}=r_{2i+1},$ for $0\le i\le \Bbbk-1.$ Next we prove (\ref{expansion of q2k-2m at origin}) inductively. When $m=1,$ note that by $\mathscr{A^*}q_{2\Bbbk-1}=q_{2\Bbbk},$ we have
	 \begin{equation*}
	 	q_{2\Bbbk-1}=\frac{1}{y^{d-3}\Lambda Q}\int_{0}^{y} q_{2\Bbbk}\Lambda Q\cdot x^{d-3}\,\mathrm{d}x+\frac{C}{y^{d-3}\Lambda Q}.
	 \end{equation*}
	 From the facts that $\int \frac{|q_{2\Bbbk-1|^2}}{y^2}\lesssim \mathscr{E}_{2\Bbbk},$ $\Lambda Q(y)\simeq y^2$ as $y\rightarrow 0,$ and $\int_{y<1} \frac{\left|\frac{1}{y^{d-3}\Lambda Q}\right|^2}{y^2}=+\infty,$ one deduces $C=0.$ That is, $r_1(y)=q_{2\Bbbk-1}(y).$ And by H\"older, we get 
	 \begin{align*}
	 	|r_1(y)|&\lesssim \frac{1}{y^{d-3}}\left(\int_{0}^y |q_{2\Bbbk}|^2x^{d-3}\,\mathrm{d}x\right)^{\frac{1}{2}}\left(\int_{0}^y x^4\cdot x^{d-3}\,\mathrm{d}x\right)^{\frac{1}{2}}\\
	 	&\lesssim y^{-\frac{d}{2}+2}\sqrt{\mathscr{E}_{2\Bbbk}},\,\,\,\text{for}\,\,\,y<1.
	 \end{align*}
	 By the definition of $\{r_i\},$ $r_2(y)=-\Lambda Q\int_{a}^y \frac{r_1}{\Lambda Q}\,\mathrm{d}x,$ then straight calculation yields 
	 \begin{equation*}
	 	|r_2(y)|\lesssim y^2\int_{a}^y \frac{x^{-\frac{d}{2}+2}\sqrt{\mathscr{E}_{2\Bbbk}}}{x^2}\,\mathrm{d}x
	 	\lesssim y^{-\frac{d}{2}+3}|\ln y|\mathscr{E}_{2\Bbbk},\,\,\,\text{for}\,\,\,y<1.
	 \end{equation*} 
	 Note that $\mathscr{A}r_2=r_1=q_{2\Bbbk-1},$ then $\mathscr{L}r_2=\mathscr{A^*}q_{2\Bbbk-1}=q_{2\Bbbk}=\mathscr{L}q_{2\Bbbk-2},$ and it tells that 
	 \begin{equation}\label{expansion of q2k-2 at origin}
	 	q_{2\Bbbk-2}=c\Lambda Q+r_2,
	 \end{equation}  
	 where we also used the facts that $\int \frac{|q_{2\Bbbk-2}|^2}{y^4}\lesssim \mathscr{E}_{2\Bbbk}$ and $\int_{y<1} \frac{|\Gamma|^2}{y^4}=+\infty.$
	 Then substituting the point $a$ into (\ref{expansion of q2k-2 at origin}), we derive $|c|\lesssim \sqrt{\mathscr{E}_{2\Bbbk}},$ then it follows that
	 \begin{equation*}
	 	|q_{2\Bbbk-2}|\lesssim y^{-\frac{d}{2}+3}|\ln y|\mathscr{E}_{2\Bbbk},\,\,\,\text{for}\,\,\,y<1.
	 \end{equation*}
	 Next, applying the definition of $\mathscr{A}$ and asymptotics of $V,$ we have 
	 \begin{equation*}
	 	|\partial_y r_2|\lesssim |r_1|+\left|\frac{r_2}{y}\right|\lesssim y^{-\frac{d}{2}+2}|\ln y|\mathscr{E}_{2\Bbbk},\,\,\,\text{for}\,\,\,y<1.
	 \end{equation*}  
	 This concludes the proof when $m=1.$ Now assuming the cases for ``$\le m$'' are true, and we prove that the case for $m+1$ holds true. Note that one can readily verify the inequality
	 \begin{equation*}
	 	\int_{0}^y |\ln x|^n\,\mathrm{d}x\lesssim y|\ln y|^n,
	 \end{equation*}
	 for all positive integers $n.$ Then it combined with induction hypothesis gives 
	 \begin{align*}
	 	|r_{2m+1}|&=\left|\frac{1}{y^{d-3}}\int_{0}^y r_{2m}\Lambda Q\cdot x^{d-3}\,\mathrm{d}x\right|\\
	 	&\lesssim \frac{1}{y^{d-1}}\sqrt{\mathscr{E}_{2\Bbbk}}\int_{0}^y x^{2m+\frac{d}{2}}|\ln x|^m\,\mathrm{d}x\\
	 	&\lesssim y^{2m+2-\frac{d}{2}}|\ln y|^m\sqrt{\mathscr{E}_{2\Bbbk}}.
	 \end{align*} 
	 Hence,
	 \begin{equation*}
	 	|r_{2m+2}|=\left|-\Lambda Q\int_{a}^y \frac{r_{2m+1}}{\Lambda Q}\,\mathrm{d}x\right|\lesssim y^2 \int_y^a \frac{x^{2m+2-\frac{d}{2}}\sqrt{\mathscr{E}_{2\Bbbk}}|\ln x|^m}{x^2}\,\mathrm{d}x,
	 \end{equation*}
	 and we further estimate it into two cases. If $2m+1-\frac{d}{2}<0,$ then 
	 \begin{equation*}
	 	|r_{2m+2}|\lesssim y^2\cdot y^{2m+1-\frac{d}{2}}\sqrt{\mathscr{E}_{2\Bbbk}}\int_y^a \frac{|\ln x|^m}{x}\,\mathrm{d}x\lesssim y^{2m+3-\frac{d}{2}}|\ln y|^{m+1}\sqrt{\mathscr{E}_{2\Bbbk}}.
	 \end{equation*}
	 If $2m+1-\frac{d}{2}\ge 0,$ note that using induction one can readily verify that 
	 \begin{equation*}
	 	\int_y^a x^{2m+1-\frac{d}{2}}\cdot \frac{|\ln x|^m}{x}\,\mathrm{d}x\lesssim y^{2m+1-\frac{d}{2}}|\ln y|^m,
	 \end{equation*}
	 then it immediately follows 
	 \begin{equation*}
	 	|r_{2m+2}|\lesssim y^{2m+3-\frac{d}{2}}|\ln y|^{m}\sqrt{\mathscr{E}_{2\Bbbk}}.
	 \end{equation*}
	 Combining above two cases together, we get
	 \begin{equation*}
	 	|r_{2m+2}|\lesssim y^{2m+3-\frac{d}{2}}|\ln y|^{m+1}\sqrt{\mathscr{E}_{2\Bbbk}}.
	 \end{equation*}
	 Since $\mathscr{A}r_{2m+2}=r_{2m+1},$ $\mathscr{L}r_{2m+2}=r_{2m},$ one sees that 
	 \begin{equation*}
	 	\mathscr{L}q_{2\Bbbk-2(m+1)}=q_{2\Bbbk-2m}=\sum\limits_{i=1}^m c_{i,m}T_{m-i}+r_{2m}=\sum\limits_{i=1}^m -c_{i,m}\mathscr{L}T_{m+1-i}+\mathscr{L}r_{2m+2}.
	 \end{equation*}
     Due to the facts that $\int_{y\le 1}\frac{|q_{2\Bbbk-2(m+1)}|^2}{y^4}\lesssim \mathscr{E}_{2\Bbbk}$ and $\int_{y<1} \frac{|\Gamma|^2}{y^4}=+\infty,$ we deduce 
     \begin{equation*}
     	q_{2\Bbbk-2(m+1)}=\sum\limits_{i=1}^m -c_{i,m}T_{m+1-i}+r_{2m+2}+\tilde{c}\Lambda Q.
     \end{equation*}
     Then substituting the point $a$ into above expansion gives $|\tilde{c}|\lesssim \sqrt{\mathscr{E}_{2\Bbbk}}.$ And by induction hypothesis and direct computations, we have
     \begin{equation*}
     	|\partial_y^j r_{2m+2}|\lesssim \sum\limits_{i=0}^j \frac{|r_{2m+2-2i}|}{y^{j-i}}\lesssim \frac{\sqrt{\mathscr{E}_{2\Bbbk}}\sum\limits_{i=0}^j y^{2m+3-i-\frac{d}{2}}|\ln y|^{m+1}}{y^{j-i}}\lesssim y^{2m+3-\frac{d}{2}-j}|\ln y|^{m+1}\sqrt{\mathscr{E}_{2\Bbbk}}.
     \end{equation*}
      This concludes the proof of (\ref{expansion of q2k-2m at origin}).\par 
      Proof of \textnormal{(\romannumeral2)}. Using (\ref{expansion of q2k-2m at origin}) for $m=\Bbbk-i,$ we get 
      \begin{equation*}
      	|q_{2i}|\lesssim \sum\limits_{i=1}^m \sqrt{\mathscr{E}_{2\Bbbk}}y^{2+2(m-i)}+\sqrt{\mathscr{E}_{2\Bbbk}}y^{2m+1-\frac{d}{2}}|\ln y|^m\lesssim y^{-\frac{d}{2}+3}|\ln y|^m \sqrt{\mathscr{E}_{2\Bbbk}},
      \end{equation*}
      for $0\le i\le \Bbbk-1.$ And the second inequality follows from the relation $\mathscr{A}q_{2i-2}=q_{2i-1},$ for $1\le i\le \Bbbk.$\par 
      Proof of \textnormal{(\romannumeral3)}. Note that $|\partial_y^i q|\lesssim \sum\limits_{j=0}^i \frac{|q_j|}{y^{i-j}},$ then by Lemma \ref{coercivity of iterate of L} and \textnormal{(\romannumeral2)}, we have
      \begin{align*}
      	\sum\limits_{i=0}^{2m} \int \frac{|\partial_y^i q|^2}{1+y^{4m-2i}}&\lesssim \mathscr{E}_{2m}+\sum\limits_{i=0}^{2m-1}\int_{y<1} |\partial_y^i q|^2+\sum\limits_{i=0}^{2m-1}\int_{y>1}\frac{|\partial_y^i q|^2}{y^{4m-2i}}\\
      	&\lesssim \mathscr{E}_{2m}+\mathscr{E}_{2\Bbbk}\int_{y<1} y^3 |\ln y|^{2\Bbbk}\,\mathrm{d}y+\sum\limits_{i=0}^{2m-1}\sum\limits_{j=0}^i \int_{y>1} \frac{|q_j|^2}{y^{4m-2j}}\lesssim \mathscr{E}_{2m}.
      \end{align*} 
      And consequently, one derives that if $i+j=2m$ with $1\le m\le \Bbbk,$ 
      \begin{equation*}
      	\int \frac{|\partial_y^i q|^2}{1+y^{2j}}=\int \frac{|\partial_y^i q|^2}{1+y^{4m-2i}}\lesssim \mathscr{E}_{2m}.
      \end{equation*}
      If $i+j=2m+1$ with $1\le m\le \Bbbk-1,$
      \begin{equation*}
      	\int \frac{|\partial_y^i q|^2}{1+y^{2j}}=\int \frac{|\partial_y^i q|^2}{1+y^{4m+2-2i}}\lesssim \left(\int \frac{|\partial_y^i q|^2}{1+y^{4m-2i}}\right)^{\frac{1}{2}}\left(\int \frac{|\partial_y^i q|^2}{1+y^{4m+4-2i}}\right)^{\frac{1}{2}}\lesssim \sqrt{\mathscr{E}_{2m}}\sqrt{\mathscr{E}_{2(m+1)}}.
      \end{equation*} \par
      Proof of \textnormal{(\romannumeral4)}. By direct calculation and \textnormal{(\romannumeral3)}, we derive that if $1\le i+j\le 2\Bbbk-1$ and $y\ge 1,$ 
      \begin{align*}
      	\left|\frac{\partial_y^i q}{y^j}\right|^2&\lesssim \left|\int_y^{+\infty}\partial_x \left(\frac{(\partial_x q)^2}{x^{2j}}\right)\,\mathrm{d}x\right|\\
      	&\lesssim \int_y^{+\infty} \frac{|\partial_x^i q|\cdot |\partial_x^{i+1}q|}{x^{d-3}\cdot x^{2j}}x^{d-3}\,\mathrm{d}x+\int_y^{+\infty} \frac{|\partial_x^i q|^2}{x^{d-3}\cdot x^{2j+1}}x^{d-3}\,\mathrm{d}x\\
      	&\lesssim \frac{1}{y^{d-4}}\left(\int_y^{+\infty} \frac{|\partial_x^i q|^2}{x^{2j+2}}+\int_y^{+\infty}\frac{|\partial_x^{i+1}q|^2}{x^{2j}}\right)\\
      	&\lesssim  \frac{1}{y^{d-4}}\cdot \left\{\begin{aligned}
      		&\mathscr{E}_{2m},\,\,\,\text{for}\,\,\, i+j+1=2m,\,\,\,1\le m\le \Bbbk.\\
      		&\sqrt{\mathscr{E}_{2m}}\sqrt{\mathscr{E}_{2(m+1)}},\,\,\,\text{for}\,\,\, i+j=2m,\,\,\,1\le m\le \Bbbk-1.
      	\end{aligned}\right.
      \end{align*}
\end{pf} 

\vspace{\baselineskip}
\section{Description of initial data and bootstrap assumption}\label{Description of initial data and bootstrap assumption}
Assumptions on initial data are as follows.
\begin{definition}\label{def of initial data}
	Denoting $\mathcal{V}:=P_l \mathcal{U}.$ Let $s_0\ge 1,$ Assuming initially\\
\textnormal{(\romannumeral1)} (Initial unstable modes)
\begin{equation}\label{assump on V2 to Vl initially}
	s_0^{\frac{\eta}{2}(1-\delta)}(\mathcal{V}_2(s_0),\cdots,\mathcal{V}_l(s_0))\in \mathcal{B}_{l-1}(0,1).
\end{equation}
\textnormal{(\romannumeral2)} (Initial stable modes) 
\begin{equation}\label{assump on V1 and bk from l+1 to L initially}
	|s_0^{\frac{\eta}{2}(1-\delta)}\mathcal{V}_1(s_0)|<1,\,\,\, |b_k(s_0)|<s_0^{-k-\frac{5l(2k-\gamma)}{2l-\gamma}}\,\,\,\text{for}\,\,\,l+1\le k
	\le L.
\end{equation}
\textnormal{(\romannumeral3)} (Initial energy)
\begin{equation}\label{assump on E2k initially}
	\sum\limits_{k=\hbar+2}^{\Bbbk} \mathscr{E}_{2k}(s_0)<s_0^{-\frac{10Ll}{2l-\gamma}}.
\end{equation}
\textnormal{(\romannumeral4)} up to a fixed rescaling, we may assume
\begin{equation}\label{assump on lamba initially}
	\lambda(s_0)=1.
\end{equation}	
\end{definition}

Next, we set up the bootstrap assumption as follows.
\begin{definition}\label{bootstrap assump}
Let $K\gg 1$ denote some large enough universal constant to be chosen later and $s\ge 1.$ Defining $\mathcal{S}_K(s)$ as the set of all $(b_1(s),\cdots,b_L(s),q(s))$ such that\\ 
\textnormal{(\romannumeral1)} (Control of unstable modes)
\begin{equation}\label{bootstrap assump on V2 to Vl}
	s^{\frac{\eta}{2}(1-\delta)}(\mathcal{V}_2(s),\cdots,\mathcal{V}_l(s))\in \mathcal{B}_{l-1}(0,1).
\end{equation}
\textnormal{(\romannumeral2)} (Control of stable modes)
\begin{equation}\label{bootstrap assump on V1 and bk from l+1 to L}
	|s^{\frac{\eta}{2}(1-\delta)}\mathcal{V}_1(s)|\le 10,\,\,\, |b_k(s)|\le \frac{10}{s^k}\,\,\,\text{for}\,\,\, l+1\le k\le L.
\end{equation}
\textnormal{(\romannumeral3)} (Control of highest order energy)
\begin{equation}\label{boootstrap assump on E2Bbbk}
	\mathscr{E}_{2\Bbbk}(s)\le Ks^{-[2L+2(1-\delta)(1+\eta)]}.
\end{equation}
\textnormal{(\romannumeral4)} (Control of lower order energy)
\begin{equation}\label{boootstrap assump on E2k lower}
	\mathscr{E}_{2m}(s)\le\left\{\begin{aligned}
		&Ks^{-\frac{l}{2l-\gamma}(4m-d+2)},\,\,\, \hbar+2\le m\le l+\hbar.\\
		&Ks^{-[2(m-\hbar-1)+2(1-\delta)-K\eta]},\,\,\, l+\hbar+1\le m\le \Bbbk-1.
	\end{aligned}\right.
\end{equation}   	
\end{definition}
\begin{rk}
	\normalfont
	Note that we have\\
     \textnormal{(\romannumeral1)} If $s_0$ is large enough, initial data $(b(s_0),q(s_0))\in \mathcal{S}_K(s_0).$\\
	\textnormal{(\romannumeral2)} If $(b(s),q(s))\in \mathcal{S}_K(s),$ then 
	\begin{equation*}
		b_k(s)\simeq \frac{c_k}{s^k},\,\,\, 1\le k\le l.
	\end{equation*}  
	In particular,
	\begin{equation*}
		b_1(s)\simeq \frac{c_1}{s},\,\,\, |b_k(s)|\lesssim b_1^k(s),\,\,\, 1\le k\le L.
	\end{equation*}
	Thus the assumptions in Proposition \ref{first approximation} are justified.
\end{rk}

\vspace{\baselineskip}
\section{Modulation estimates and improved bounds}
\begin{prop}\label{modulation estimates}
	For $K\gg 1$ some universal large constant, assuming there is $s_0(K)\gg 1$ such that $(b(s),q(s))\in \mathcal{S}_K(s)$ on $s\in [s_0,s_1]$ for some $s_1\ge s_0.$ Then for $s\in [s_0,s_1],$ we have
	\begin{equation}\label{modu esti for lambda and bk from 1 to L-1}
		\sum\limits_{k=1}^{L-1}|(b_k)_s+(2k-\gamma)b_1b_k-b_{k+1}|+\left|b_1+\frac{\lambda_s}{\lambda}\right|\lesssim b_1^{L+1+(1-\delta)(1+\eta)},
	\end{equation}
and
\begin{equation}\label{modu esti for bL}
	|(b_L)_s+(2L-\gamma)b_1b_L|\lesssim \frac{\sqrt{\mathscr{E}_{2\Bbbk}}}{M^{2\delta}}+b_1^{L+1+(1-\delta)(1+\eta)}.
\end{equation}
\end{prop}
\begin{pf}
	\normalfont
	The idea is to take inner product of (\ref{eq: eq of q}) with $\mathscr{L}^k \Phi_M$ for $1\le k\le L,$ make use of orthogonal condition (\ref{orthogonal condition}), apply H\"older and coercivity property of $\mathscr{L}.$\par 
	Denoting $D(t):=|b_1+\frac{\lambda_s}{\lambda}|+\sum\limits_{k=1}^L |(b_k)_s+(2k-\gamma)b_1b_k-b_{k+1}|.$ Now we take inner product of (\ref{eq: eq of q}) with $\mathscr{L}^L\Phi_M,$ and in view of (\ref{orthogonal condition}) we have
	\begin{equation}\label{innner product of eq of q with highest order of L}
		\langle \widehat{Mod}(t),\mathscr{L}^L\Phi_M\rangle=-\langle \mathscr{L}^L \widetilde{\Psi}_b,\Phi_M\rangle-\langle \mathscr{L}^{L+1}q,\Phi_M\rangle-\langle -\frac{\lambda_s}{\lambda}\Lambda q-\mathcal{H}(q)+\mathcal{N}(q),\mathscr{L}^L\Phi_M\rangle.
	\end{equation}
Then we estimate every term in (\ref{innner product of eq of q with highest order of L}). By (\ref{def of widehat Mod}),
\begin{equation}\label{inner prod of widehatmod and LLPhiM}
	\langle \widehat{Mod}(t),\mathscr{L}^L \Phi_M\rangle=O(D(t))\langle \Lambda \widetilde{Q}_b,\mathscr{L}^L \Phi_M\rangle+\langle Mod(t),\mathscr{L}^L \Phi_M\rangle.
\end{equation}
By the degrees of $T_k$ and $S_k$, we see that
\begin{align*}
	\langle \Lambda \widetilde{Q}_b,\mathscr{L}^L \Phi_M\rangle&=\sum\limits_{k=1}^L \langle b_k\Lambda T_k,\mathscr{L}^L\Phi_M\rangle+\sum\limits_{k=2}^{L+2}\langle \Lambda S_k,\mathscr{L}^L\Phi_M\rangle\\ &\lesssim \sum\limits_{k=1}^L b_1^kM^{d-2-2\gamma-2(L-k)}+\sum\limits_{k=2}^{L+2} b_1^kM^{d-2\gamma-2-2(L+2-k)}\lesssim b_1M^C,
\end{align*}
where we also have used the fact that $|c_{j,M}|\lesssim M^{2j},$ $0\le j\le L,$ which can be verified by induction on $j.$ By \textnormal{(\romannumeral1)} of Proposition \ref{first approximation} and (\ref{more property of Phi}),
\begin{equation*}
	\langle Mod(t),\mathscr{L}^L\Phi_M\rangle=(-1)^L\langle \chi_M\Lambda Q,\Lambda Q\rangle[(b_L)_s+(2L-\gamma)b_1b_L]+O\left(D(t)\sum\limits_{k=1}^L\sum\limits_{j=k+1}^{L+2}\left\langle \frac{\partial S_j}{\partial b_k},\mathscr{L}^L\Phi_M\right\rangle\right),
\end{equation*}
where $\left\langle \frac{\partial S_j}{\partial b_k},\mathscr{L}^L\Phi_M\right\rangle\lesssim b_1^{j-k}M^{d-2\gamma-2(L+2-j)}\lesssim b_1M^C,$ for $k+1\le j\le L+2,$ $1\le k\le L.$ Applying above estimates into (\ref{inner prod of widehatmod and LLPhiM}), we get
\begin{equation}\label{esti of inner prod of widehatmod and LLPhiM}
	\langle \widehat{Mod}(t),\mathscr{L}^L \Phi_M\rangle=(-1)^L\langle \chi_M\Lambda Q,\Lambda Q\rangle[(b_L)_s+(2L-\gamma)b_1b_L]+O(M^Cb_1D(t)).
\end{equation}
Using (\ref{3 esti of tilde Psib}) with $m=L-\hbar-1$ and H\"older, we have 
\begin{equation}\label{inner prod of LLwidetildePsib and PhiM}
	\langle \mathscr{L}^L \widetilde{\Psi}_b,\Phi_M\rangle\lesssim \left(\int_{y\le 2M}|\mathscr{L}^L\widetilde{\Psi}_b|^2\right)^{\frac{1}{2}}\left(\int_{y\le 2M} |\Phi_M|^2\right)^{\frac{1}{2}}\lesssim M^Cb_1^{L+3}.
\end{equation}
For the term $\langle \mathscr{L}^{L+1}q,\Phi_M\rangle,$ by Lemma \ref{coercivity of iterate of L},\begin{equation*}
	\mathscr{E}_{2\Bbbk}(q)\gtrsim \int \frac{|\mathscr{L}^{L+1}q|^2}{y^4(1+y^{4(\hbar-1)})}\gtrsim \int_{y\le 2M}\frac{|\mathscr{L}^{L+1}q|^2}{(1+y^{4\hbar})},
\end{equation*}
then by H\"older it follows that
\begin{equation}\label{inner prod of LL+1q and PhiM}
	|\langle \mathscr{L}^{L+1}q,\Phi_M\rangle|\lesssim M^{2\hbar}\left(\int_{y\le 2M}\frac{|\mathscr{L}^{L+1}q|^2}{(1+y^{4\hbar})}\right)^{\frac{1}{2}}\left(\int_{y\le 2M}|\Phi_M|^2\right)^{\frac{1}{2}}\lesssim M^{2\hbar+\frac{d-2}{2}-\gamma}\sqrt{\mathscr{E}_{2\Bbbk}(q)}.
\end{equation}
By triangle inequality, 
\begin{equation*}
	\left\langle -\frac{\lambda_s}{\lambda}\Lambda q,\mathscr{L}^L\Phi_M\right\rangle\lesssim (D(t)+b_1)\langle \Lambda q,\mathscr{L}^L\Phi_M\rangle.
\end{equation*}
Note that by H\"older,\begin{equation*}
	\langle \Lambda q,\mathscr{L}^L\Phi_M\rangle\lesssim \left\| \frac{\partial_y q}{y^2(1+y^{2(\Bbbk-2)+1})}\right\|_{L^2(y\le 2M)} \| y^3(1+y^{2(\Bbbk-2)+1})\mathscr{L}^L\Phi_M\|_{L^2(y\le 2M)}\lesssim M^C\sqrt{\mathscr{E}_{2\Bbbk}},
\end{equation*}
where we also have used the fact that by Lemma \ref{coercivity of iterate of L} and Lemma \ref{coercivity of L},
\begin{equation}\label{coercivity used in modulation estimates}
	\mathscr{E}_{2\Bbbk}(q)\gtrsim \int \frac{|\mathscr{L}q|^2}{y^4{1+y^{4(\Bbbk-2)}}}\gtrsim \int \frac{|\partial_y q|^2}{y^4(1+y^{4(\Bbbk-2)+2})}+\int \frac{q^2}{y^6(1+y^{4(\Bbbk-2)+2})}.
\end{equation} 
Therefore,\begin{equation}\label{inner prod of lamdaslambdaLambdaq and LLPhiM}
	\left\langle -\frac{\lambda_s}{\lambda}\Lambda q,\mathscr{L}^L\Phi_M\right\rangle\lesssim (D(t)+b_1)M^C\sqrt{\mathscr{E}_{2\Bbbk}}.
\end{equation}
Again by H\"older and (\ref{coercivity used in modulation estimates}), we derive
\begin{align}
	\langle -\mathcal{H}(q),\mathscr{L}^L \Phi_M\rangle&\lesssim \langle \frac{q}{y^2}|Q\Theta_{b}+\Theta_{b}+(\Theta_{b})^2|,\mathscr{L}^L \Phi_M\rangle\notag \\ &\lesssim \left\|\frac{q}{y^3(1+y^{2(\Bbbk-2)+1})}\right\|_{L^2}\|y(1+y^{2(\Bbbk-2)+1})(\Theta_{b}+\Theta_{b}Q+\Theta_{b}^2)\mathscr{L}^L\Phi_M\|_{L^2(y\le 2M)}\notag \\&\lesssim M^C b_1\sqrt{\mathscr{E}_{2\Bbbk}}.\label{inner prod of Hq and LLPhiM}
\end{align}
By (\ref{coercivity used in modulation estimates}), we obtain
\begin{align}
	\langle \mathcal{N}(q),\mathscr{L}^L \Phi_M\rangle&\lesssim \langle \frac{q^2}{y^2}(|Q-1|+|\Theta_{b}|+|q|),\mathscr{L}^L \Phi_M\rangle\notag\\ &\lesssim \int \frac{q^2}{y^6(1+y^{4(\Bbbk-2)+2})}\left\|y^4(1+y^{4(\Bbbk-2)+2})(|Q-1|+|\Theta_{b}|+|q|)\mathscr{L}^L \Phi_M\right\|_{L^{\infty}(y\le 2M)}\notag\\ &\lesssim M^C \mathscr{E}_{2\Bbbk}.\label{inner prod of Nq and LLPhiM}
\end{align}
Substituting (\ref{esti of inner prod of widehatmod and LLPhiM}), (\ref{inner prod of LLwidetildePsib and PhiM}), (\ref{inner prod of LL+1q and PhiM}), (\ref{inner prod of lamdaslambdaLambdaq and LLPhiM}), (\ref{inner prod of Hq and LLPhiM}) and (\ref{inner prod of Nq and LLPhiM}) into (\ref{innner product of eq of q with highest order of L}), we have
\begin{equation}\label{in pf modu esti on bL}
	|(b_L)_s+(2L-\gamma)b_1b_L|\lesssim \frac{\sqrt{\mathscr{E}_{2\Bbbk}}}{M^{2\delta}}+b_1^{L+3}+M^C b_1D(t).
\end{equation}\par 
Next we take inner product of (\ref{eq: eq of q}) with $\mathscr{L}^k\Phi_M$ for $1\le k\le L-1,$ we get 
\begin{equation*}
	\langle \widehat{Mod}(t),\mathscr{L}^k \Phi_M\rangle=-\langle \mathscr{L}^k\widetilde{\Psi}_b,\Phi_M\rangle-\langle -\frac{\lambda_s}{\lambda}\Lambda q-\mathcal{H}(q)+\mathcal{N}(q),\mathscr{L}^k \Phi_M\rangle,
\end{equation*}
then similar to the derivation of (\ref{in pf modu esti on bL}), we have
\begin{equation}\label{in pf modu esti on bk}
	|(b_k)_s+(2k-\gamma)b_1b_k-b_{k+1}|\lesssim b_1^{L+3}+M^Cb_1(\sqrt{\mathscr{E}_{2\Bbbk}}+D(t)).
\end{equation}\par 
Then we take inner product of (\ref{eq: eq of q}) with $\Phi_M,$ getting
\begin{equation*}
	\langle \widehat{Mod}(t),\Phi_M\rangle=-\langle \widetilde{\Psi}_b,\Phi_M\rangle-\langle -\frac{\lambda_s}{\lambda}\Lambda q-\mathcal{H}(q)+\mathcal{N}(q),\Phi_M\rangle.
\end{equation*} 
Note that by above estimates, \begin{equation*}
		\langle \widehat{Mod}(t),\Phi_M\rangle=-(\frac{\lambda_s}{\lambda}+b_1)\langle \Lambda Q,\chi_M\Lambda Q\rangle+O(M^Cb_1D(t)).
\end{equation*}
Thus again similar to the derivation of (\ref{in pf modu esti on bL}), we see that 
\begin{equation}\label{in pf mod esti on lambdaslambda+b1}
	\left|\frac{\lambda_s}{\lambda}+b_1\right|\lesssim b_1^{L+3}+M^C b_1(\sqrt{\mathscr{E}_{2\Bbbk}}+D(t)).
\end{equation}\par 
Now we sum up (\ref{in pf modu esti on bL})-(\ref{in pf mod esti on lambdaslambda+b1}) and apply (\ref{boootstrap assump on E2Bbbk}) to get
\begin{equation}\label{esti of Dt}
	D(t)\lesssim \frac{\sqrt{\mathscr{E}_{2\Bbbk}}}{M^{2\delta}}+b_1^{L+1+(1-\delta)(1+\eta)},
\end{equation}
then substituting (\ref{esti of Dt}) back into (\ref{in pf modu esti on bL})-(\ref{in pf mod esti on lambdaslambda+b1}), we get (\ref{modu esti for lambda and bk from 1 to L-1}) and (\ref{modu esti for bL}), this concludes the proof.
\end{pf}
\begin{rk}
	\normalfont
	We remark that by (\ref{modu esti for lambda and bk from 1 to L-1}) and (\ref{assump on b_k}), one has $|(b_1)_s|\lesssim b_1^2,$ this justifies the additional assumption in Proposition \ref{localized approximation}.
\end{rk}

For the $L$-th equation in the $b$-system, one needs improved modulation estimate in order to close bootstrap concerning stable modes, and we establish it as follows. 
\begin{prop}\label{improved bound for bL prop}
Under the assumptions in Proposition \ref{modulation estimates}, we have for all $s\in [s_0,s_1],$
\begin{equation}\label{improved estimate for b_L}
\left|(b_L)_s+(2L-\gamma)b_1b_L+(-1)^L \partial_s \frac{\langle \mathscr{L}^L q,\chi_{B_0}\Lambda Q\rangle}{\langle \chi_{B_0}\Lambda Q,\Lambda Q\rangle}\right|\lesssim \frac{1}{B_0^{2\delta}}\Big(\sqrt{\mathscr{E}_{2\Bbbk}}+b_1^{L+1+(1-\delta)-c_L\eta}\Big).
\end{equation}   
	For certainty, we assume $L\gg 1$ is an even integer.   
\end{prop}
\begin{pf}
	\normalfont
Now we take inner product of (\ref{eq: eq of q}) with $\mathscr{L}^L(\chi_{B_0}\Lambda Q),$ and get 
\begin{align}
&\quad\,\, \langle \chi_{B_0}\Lambda Q,\Lambda Q\rangle\left\{\frac{\mathrm{d}}{\mathrm{d}s}\left[\frac{\langle \mathscr{L}^L q,\chi_{B_0}\Lambda Q\rangle}{\langle\chi_{B_0}\Lambda Q,\Lambda Q\rangle}\right]-\langle \mathscr{L}^L q,\chi_{B_0}\Lambda Q\rangle \frac{\mathrm{d}}{\mathrm{d}s}\left[\frac{1}{\langle \Lambda Q,\chi_{B_0}\Lambda Q\rangle}\right]\right\}\notag \\ 	&=\langle \mathscr{L}^L q,\Lambda Q\partial_s \chi_{B_0}\rangle-\langle \mathscr{L}^{L+1}q,\chi_{B_0}\Lambda Q\rangle+\frac{\lambda_s}{\lambda}\langle \mathscr{L}^L\Lambda q,\chi_{B_0}\Lambda Q\rangle\notag \\ &\quad\,\, -\langle \mathscr{L}^L\widetilde{\Psi}_b,\chi_{B_0}\Lambda Q\rangle-\langle \mathscr{L}^L\widehat{Mod}(t),\chi_{B_0}\Lambda Q\rangle+\langle \mathscr{L}^L(\mathcal{H}(q)-\mathcal{N}(q)),\chi_{B_0}\Lambda Q\rangle.\label{eq of inner prod of eq of q and LLchiB0LambdaQ}
\end{align}
By H\"older and Lemma \ref{coercivity of iterate of L},\begin{equation*}
	|\langle \mathscr{L}^L q,\chi_{B_0}\Lambda Q\rangle|\lesssim \left(\int_{y\le 2B_0} \frac{|\mathscr{L}^L q|^2}{y^4+y^{4\hbar+4}}\right)^{\frac{1}{2}}\|\chi_{B_0}\Lambda Q\|_{L^2}B_0^{2\hbar+2}\lesssim B_0^{\frac{d-2}{2}-\gamma+2\hbar+2}\sqrt{\mathscr{E}_{2\Bbbk}},
\end{equation*}
where we also used the fact that $B_0^{d-2-2\gamma}\lesssim \langle \Lambda Q,\chi_{B_0}\Lambda Q\rangle\lesssim B_0^{d-2-2\gamma}.$ Note also that $|\partial_s \chi_{B_0}|=\left|\frac{y\partial_s b_1}{2B_0b_1}\chi'(\frac{y}{B_0})\right|\lesssim 1_{B_0\le y\le 2B_0}b_1,$ then it follows 
\begin{align}
	\left|\langle \mathscr{L}^L q,\chi_{B_0}\Lambda Q\rangle \frac{\mathrm{d}}{\mathrm{d}s}\left[\frac{1}{\langle \Lambda Q,\chi_{B_0}\Lambda Q\rangle}\right]\right|&=\left|-\frac{\langle \mathscr{L}^L q,\chi_{B_0}\Lambda Q\rangle}{{\langle \Lambda Q,\chi_{B_0}\Lambda Q\rangle}^2}\langle \Lambda Q,\Lambda Q\partial_s \chi_{B_0}\rangle\right| \notag \\ &\lesssim \frac{B_0^{\frac{d-2}{2}-\gamma+2\hbar+2}\sqrt{\mathscr{E}_{2\Bbbk}}}{B_0^{2(d-2)-4\gamma}} b_1\int_{B_0\le y\le 2B_0} y^{-2\gamma+d-3}\,\mathrm{d}y \notag\\&\lesssim \frac{\sqrt{\mathscr{E}_{2\Bbbk}}}{B_0^{2\delta}}.\label{esti for the lhs in the pf of improved bd}
\end{align}
Again by H\"older and Lemma \ref{coercivity of iterate of L}, we estimate
\begin{align}
	|\langle \mathscr{L}^L q,\Lambda Q\partial_s \chi_{B_0}\rangle|&\lesssim \left(\int \frac{|\mathscr{L}^L q|^2}{y^4+y^{4\hbar+4}}\right)^{\frac{1}{2}}\left(\int_{B_0\le y\le 2B_0} (y^4+y^{4\hbar+4})|\Lambda Q|^2\right)^{\frac{1}{2}}\left|\frac{\partial_s b_1}{b_1}\right|\notag\\ &\lesssim B_0^{\frac{d-2}{2}-\gamma+2\hbar}\sqrt{\mathscr{E}_{2\Bbbk}},\label{esti of inner product of LLq and LambdaQpartialschiB0}\end{align}
	and
\begin{align}
   |\langle \mathscr{L}^{L+1}q,\chi_{B_0}\Lambda Q\rangle|&\lesssim \left(\int \frac{|\mathscr{L}^{L+1}q|^2}{y^4+y^{4\hbar}}\right)^{\frac{1}{2}}\left(\int_{y\le 2 B_0} (y^4+y^{4\hbar})|\chi_{B_0}\Lambda Q|^2\right)^{\frac{1}{2}}\notag\\ &\lesssim B_0^{\frac{d-2}{2}-\gamma+2\hbar}\sqrt{\mathscr{E}_{2\Bbbk}}.\label{esti of inner product of LL+1q and chiB0LambdaQ}
\end{align}
By H\"older, (\ref{modu esti for lambda and bk from 1 to L-1}) and (\ref{coercivity used in modulation estimates}), we see that
\begin{align}
	\left|\frac{\lambda_s}{\lambda}\langle \mathscr{L}^L\Lambda q,\chi_{B_0}\Lambda Q\rangle\right|&\lesssim b_1\left(\int \frac{|\partial_y q|^2}{y^4+y^{4(L+\hbar)+2}}\right)^{\frac{1}{2}}\left(\int (y^6+y^{4(L+\hbar)+4})|\mathscr{L}^L(\chi_{B_0}\Lambda Q)|^2\right)^{\frac{1}{2}}\notag\\ &\lesssim B_0^{\frac{d-2}{2}-\gamma+2\hbar}\sqrt{\mathscr{E}_{2\Bbbk}}.\label{esti of inner product of lambdaslambdaLLLambdaq and chiB0LambdaQ}
\end{align}
By H\"older and (\ref{4 esti of tilde Psib}) with $m=L,$ we have
\begin{align}
	|\langle \mathscr{L}^L\widetilde{\Psi}_b,\chi_{B_0}\Lambda Q\rangle|&\lesssim \left(\int_{y\le 2 B_0} \frac{|\widetilde{\Psi}_b|^2}{1+y^{4(\hbar+L+1)}}\right)^{\frac{1}{2}}\left(\int (1+y^{4(\hbar+L+1)})|\mathscr{L}^L (\chi_{B_0}\Lambda Q)|^2\right)^{\frac{1}{2}}\notag\\ &\lesssim B_0^{\frac{d-2}{2}-\gamma+2\hbar}b_1^{L+1+(1-\delta)-c_L \eta}.\label{esti of inner product of LLwidetildePsib and chiB0LambdaQ}
\end{align}
Next we estimate the term $\langle \mathscr{L}^L(\mathcal{H}(q)-\mathcal{N}(q)),\chi_{B_0}\Lambda Q\rangle.$ Considering under the condition $y\le 2B_0,$ we get $|f'(\widetilde{Q}_b)-f'(Q)|\lesssim |Q\Theta_{b}+\Theta_{b}+\Theta_{b}^2|,$ where $|\Theta_{b}|\lesssim \sum\limits_{k=1}^L b_1^ky^{2k-\gamma}+\sum\limits_{k=2}^{L+2}b_1^ky^{2(k-1)-\gamma}\lesssim b_1^{\frac{\gamma}{2}}.$ Then $|f'(\widetilde{Q}_b)-f'(Q)|\lesssim b_1^{\frac{\gamma}{2}}\ll 1,$ and hence \begin{equation}\label{rougn bd on Hq when yle2B0}
	|\mathcal{H}(q)|\lesssim \frac{|q|}{y^2}.\,\,\,\text{(a rough bound is enough here)}
\end{equation}
Note that by \textnormal{(\romannumeral4)} of Lemma \ref{coercivity-determined esti on q}, when $1\le y\le 2B_0,$ $|q|\lesssim y^{2L+2-2\delta-\gamma}\mathscr{E}_{2\Bbbk}^{\frac{1}{2}}\lesssim b_1^{\frac{\gamma}{2}+(1-\delta)\eta}\ll 1,$ then
\begin{align*}
	|f(\widetilde{Q}_b+q)-f(\widetilde{Q}_b)-f'(\widetilde{Q}_b)q|&=\left|q^2[3(\widetilde{Q}_b-1)+q]\right|\\ &\lesssim q^2\left(|Q-1|+|\Theta_{b}|+|q|\right)\lesssim |q|^2.
\end{align*}
Therefore, \begin{equation}\label{rough bd on Nq when yle2B0}
	|\mathcal{N}(q)|\lesssim \frac{|q|^2}{y^2}\lesssim \frac{|q|}{y^2},\,\,\,\text{for}\,\,\,y\le 2B_0.\,\,\,\text{(a rough bound is enough here)}
\end{equation}
By (\ref{rougn bd on Hq when yle2B0}), (\ref{rough bd on Nq when yle2B0}), H\"older and (\ref{coercivity used in modulation estimates}), we see that 
\begin{align}
	|\langle \mathscr{L}^L(\mathcal{H}(q)-\mathcal{N}(q)),\chi_{B_0}\Lambda Q\rangle|&\lesssim \int \frac{|q|}{y^2}\left|\mathscr{L}^L(\chi_{B_0}\Lambda Q)\right|\notag\\ &\lesssim \left(\int \frac{|q|^2}{y^6+y^{4\Bbbk}}\right)^{\frac{1}{2}}\left(\int (y^2+y^{4\Bbbk-4})\left|\mathscr{L}^L(\chi_{B_0}\Lambda Q)\right|^2\right)^{\frac{1}{2}}\notag\\ &\lesssim B_0^{\frac{d-2}{2}-\gamma+2\hbar}\sqrt{\mathscr{E}_{2\Bbbk}}.\label{esti of inner product of LLHq-Nq and chiB0LambdaQ}
\end{align}
It remains to estimate $\langle \mathscr{L}^L\widehat{Mod}(t),\chi_{B_0}\Lambda Q\rangle,$ direct computation gives 
\begin{equation}\label{expression of inner product of LLwidehatMod and chiB0LambdaQ}
	\begin{aligned}
		\langle \mathscr{L}^L\widehat{Mod}(t),\chi_{B_0}\Lambda Q\rangle&=(-1)^L\langle \Lambda Q,\chi_{B_0}\Lambda Q\rangle[(b_L)_s+(2L-\gamma)b_1b_L]\\ &\quad\,-\left(\frac{\lambda_s}{\lambda}+b_1\right)\langle \Lambda \Theta_{b},\mathscr{L}^L(\chi_{B_0}\Lambda Q)\rangle\\ &\quad\,+\left\langle \sum_{k=1}^L[(b_k)_s+(2k-\gamma)b_1b_k-b_{k+1}]\sum\limits_{j=k+1}^{L+2}\frac{\partial S_j}{\partial b_k},\mathscr{L}^L(\chi_{B_0}\Lambda Q)\right\rangle.
	\end{aligned}
\end{equation}
When $y\le 2B_0,$ we have $b_1y^2\lesssim 1,$ then
\begin{equation*}
	\left\{\begin{aligned}
		\sum_{k=1}^L|b_k\Lambda T_k|&\lesssim \sum_{k=1}^Lb_1^ky^{2k-\gamma}\lesssim b_1y^{2-\gamma}\\ \sum_{k=2}^{L+2}|\Lambda S_k|&\lesssim \sum_{k=2}^{L+2}b_1^ky^{2(k-1)-\gamma}\lesssim b_1^2y^{2-\gamma}
	\end{aligned}\right. \Longrightarrow |\Lambda \Theta_{b}|\lesssim b_1y^{2-\gamma},
\end{equation*}
and\begin{equation*}
	\sum\limits_{j=k+1}^{L+2}\frac{\partial S_j}{\partial b_k}\lesssim \sum\limits_{j=k+1}^{L+2}b_1^{j-k}y^{2(j-1)-\gamma}\lesssim b_1y^{2k-\gamma}.
\end{equation*}
Then by Proposition \ref{modulation estimates},
\begin{align}
	\left|\frac{\lambda_s}{\lambda}+b_1\right|\left|\langle \Lambda \Theta_{b},\mathscr{L}^L(\chi_{B_0}\Lambda Q)\rangle\right|&\lesssim b_1^{L+1+(1-\delta)(1+\eta)}\int_{0}^{2B_0}b_1y^{2-\gamma}\cdot y^{-\gamma-2L+d-3}\,\mathrm{d}y\notag\\ &\lesssim b_1^{2L+1+(1-\delta)(1+\eta)}B_0^{d-2-2\gamma},\label{esti of inner product of lambdaslambdab1LambdaThetab and LLchiB0LambdaQ}
\end{align}	 
and 
\begin{align}
	&\quad\,\sum_{k=1}^L\left|(b_k)_s+(2k-\gamma)b_1b_k-b_{k+1}\right|\left|\left\langle\sum\limits_{j=k+1}^{L+2}\frac{\partial S_j}{\partial b_k},\mathscr{L}^L(\chi_{B_0}\Lambda Q)\right\rangle\right|\notag\\&\lesssim \left(\frac{\sqrt{\mathscr{E}_{2\Bbbk}}}{M^{2\delta}}+b_1^{L+1+(1-\delta)(1+\eta)}\right)\int_{0}^{2B_0}b_1y^{2L-\gamma}\cdot y^{-\gamma-2L+d-3}\,\mathrm{d}y\notag\\ &\lesssim \left(\frac{\sqrt{\mathscr{E}_{2\Bbbk}}}{M^{2\delta}}+b_1^{L+1+(1-\delta)(1+\eta)}\right)b_1B_0^{d-2-2\gamma}.\label{esti of inner product of bks+2k-gammab1bk-bk+1partialSjpartialbk and LLchiB0LambdaQ}
\end{align}
Now substituting (\ref{esti of inner product of lambdaslambdab1LambdaThetab and LLchiB0LambdaQ}) and (\ref{esti of inner product of bks+2k-gammab1bk-bk+1partialSjpartialbk and LLchiB0LambdaQ}) into (\ref{expression of inner product of LLwidehatMod and chiB0LambdaQ}), then gathering the estimates (\ref{esti for the lhs in the pf of improved bd})-(\ref{esti of inner product of LLwidetildePsib and chiB0LambdaQ}), (\ref{esti of inner product of LLHq-Nq and chiB0LambdaQ}) and (\ref{expression of inner product of LLwidehatMod and chiB0LambdaQ}) into (\ref{eq of inner prod of eq of q and LLchiB0LambdaQ}) and dividing $(-1)^L\langle \Lambda Q,\chi_{B_0}\Lambda Q\rangle,$ we would get (\ref{improved estimate for b_L}), this concludes the proof.
\end{pf}

\vspace{\baselineskip}
\section{Energy estimates}
\begin{prop}\label{energy estimates}
	Under the assumptions in Proposition \ref{modulation estimates}, there are monotonicity formulas as follows. For highest order energy, we have 
	\begin{equation}\label{monotonicity for E2Bbbk}
		\frac{\mathrm{d}}{\mathrm{d}t}\left\{\frac{\mathscr{E}_{2\Bbbk}}{\lambda^{4\Bbbk-d+2}}[1+O(b_1^{\eta(1-\delta)})]\right\}\lesssim \frac{b_1}{\lambda^{4\Bbbk-d+4}}\left[b_1^{L+(1-\delta)(1+\eta)}\sqrt{\mathscr{E}_{2\Bbbk}}+\frac{\mathscr{E}_{2\Bbbk}}{M^{2\delta}}+b_1^{2L+2(1-\delta)(1+\eta)}\right].
	\end{equation}
And for lower order energies, for all $\hbar+2\le m\le \Bbbk-1,$ we have
\begin{equation}\label{monotoncity for E2m}
	\frac{\mathrm{d}}{\mathrm{d}t} \left\{\frac{\mathscr{E}_{2m}}{\lambda^{4m-d+2}}[1+O(b_1)]\right\}\lesssim \frac{b_1}{\lambda^{4m-d+4}}\left[b_1^{m-\hbar-1+(1-\delta)-C\eta}\sqrt{\mathscr{E}_{2m}}+b_1^{2(m-\hbar-1)+2(1-\delta)-C\eta}\right].
\end{equation}
\end{prop}
\begin{pf}
	\normalfont
	Some aspects of this proof are parallel to the proof of Proposition 4.4 in \cite{ghoul2018stability}, thus we shall omit some details. For simplicity, we shall only prove (\ref{monotonicity for E2Bbbk}). For convenience, we abuse notation by abbreviating $v_k^*$ as $v_k$ for $k\in \mathbb{N}.$\par 
	Firstly, we set up the energy identity. Acting $\mathscr{L}_{\lambda}^{\Bbbk-1}$ on equation (\ref{eq: eq of v}), we have
\begin{equation}\label{eq:act LlambdaBbbk-1 on eq of v}
	\partial_t v_{2\Bbbk-2}+\mathscr{L}_{\lambda}v_{2\Bbbk-2}=[\partial_t,\mathscr{L}_{\lambda}^{\Bbbk-1}]v+\mathscr{L}_{\lambda}^{\Bbbk-1}\left(\frac{1}{\lambda^2}\mathcal{F}_{\lambda}\right).	
\end{equation}
Then acting $\mathscr{A}_{\lambda}$ on (\ref{eq:act LlambdaBbbk-1 on eq of v}), we get
\begin{equation}\label{eq:act AlambdaLlambdaBbbk-1 on eq of v}
	\partial_t v_{2\Bbbk-1}+\widetilde{\mathscr{L}}v_{2\Bbbk-1}=\frac{\partial_t V_{\lambda}}{r}v_{2\Bbbk-2}+\mathscr{A}_{\lambda}[\partial_t,\mathscr{L}_{\lambda}^{\Bbbk-1}]v+\mathscr{A}_{\lambda}\mathscr{L}_{\lambda}^{\Bbbk-1}\left(\frac{1}{\lambda^2}\mathcal{F}_{\lambda}\right).
\end{equation}
Making use of (\ref{eq:act LlambdaBbbk-1 on eq of v}) and (\ref{eq:act AlambdaLlambdaBbbk-1 on eq of v}), we still have the energy identity as 
\begin{align}
	&\quad\,\frac{1}{2}\frac{\mathrm{d}}{\mathrm{d}t}\left\{\frac{\mathscr{E}_{2\Bbbk}}{\lambda^{4\Bbbk-d+2}}+2\int \frac{b_1(\Lambda V)_{\lambda}}{\lambda^2 r}v_{2\Bbbk-1}v_{2\Bbbk-2}\right\}\notag\\ &=-\int |\widetilde{\mathscr{L}}_{\lambda}v_{2\Bbbk-1}|^2-\left(\frac{\lambda_s}{\lambda}+b_1\right)\int \frac{(\Lambda \widetilde{Z})_{\lambda}}{2\lambda^2 r^2}v_{2\Bbbk-1}^2-\int \frac{b_1(\Lambda V)_{\lambda}}{\lambda^2 r}v_{2\Bbbk-2}\widetilde{\mathscr{L}}_{\lambda}v_{2\Bbbk-1}\notag\\ &\quad\,\, +\int \frac{\mathrm{d}}{\mathrm{d}t}\left(\frac{b_1(\Lambda V)_{\lambda}}{\lambda^2 r}\right)v_{2\Bbbk-1}v_{2\Bbbk-2}+\int \frac{b_1(\Lambda V)_{\lambda}}{\lambda^2 r}v_{2\Bbbk-1}\left([\partial_t,\mathscr{L}_{\lambda}^{\Bbbk-1}]v+\mathscr{L}_{\lambda}^{\Bbbk-1}\left(\frac{1}{\lambda^2}\mathcal{F}_{\lambda}\right)\right)\notag\\ &\quad\,\, +\int \left(\widetilde{\mathscr{L}}_{\lambda}v_{2\Bbbk-1}+\frac{b_1(\Lambda V)_{\lambda}}{\lambda^2 r}v_{2\Bbbk-2}\right)\left(\frac{\partial_t V_{\lambda}}{r}v_{2\Bbbk-2}+\mathscr{A}_{\lambda}[\partial_t,\mathscr{L}_{\lambda}^{\Bbbk-1}]v+\mathscr{A}_{\lambda}\mathscr{L}_{\lambda}^{\Bbbk-1}\left(\frac{1}{\lambda^2}\mathcal{F}_{\lambda}\right)\right).\label{energy identity}
\end{align}\par
Now we estimate terms in (\ref{energy identity}). Note that by Lemma \ref{coercivity of iterate of L}, we get
\begin{equation}\label{coercivity used in energy estimate}
	\mathscr{E}_{2\Bbbk}(q)\gtrsim \int \frac{|q_{2\Bbbk-1}|^2}{y^2}+\sum\limits_{j=0}^{\Bbbk-1}\int \frac{|q_{2j}|^2}{y^4(1+y^{4(\Bbbk-1-j)})}+\sum\limits_{j=0}^{\Bbbk-2}\int \frac{|q_{2j+1}|^2}{y^6(1+y^{4(\Bbbk-2-j)})}.
\end{equation}\par 
On the second term on the LHS of (\ref{energy identity}). Note that by (\ref{def: def of V}),
\begin{equation*}
	|\Lambda V(y)|\lesssim \left\{\begin{aligned}
		&y^2\,\,\,\text{as}\,\,\,y\rightarrow 0\\ &y^{-2\gamma}\lesssim y^{-2}\,\,\,\text{as}\,\,\,y\rightarrow \infty
	\end{aligned}\right. \Longrightarrow |\Lambda V(y)|\lesssim \frac{y^2}{1+y^4}.
\end{equation*}
Then by H\"older and (\ref{coercivity used in energy estimate}), we estimate
\begin{align}
	\left|\int \frac{b_1(\Lambda V)_{\lambda}}{\lambda^2 r}v_{2\Bbbk-1}v_{2\Bbbk-2}\right|&=\frac{1}{\lambda^{4\Bbbk-d+2}}\left|\int \frac{b_1\Lambda V}{y}q_{2\Bbbk-1}q_{2\Bbbk-2}\right|\notag\\ &\lesssim \frac{b_1}{\lambda^{4\Bbbk-d+2}}\left(\int \frac{|q_{2\Bbbk-1}|^2}{y^2} \right)^{\frac{1}{2}}\left(\int \frac{|q_{2\Bbbk-2}|^2}{1+y^4}\right)^{\frac{1}{2}}\lesssim \frac{b_1}{\lambda^{4\Bbbk-d+2}}\mathscr{E}_{2\Bbbk}(q).\label{in energy identity lhs esti}
\end{align}\par   
Note that by (\ref{def: def of V}),\begin{equation*}
	\Lambda \widetilde{Z}=2V\Lambda V+(d-2)\Lambda V-\Lambda^2 V\lesssim \left\{\begin{aligned}
		&y^2\,\,\,\text{as}\,\,\,y\rightarrow 0\\ &y^{-2\gamma}\,\,\,\text{as}\,\,\,y\rightarrow \infty
	\end{aligned}\right. \lesssim \frac{y^2}{1+y^4}. 
\end{equation*}
Then by (\ref{modu esti for lambda and bk from 1 to L-1}) and (\ref{coercivity used in energy estimate}), we have 
\begin{align}
	\left|\left(\frac{\lambda_s}{\lambda}+b_1\right)\int \frac{(\Lambda \widetilde{Z})_{\lambda}}{\lambda^2 r^2}v_{2\Bbbk-1}^2\right|&=\left|\left(\frac{\lambda_s}{\lambda}+b_1\right)\frac{1}{\lambda^{4\Bbbk-d+4}}\int \frac{\Lambda \widetilde{Z}}{y^2}q_{2\Bbbk-1}^2\right|\notag\\ &\lesssim \frac{b_1^{L+1+(1-\delta)(1+\eta)}}{\lambda^{4\Bbbk-d+4}}\int \frac{q_{2\Bbbk-1}^2}{y^2}\lesssim \frac{b_1^{L+1+(1-\delta)(1+\eta)}}{\lambda^{4\Bbbk-d+4}}\mathscr{E}_{2\Bbbk}(q).\label{in energy identity rhs 1 esti}
\end{align}\par 
Again by (\ref{def: def of V}) and (\ref{coercivity used in energy estimate}), we get
\begin{align}
	\left|\int \frac{b_1(\Lambda V)_{\lambda}}{\lambda^2 r}v_{2\Bbbk-2}\widetilde{\mathscr{L}}_{\lambda}v_{2\Bbbk-1}\right|&\le \frac{1}{4}\int |\widetilde{\mathscr{L}}_{\lambda}v_{2\Bbbk-1}|^2+\int \frac{b_1^2|(\Lambda V)_{\lambda}|^2}{\lambda^4r^2}v_{2\Bbbk-2}^2\notag\\ &=\frac{1}{4}\int |\widetilde{\mathscr{L}}_{\lambda}v_{2\Bbbk-1}|^2+\frac{b_1^2}{\lambda^{4\Bbbk-d+4}}\int \frac{|\Lambda V(y)|^2}{y^2}q_{2\Bbbk-2}^2\notag\\ &\le \frac{1}{4}\int |\widetilde{\mathscr{L}}_{\lambda}v_{2\Bbbk-1}|^2+\frac{Cb_1^2}{\lambda^{4\Bbbk-d+4}}\mathscr{E}_{2\Bbbk}.\label{in energy identity rhs 2 esti}
\end{align}\par 
Note that by (\ref{modu esti for lambda and bk from 1 to L-1}), we derive
\begin{align*}
	\left|\frac{\mathrm{d}}{\mathrm{d}t}\left(\frac{b_1(\Lambda V)_{\lambda}}{\lambda^2}\right)\right|&=\left|\frac{(b_1)_s}{\lambda}(\Lambda V)_{\lambda}-\frac{b_1}{\lambda^4}\frac{\lambda_s}{\lambda}(\Lambda^2 V)_{\lambda}-\frac{2b_1(\Lambda V)_{\lambda}}{\lambda^4}\frac{\lambda_s}{\lambda}\right|\\ &\lesssim \frac{b_1^2}{\lambda^4}\left(|(\Lambda V)_{\lambda}|+|(\Lambda^2 V)_{\lambda}|\right).
\end{align*}
Then again by (\ref{def: def of V}), H\"older and (\ref{coercivity used in energy estimate}), we have 
\begin{align}
	\left|\int \frac{\mathrm{d}}{\mathrm{d}t}\left(\frac{b_1(\Lambda V)_{\lambda}}{\lambda^2 r}\right)v_{2\Bbbk-1}v_{2\Bbbk-2}\right|&\lesssim \frac{b_1^2}{\lambda^{4\Bbbk-d+4}}\int \frac{|\Lambda V|+|\Lambda^2 V|}{y}q_{2\Bbbk-1}q_{2\Bbbk-2}\notag\\ &\lesssim \frac{b_1^2}{\lambda^{4\Bbbk-d+4}}\left(\frac{|q_{2\Bbbk-1}|^2}{y^2}\right)^{\frac{1}{2}}\left(\frac{|q_{2\Bbbk-2}|^2}{y^4}\right)^{\frac{1}{2}}\lesssim \frac{b_1^2}{\lambda^{4\Bbbk-d+4}}\mathscr{E}_{2\Bbbk}.\label{in energy identity rhs 3 esti}
\end{align}\par 
Note that by (\ref{modu esti for lambda and bk from 1 to L-1}),
\begin{equation*}
	\partial_t V_{\lambda}=-\frac{\lambda_s}{\lambda}\frac{1}{\lambda^2}(\Lambda V)_{\lambda}\Longrightarrow \left|\frac{\partial_t V_{\lambda}}{r}\right|\lesssim \frac{b_1|(\Lambda V)_{\lambda}|}{\lambda^2r}.
\end{equation*}
Then again using (\ref{def: def of V}) and (\ref{coercivity used in energy estimate}), we compute
\begin{align}
	&\quad\,\,\left|\int \left(\widetilde{\mathscr{L}}_{\lambda}v_{2\Bbbk-1}+\frac{b_1(\Lambda V)_{\lambda}}{\lambda^2 r}v_{2\Bbbk-2}\right)\frac{\partial_t V_{\lambda}}{r}v_{2\Bbbk-2}\right|\notag\\&\le \frac{1}{4}\int |\widetilde{\mathscr{L}}_{\lambda}v_{2\Bbbk-1}|^2+C\int \left(\frac{b_1(\Lambda V)_{\lambda}}{\lambda^2r}\right)^2v_{2\Bbbk-2}^2\notag\\ &=\frac{1}{4}\int |\widetilde{\mathscr{L}}_{\lambda}v_{2\Bbbk-1}|^2+C\frac{b_1^2}{\lambda^{4\Bbbk-d+4}}\int \frac{|\Lambda V|^2}{y^2}q_{2\Bbbk-2}^2\notag\\ &\le \frac{1}{4}\int |\widetilde{\mathscr{L}}_{\lambda}v_{2\Bbbk-1}|^2+C\frac{b_1^2}{\lambda^{4\Bbbk-d+4}}\int \frac{q_{2\Bbbk-2}^2}{y^4}\notag\\ &\le \frac{1}{4}\int |\widetilde{\mathscr{L}}_{\lambda}v_{2\Bbbk-1}|^2+C\frac{b_1^2}{\lambda^{4\Bbbk-d+4}}\mathscr{E}_{2\Bbbk}(q).\label{in energy identity rhs 4 esti}
\end{align}\par 
Similar to the estimate of (\ref{in energy identity rhs 4 esti}), we have
\begin{align}
	&\quad\,\,\left|\int \frac{b_1(\Lambda V)_{\lambda}}{\lambda^2 r}v_{2\Bbbk-1}[\partial_t,\mathscr{L}_{\lambda}^{\Bbbk-1}]v\right|+\left|\int \left(\widetilde{\mathscr{L}}_{\lambda}v_{2\Bbbk-1}+\frac{b_1(\Lambda V)_{\lambda}}{\lambda^2 r}v_{2\Bbbk-2}\right)\mathscr{A}_{\lambda}[\partial_t,\mathscr{L}_{\lambda}^{\Bbbk-1}]v\right|\notag\\ &\le C\Bigg[\int \frac{b_1^2}{\lambda^2}\frac{v_{2\Bbbk-1}^2}{r^2}+\int \frac{|[\partial_t,\mathscr{L}_{\lambda}^{\Bbbk-1}]v|^2|(\Lambda V)_{\lambda}|^2}{\lambda^2}+\int |\mathscr{A}_{\lambda}[\partial_t,\mathscr{L}_{\lambda}^{\Bbbk-1}]v|^2\notag\\ &\quad\,\,+\int \left(\frac{b_1(\Lambda V)_{\lambda}}{\lambda^2r}\right)^2v_{2\Bbbk-2}^2\Bigg]+\frac{1}{4}\int |\widetilde{\mathscr{L}}_{\lambda}v_{2\Bbbk-1}|^2\notag\\ &\le \frac{1}{4}\int |\widetilde{\mathscr{L}}_{\lambda}v_{2\Bbbk-1}|^2+C\left(\frac{b_1^2}{\lambda^{4\Bbbk-d+4}}\mathscr{E}_{2\Bbbk}+\int \frac{|[\partial_t,\mathscr{L}_{\lambda}^{\Bbbk-1}]v|^2}{\lambda^2(1+y^4)}+\int |\mathscr{A}_{\lambda}[\partial_t,\mathscr{L}_{\lambda}^{\Bbbk-1}]v|^2\right).\label{in energy identity rhs 5 esti}
\end{align}
In (\ref{in energy identity rhs 5 esti}), we claim that 
\begin{equation}\label{commutator related estimate}
	\int \frac{|[\partial_t,\mathscr{L}_{\lambda}^{\Bbbk-1}]v|^2}{\lambda^2(1+y^2)}+\int |\mathscr{A}_{\lambda}[\partial_t,\mathscr{L}_{\lambda}^{\Bbbk-1}]v|^2\lesssim \frac{b_1^2}{\lambda^{4\Bbbk-d+4}}\mathscr{E}_{2\Bbbk}.
\end{equation}
We changed $y^4$ into $y^2$ so that the two integrals in (\ref{commutator related estimate}) are of the same order. Let us prove this claim. Note that 
\begin{equation}\label{reduction of the commutator}
	[\partial_t,\mathscr{L}_{\lambda}^{k-1}]g=\sum_{m=0}^{k-2}\mathscr{L}_{\lambda}^m[\partial_t,\mathscr{L}_{\lambda}]\mathscr{L}_{\lambda}^{k-2-m}g,
\end{equation}
for any $k\ge 2$ and any smooth radial function $g.$ One can prove it by induction on $k,$ we shall omit the details. Then 
\begin{equation*}
	[\partial_t,\mathscr{L}_{\lambda}^{\Bbbk-1}]v=\sum_{m=0}^{\Bbbk-2}\mathscr{L}_{\lambda}^m\left(\frac{\partial_t Z_{\lambda}}{r^2}\mathscr{L}_{\lambda}^{\Bbbk-2-m}v\right),\,\,\,\text{with}\,\,\,\frac{\partial_t Z_{\lambda}}{r^2}=-\frac{\lambda_s}{\lambda}\frac{(\Lambda Z)_{\lambda}}{\lambda^2 r^2}.
\end{equation*}
By (\ref{modu esti for lambda and bk from 1 to L-1}), we deduce that
\begin{equation*}
	\int \frac{|[\partial_t,\mathscr{L}_{\lambda}^{\Bbbk-1}]v|^2}{\lambda^2(1+y^2)}\lesssim \frac{b_1^2}{\lambda^{4\Bbbk-d+4}}\sum_{m=0}^{\Bbbk-2}\int \frac{1}{1+y^2}\left|\mathscr{L}^m\left(\frac{\Lambda Z}{y^2}\mathscr{L}^{\Bbbk-2-m}q\right)\right|^2.
\end{equation*}
When $m=0,$ note that 
\begin{equation*}
	\left|\frac{\Lambda Z}{y^2}\right|=\left|\frac{(d-2)f''(Q)\Lambda Q}{y^2}\right|\lesssim \frac{1}{y^{2\gamma+2}}\lesssim \frac{1}{1+y^4}.
\end{equation*}
Then by (\ref{coercivity used in energy estimate}),
\begin{equation*}
	\int \frac{1}{1+y^2}\left|\frac{\Lambda Z}{y^2}\mathscr{L}^{\Bbbk-2}q\right|^2\lesssim \int \frac{q_{2\Bbbk-4}^2}{1+y^{10}}\lesssim \mathscr{E}_{2\Bbbk}(q).
\end{equation*}
When $1\le m\le \Bbbk-2,$ by (\ref{leibniz for iterate of L}) with $\phi=\frac{\Lambda Z}{y^2},$ $g=\mathscr{L}^{\Bbbk-2-m}q,$ we have
\begin{align*}
\mathscr{L}^m\left(\frac{\Lambda Z}{y^2}\mathscr{L}^{\Bbbk-2-m}q\right)&=\sum\limits_{i=0}^m \mathscr{L}^{\Bbbk-2-(m-i)}q\left(\frac{\Lambda Z}{y^2}\right)_{2m,2i}+\sum\limits_{i=0}^{m-1}\mathscr{A}\mathscr{L}^{\Bbbk-2-(m-i)}q\left(\frac{\Lambda Z}{y^2}\right)_{2m,2i+1}\\&\lesssim \sum\limits_{i=0}^m \frac{q_{2(\Bbbk-2-(m-i))}}{1+y^{2\gamma+2+2(m-i)}}+\sum\limits_{i=0}^{m-1}\frac{q_{2(\Bbbk-2-(m-i))+1}}{1+y^{2\gamma+2+2(m-i)-1}},
\end{align*}
where we also used the fact that \begin{equation*}\left|\left(\frac{\Lambda Z}{y^2}\right)_{2m,i}\right|\lesssim \frac{1}{y^{2\gamma+2+2m-i}},\,\,\, \text{for}\,\,\, 0\le i\le 2m.\end{equation*}
Then by (\ref{coercivity used in energy estimate}), we see that
\begin{align*}
	\int \frac{1}{1+y^2}\left|\mathscr{L}^m\left(\frac{\Lambda Z}{y^2}\mathscr{L}^{\Bbbk-2-m}q\right)\right|^2&\lesssim \sum\limits_{i=0}^m \int \frac{q_{2(\Bbbk-2-(m-i))}^2}{1+y^{4\gamma+6+4(m-i)}}+\sum\limits_{i=0}^{m-1} \int \frac{q_{2(\Bbbk-2-(m-i))+1}^2}{1+y^{4\gamma+4+4(m-i)}}\\ &\lesssim \sum\limits_{i=0}^m \int \frac{q_{2(\Bbbk-2-(m-i))}^2}{y^4(1+y^{4(1+m-i)})}+\sum\limits_{i=0}^{m-1}\int \frac{q_{2(\Bbbk-2-(m-i))+1}^2}{y^6(1+y^{4(m-i)})}\lesssim \mathscr{E}_{2\Bbbk}(q).
\end{align*}
This concludes the proof of (\ref{commutator related estimate}).\par 
Again by (\ref{def: def of V}), H\"older and (\ref{coercivity used in energy estimate}), we get
\begin{align}
	\left|\int \frac{b_1(\Lambda V)_{\lambda}}{\lambda^2 r}v_{2\Bbbk-1}\mathscr{L}_{\lambda}^{\Bbbk-1}\left(\frac{1}{\lambda^2}\mathcal{F}_{\lambda}\right)\right|&=\frac{b_1}{\lambda^{4\Bbbk-d+4}}\int \frac{\Lambda V}{y}q_{2\Bbbk-1}\mathscr{L}^{\Bbbk-1}\mathcal{F}\notag\\ &\lesssim \frac{b_1}{\lambda^{4\Bbbk-d+4}}\left(\int \frac{q_{2\Bbbk-1}^2}{y^2}\right)^{\frac{1}{2}}\left(\int \frac{|\mathscr{L}^{\Bbbk-1}\mathcal{F}|^2}{1+y^4}\right)^{\frac{1}{2}}\notag\\ &\lesssim \frac{b_1}{\lambda^{4\Bbbk-d+4}}\sqrt{\mathscr{E}_{2\Bbbk}}\left(\int \frac{|\mathscr{L}^{\Bbbk-1}\mathcal{F}|^2}{1+y^4}\right)^{\frac{1}{2}}.\label{in energy identity rhs 6 esti}
\end{align}\par 
Similar to the estimate (\ref{in energy identity rhs 6 esti}), we have
\begin{equation}\label{in energy identity rhs 7 esti}
	\left|\int \frac{b_1(\Lambda V)_{\lambda}}{\lambda^2 r}v_{2\Bbbk-2}\mathscr{A}_{\lambda}\mathscr{L}_{\lambda}^{\Bbbk-1}\left(\frac{1}{\lambda^2}\mathcal{F}_{\lambda}\right)\right|\lesssim \frac{b_1}{\lambda^{4\Bbbk-d+4}}\sqrt{\mathscr{E}_{2\Bbbk}}\left(\int \frac{|\mathscr{A}\mathscr{L}^{\Bbbk-1}\mathcal{F}|^2}{1+y^2}\right)^{\frac{1}{2}}.
\end{equation}\par 
On the term $\int \widetilde{\mathscr{L}}_{\lambda}v_{2\Bbbk-1}\mathscr{A}_{\lambda}\mathscr{L}_{\lambda}^{\Bbbk-1}\left(\frac{1}{\lambda^2}\mathcal{F}_{\lambda}\right).$ Denoting \begin{equation*}
	\xi_L:=\frac{\langle \mathscr{L}^L q,\chi_{B_0}\Lambda Q\rangle}{\langle \chi_{B_0}\Lambda Q,\Lambda Q\rangle}\widetilde{T}_L,
\end{equation*}
we set up the decomposition 
\begin{equation*}
	\mathcal{F}=:\partial_s \xi_L+\mathcal{F}_0+\mathcal{F}_1,\,\,\,\text{where}\,\,\,\mathcal{F}_0:=-\widetilde{\Psi}_b-\widehat{Mod}-\partial_s \xi_L\,\,\,\text{and}\,\,\,\mathcal{F}_1:=\mathcal{H}(q)-\mathcal{N}(q).
\end{equation*}
Then by H\"older, we obtain
\begin{align}
	&\quad\,\,\int \widetilde{\mathscr{L}}_{\lambda}v_{2\Bbbk-1}\mathscr{A}_{\lambda}\mathscr{L}_{\lambda}^{\Bbbk-1}\left(\frac{1}{\lambda^2}\mathcal{F}_{\lambda}\right)\notag\\ &=\frac{1}{\lambda^{4\Bbbk-d+4}}\int \widetilde{\mathscr{L}}q_{2\Bbbk-1}\mathscr{A}\mathscr{L}^{\Bbbk-1}(\partial_s \xi_L+\mathcal{F}_0+\mathcal{F}_1)\notag\\ &=\frac{1}{\lambda^{4\Bbbk-d+4}}\left(\int \mathscr{A}^*q_{2\Bbbk-1}\mathscr{L}^{\Bbbk}(\partial_s \xi_L)+\int \mathscr{A}^*q_{2\Bbbk-1}\mathscr{L}^{\Bbbk}\mathcal{F}_0+\int \widetilde{\mathscr{L}}q_{2\Bbbk-1}\mathscr{A}\mathscr{L}^{\Bbbk-1}\mathcal{F}_1\right)\notag\\ &\le \frac{1}{\lambda^{4\Bbbk-d+4}}\int \mathscr{L}^{\Bbbk}q\mathscr{L}^{\Bbbk}(\partial_s \xi_L)+\frac{1}{\lambda^{4\Bbbk-d+4}}\left(\int |\mathscr{L}^{\Bbbk}q|^2\right)^{\frac{1}{2}}\left(\int |\mathscr{L}^{\Bbbk}\mathcal{F}_0|^2\right)^{\frac{1}{2}}\notag\\&\quad\,\,+\frac{1}{8}\frac{1}{\lambda^{4\Bbbk-d+4}}\int |\widetilde{\mathscr{L}}q_{2\Bbbk-1}|^2+\frac{2}{\lambda^{4\Bbbk-d+4}}\int |\mathscr{A}\mathscr{L}^{\Bbbk-1}\mathcal{F}_1|^2\notag\\ &\le \frac{1}{\lambda^{4\Bbbk-d+4}}\int \mathscr{L}^{\Bbbk}q\mathscr{L}^{\Bbbk}(\partial_s \xi_L)+\frac{1}{8}\int |\widetilde{\mathscr{L}}_{\lambda}v_{2\Bbbk-1}|^2+\frac{C}{\lambda^{4\Bbbk-d+4}}\left(\sqrt{\mathscr{E}_{2\Bbbk}}\|\mathscr{L}^{\Bbbk}\mathcal{F}_0\|_{L^2}+\|\mathscr{A}\mathscr{L}^{\Bbbk-1}\mathcal{F}_1\|_{L^2}^2\right).\label{in energy identity rhs 8 esti}
\end{align}\par 
Now substituting (\ref{commutator related estimate}) into (\ref{in energy identity rhs 5 esti}), then gathering estimates (\ref{in energy identity rhs 1 esti})-(\ref{in energy identity rhs 5 esti}) and (\ref{in energy identity rhs 6 esti})-(\ref{in energy identity rhs 8 esti}) into (\ref{energy identity}), we have 
\begin{align}
	\frac{1}{2}\frac{\mathrm{d}}{\mathrm{d}t}\left\{\frac{\mathscr{E}_{2\Bbbk}}{\lambda^{4\Bbbk-d+2}}(1+O(b_1))\right\}&\le -\frac{1}{8}\int |\widetilde{\mathscr{L}}_{\lambda}v_{2\Bbbk-1}|^2+\frac{1}{\lambda^{4\Bbbk-d+4}}\int \mathscr{L}^{\Bbbk}q\mathscr{L}^{\Bbbk}(\partial_s \xi_L)+\frac{Cb_1^2}{\lambda^{4\Bbbk-d+4}}\mathscr{E}_{2\Bbbk}\notag\\ &\quad\,\, +\frac{Cb_1}{\lambda^{4\Bbbk-d+4}}\sqrt{\mathscr{E}_{2\Bbbk}}\left[\left(\int \frac{|\mathscr{A}\mathscr{L}^{\Bbbk-1}\mathcal{F}|^2}{1+y^2}\right)^{\frac{1}{2}}+\left(\int \frac{|\mathscr{L}^{\Bbbk-1}\mathcal{F}|^2}{1+y^4}\right)^{\frac{1}{2}}\right]\notag\\ &\quad\,\,+\frac{C}{\lambda^{4\Bbbk-d+4}}\left(\sqrt{\mathscr{E}_{2\Bbbk}}\|\mathscr{L}^{\Bbbk}\mathcal{F}_0\|_{L^2}+\|\mathscr{A}\mathscr{L}^{\Bbbk-1}\mathcal{F}_1\|_{L^2}^2\right).\label{energy identity first sorted}
\end{align}
Note that in (\ref{energy identity first sorted}), integrals containing $\mathcal{F}$ can be controlled by correponding integrals containing $\widetilde{\Psi}_b,$ $\widehat{Mod},$ $\mathcal{H}(q),$ $\mathcal{N}(q).$ The integral containing $\mathcal{F}_0$ can be controlled by corresponding integral containing $\widetilde{\Psi}_b,$ $\widetilde{Mod}:=\widehat{Mod}+\partial_s \xi_L.$ The integral containing $\mathcal{F}_1$ can be controlled by corresponding integral containing $\mathcal{H}(q),$ $\mathcal{N}(q).$ Then we shall proceed to estimate terms in (\ref{energy identity first sorted}) with such further decompositions.\par 
Next we estimate $\widetilde{\Psi}_b$ term in (\ref{energy identity first sorted}). Applying (\ref{2 esti of tilde Psib}), we see that
\begin{equation}\label{in energy id esti of widetildePsib term}
	\left(\int \frac{|\mathscr{A}\mathscr{L}^{\Bbbk-1}\widetilde{\Psi}_b|^2}{1+y^2}\right)^{\frac{1}{2}}+\left(\int \frac{|\mathscr{L}^{\Bbbk-1}\widetilde{\Psi}_b|^2}{1+y^4}\right)^{\frac{1}{2}}+\|\mathscr{L}^{\Bbbk}\widetilde{\Psi}_b\|_{L^2}\lesssim b_1^{L+1+(1-\delta)(1+\eta)}.
\end{equation}\par 
Then we estimate $\widehat{Mod}$ term in (\ref{energy identity first sorted}), we claim
\begin{equation}\label{in energy id esti of widehatMod term}
	\left(\int \frac{|\mathscr{A}\mathscr{L}^{\Bbbk-1}\widehat{Mod}|^2}{1+y^2}\right)^{\frac{1}{2}}+\left(\int \frac{|\mathscr{L}^{\Bbbk-1}\widehat{Mod}|^2}{1+y^4}\right)^{\frac{1}{2}}\lesssim b_1^{(1-\delta)(1+\eta)}\left(\frac{\sqrt{\mathscr{E}_{2\Bbbk}}}{M^{2\delta}}+b_1^{L+1+(1-\delta)(1+\eta)}\right).
\end{equation}
It suffices to estimate the second term on the LHS, we shall omit the proof since it is a straight consequence of Proposition \ref{modulation estimates} and degrees of $T_i$ and $S_i.$\par 
Next we estimate $\widetilde{Mod}$ term in (\ref{energy identity first sorted}). Claim that
\begin{equation}\label{in energy id esti of widetildeMod term}
	\left(\int |\mathscr{L}^{\Bbbk}\widetilde{Mod}|^2\right)^{\frac{1}{2}}\lesssim b_1\left(\frac{\sqrt{\mathscr{E}_{2\Bbbk}}}{M^{2\delta}}+b_1^{\eta (1-\delta)}\sqrt{\mathscr{E}_{2\Bbbk}}+b_1^{L+1+(1-\delta)(1+\eta)-c_L\eta}\right).
\end{equation}
We further write 
\begin{align*}
	\widetilde{Mod}&=-\left(\frac{\lambda_s}{\lambda}+b_1\right)\Lambda \widetilde{Q}_b+\sum\limits_{i=1}^{L-1}[(b_i)_s+(2i-\gamma)b_1b_i-b_{i+1}]\widetilde{T}_i\\ &\quad\,\,+\sum\limits_{i=1}^{L}[(b_i)_s+(2i-\gamma)b_1b_i-b_{i+1}]\sum\limits_{j=i+1}^{L+2}\chi_{B_1}\frac{\partial S_j}{\partial b_i}\\ &\quad\,\,+\left[(b_L)_s+(2L-\gamma)b_1b_L+\partial_s\left\{\frac{\langle \mathscr{L}^L q,\chi_{B_0}\Lambda Q\rangle}{\langle \Lambda Q,\chi_{B_0}\Lambda Q\rangle}\right\}\right]\widetilde{T}_L+\frac{\langle \mathscr{L}^L q,\chi_{B_0}\Lambda Q\rangle}{\langle \Lambda Q,\chi_{B_0}\Lambda Q\rangle}\partial_s \widetilde{T}_L.
\end{align*}
Note that straight calculation yields 
\begin{align*}
	&\int |\mathscr{L}^{\Bbbk}\Lambda \widetilde{Q}_b|^2+\sum\limits_{i=1}^L\sum\limits_{j=i+1}^{L+2}\int \left|\mathscr{L}^{\Bbbk}\left(\chi_{B_1}\frac{\partial S_j}{\partial b_i}\right)\right|^2\lesssim b_1^2,\\ &\sum\limits_{i=1}^{L-1}\int |\mathscr{L}^{\Bbbk}\widetilde{T}_i|^2\lesssim b_1^{2(2-\delta)(1+\eta)},\,\,\,\int |\mathscr{L}^{\Bbbk}\widetilde{T}_L|^2\lesssim b_1^{2(1-\delta)(1+\eta)}. 
\end{align*}
And by the proof of (\ref{esti for the lhs in the pf of improved bd}), it follows that 
\begin{equation}\label{estimate of the fraction part}
	\left|\frac{\langle \mathscr{L}^L q,\chi_{B_0}\Lambda Q\rangle}{\langle \Lambda Q,\chi_{B_0}\Lambda Q\rangle}\right|\lesssim b_1^{-(1-\delta)}\sqrt{\mathscr{E}_{2\Bbbk}}.
\end{equation}
Also, in view of $|\partial_s \chi_{B_1}|=\left|\frac{(1+\eta)}{2}\frac{y}{B_1}\frac{\partial_s b_1}{b_1}\chi'\left(\frac{y}{B_1}\right)\right|\lesssim 1_{B_1\le y\le 2B_1}b_1,$ we have
\begin{equation*}
	\int |\mathscr{L}^{\Bbbk}(\partial_s \widetilde{T}_L)|^2\lesssim b_1^2\int_{B_1\le y\le 2B_1} \frac{y^{d-3}}{y^{4(\Bbbk-L)+2\gamma}}\,\mathrm{d}y\lesssim b_1^{2+2(1-\delta)(1+\eta)}.
\end{equation*}
Then above estimates combined with Proposition \ref{modulation estimates} and Proposition \ref{improved bound for bL prop} gives us (\ref{in energy id esti of widetildeMod term}).\par 
Then we estimate $\mathcal{H}(q)$ term in (\ref{energy identity first sorted}). Claim
\begin{equation}\label{in energy id esti of Hq term}
	\int \frac{|\mathscr{A}\mathscr{L}^{\Bbbk-1}\mathcal{H}(q)|^2}{1+y^2}+\int \frac{|\mathscr{L}^{\Bbbk-1}\mathcal{H}(q)|^2}{1+y^4}+\int |\mathscr{A}\mathscr{L}^{\Bbbk-1}\mathcal{H}(q)|^2\lesssim b_1^2\mathscr{E}_{2\Bbbk}.
\end{equation}
It suffices to estimate the third term on the LHS. Denoting \begin{equation*}
	\mathcal{H}(q)=:\phi q,\,\,\,\text{where}\,\,\,\phi:=\frac{-3(d-2)}{y^2}\widetilde{\Theta}_b[2(Q-1)+\widetilde{\Theta}_b]\,\,\,\text{and}\,\,\,\widetilde{\Theta}_b:=\chi_{B_1}\Theta_{b}.
\end{equation*}
Using (\ref{leibniz for A composite iterate of L}), we get
\begin{equation*}
	\mathscr{A}\mathscr{L}^{\Bbbk-1}\mathcal{H}(q)=\sum\limits_{m=0}^{\Bbbk-1}q_{2m+1}\phi_{2\Bbbk-1,2m+1}+q_{2m}\phi_{2\Bbbk-1,2m}.
\end{equation*}
Note that by direct computation, \begin{align*}
	|(Q-1)\widetilde{\Theta}_b|&=\left|\left(\sum\limits_{i=1}^L\chi_{B_1}b_iT_i+\sum\limits_{i=2}^{L+2}\chi_{B_1}S_i\right)(Q-1)\right|\\ &\lesssim 1_{y\le 2B_1}\left(\sum\limits_{i=1}^Lb_1^iy^{2i-\gamma}+\sum\limits_{i=2}^{L+2}b_1^iy^{2(i-1)-\gamma}\right)y^{\gamma}\ll 1_{y\le 2B_1}b_1y^{2-\gamma},
\end{align*}
where we also used the fact that
\begin{equation*}
	(b_1y^2)^Ny^{-\gamma}\lesssim b_1^{\frac{\gamma}{2}-\eta(N-\frac{\gamma}{2})}\ll 1,\,\,\,\text{for any integer}\,\,\,N\ge 1\,\,\,\text{and any}\,\,\,y\le 2B_1.
\end{equation*}
Then in general,\begin{equation*}
	|\phi_{k,i}|\lesssim 1_{y\le 2B_1}\frac{b_1}{1+y^{\gamma+k-i}}.
\end{equation*}
Hence combined with (\ref{coercivity used in energy estimate}), we derive
\begin{align*}
	\int |\mathscr{A}\mathscr{L}^{\Bbbk-1}\mathcal{H}(q)|^2&\lesssim \sum\limits_{m=0}^{\Bbbk-1}b_1^2\left(\int_{y\le 2B_1}\frac{|q_{2m+1}|^2}{1+y^{2\gamma+4(\Bbbk-m)-4}}+\int \frac{|q_{2m}|^2}{1+y^{2\gamma+4(\Bbbk-m)-2}}\right)\\ &\lesssim b_1^2\left(\sum\limits_{m=0}^{\Bbbk-2}\int \frac{|q_{2m+1}|^2}{y^6(1+y^{4(\Bbbk-2-m)})}+\int \frac{|q_{2\Bbbk-1}|^2}{y^2}+\sum\limits_{m=0}^{\Bbbk-1}\int \frac{|q_{2m}|^2}{y^4(1+y^{4(\Bbbk-1-m)})}\right)\\ &\lesssim b_1^2\mathscr{E}_{2\Bbbk}(q). 
\end{align*}
This concludes the proof of (\ref{in energy id esti of Hq term}).\par
Next we estimate $\mathcal{N}(q)$ term in (\ref{energy identity first sorted}). Claim
\begin{equation}\label{in energy id esti of Nq term 1}
	\int |\mathscr{A}\mathscr{L}^{\Bbbk-1}\mathcal{N}(q)|^2\lesssim b_1^{2L+1+2(1-\delta)(1+\eta)}
\end{equation}
and \begin{equation}\label{in energy id esti of Nq term 2}
	\int \frac{|\mathscr{A}\mathscr{L}^{\Bbbk-1}\mathcal{N}(q)|^2}{1+y^2}+\int \frac{|\mathscr{L}^{\Bbbk-1}\mathcal{N}(q)|^2}{1+y^4}\lesssim b_1^{2L+2+2(1-\delta)(1+\eta)}.
\end{equation}
We shall only prove (\ref{in energy id esti of Nq term 1}) since the proof of (\ref{in energy id esti of Nq term 2}) is similar. Let us estimate the integral for $y\le 1$ and $y\ge 1$ separately.\par 
When $y< 1.$ Rewriting \begin{equation*}
	\mathcal{N}(q)=\frac{q^2}{y^2}\phi,\,\,\,\text{where}\,\,\,\phi:=\phi_1+\phi_2,\,\,\,\phi_1:=(d-2)(3\widetilde{Q}_b+q),\,\,\,\phi_2:=-3(d-2).
\end{equation*}
By \textnormal{(\romannumeral1)} of Lemma \ref{coercivity-determined esti on q}, we get 
\begin{equation*}
	\frac{q^2}{y^2}=\frac{1}{y^2}\left(\sum\limits_{i=0}^{\Bbbk-1}c_iT_i(y)+r_q(y)\right)^2=\sum\limits_{i=0}^{\Bbbk-1}\widetilde{c}_iy^{4i+2}+\widetilde{r}_q
\end{equation*}
with\begin{equation*}
	|\widetilde{c}_i|\lesssim \mathscr{E}_{2\Bbbk},\,\,\,|\partial_y^j \widetilde{r}_q|\lesssim y^{2\Bbbk+1-\frac{d}{2}-j}|\ln y|^{\Bbbk}\mathscr{E}_{2\Bbbk}.
\end{equation*}
By (\ref{asymp: ground state}), Proposition \ref{first approximation}, \textnormal{(\romannumeral1)} of Lemma \ref{coercivity-determined esti on q} and \textnormal{(\romannumeral3)} of Definition \ref{bootstrap assump}, we have
\begin{equation*}
	\phi_1=\sum\limits_{i=0}^{\Bbbk-1}\widehat{c}_iy^{2i+2}+\widehat{r}_q
\end{equation*} 
with\begin{equation*}
	|\widehat{c}_i|\lesssim 1,\,\,\,|\partial_y^j \widehat{r}_q|\lesssim y^{2\Bbbk+1-\frac{d}{2}-j}|\ln y|^{\Bbbk}.
\end{equation*}
It immediately follows \begin{equation*}
	\mathcal{N}(q)=\sum\limits_{i=0}^{\Bbbk-1}\widehat{\widetilde{c}}_iy^{2i+2}+\widehat{\widetilde{r}}_q
\end{equation*}
with\begin{equation*}
	|\widehat{\widetilde{c}}_i|\lesssim \mathscr{E}_{2\Bbbk}\,\,\,\text{and}\,\,\,|\partial_y^j \widehat{\widetilde{r}}_q|\lesssim y^{2\Bbbk+1-\frac{d}{2}-j}|\ln y|^{\Bbbk}\mathscr{E}_{2\Bbbk}.
\end{equation*}
Hence,\begin{align*}
	|\mathscr{A}\mathscr{L}^{\Bbbk-1}\mathcal{N}(q)|&=|\mathscr{A}\mathscr{L}^{\Bbbk-1}(\sum\limits_{i=0}^{\Bbbk-1}\widehat{\widetilde{c}}_iy^{2i+2})+\mathscr{A}\mathscr{L}^{\Bbbk-1}\widehat{\widetilde{r}}_q|\\ &\lesssim \sum\limits_{i=0}^{\Bbbk-1}|\widehat{\widetilde{c}}_i|y^3+\sum\limits_{i=0}^{2\Bbbk-1}\frac{|\partial_y^j \widehat{\widetilde{r}}_q|}{y^{2\Bbbk-1-i}} \lesssim y^{-\frac{d}{2}+2}|\ln y|^{\Bbbk}\mathscr{E}_{2\Bbbk}.
\end{align*}
Then by \textnormal{(\romannumeral3)} of Definition \ref{bootstrap assump}, we estimate \begin{equation}\label{esti of ALBbbk-1Nq2 integral for y less than 1}
	\int_{y<1} |\mathscr{A}\mathscr{L}^{\Bbbk-1}\mathcal{N}(q)|^2\lesssim \mathscr{E}_{2\Bbbk}^2 \int_{y<1} y|\ln y|^{2\Bbbk}\,\mathrm{d}y\lesssim \mathscr{E}_{2\Bbbk}^2\lesssim b_1^{4[L+(1-\delta)(1+\eta)]}.
\end{equation}\par
When $y\ge 1,$ rewriting 
\begin{equation*}
	\mathcal{N}(q)=Z^2\phi,\,\,\,\text{where}\,\,\,Z:=\frac{q}{y},\,\,\,\phi:=(d-2)[3(\widetilde{Q}_b-1)+q].
\end{equation*}
By Leibniz rule, we get \begin{align*}
	\int_{y\ge 1} |\mathscr{A}\mathscr{L}^{\Bbbk-1}\mathcal{N}(q)|^2&\lesssim \sum\limits_{k=0}^{2\Bbbk-1}\int_{y\ge 1}\frac{|\partial_y^k\mathcal{N}(q)|^2}{y^{4\Bbbk-2k-2}}\\ &\lesssim \sum\limits_{k=0}^{2\Bbbk-1}\sum\limits_{i=0}^k\int_{y\ge 1}\frac{|\partial_y^iZ^2|^2|\partial_y^{k-i}\phi|^2}{y^{4\Bbbk-2k-2}}\\ &\lesssim \sum\limits_{k=0}^{2\Bbbk-1}\sum\limits_{i=0}^k\sum\limits_{m=0}^i\int_{y\ge 1}\frac{|\partial_y^mZ|^2|\partial_y^{i-m}Z|^2|\partial_y^{k-i}\phi|^2}{y^{4\Bbbk-2k-2}}. 
\end{align*}
Then we focus on proving that for $0\le k\le 2\Bbbk-1,$ $0\le i\le k,$ $0\le m\le i,$ \begin{equation}\label{most thought esti in energy estimate}
	A_{k,i,m}:=\int_{y\ge 1}\frac{|\partial_y^mZ|^2|\partial_y^{i-m}Z|^2|\partial_y^{k-i}\phi|^2}{y^{4\Bbbk-2k-2}}\lesssim b_1^{2L+1+2(1-\delta)(1+\eta)},
\end{equation}
which would conclude the proof of (\ref{in energy id esti of Nq term 1}). We shall split the proof into three cases as the following three paragraphs.\par 
When $k=0.$ In this case, $k=i=m=0.$ Note that $\phi$ is bounded as $y\rightarrow \infty,$ then 
\begin{equation*}
	A_{0,0,0}=\int_{y\ge 1} \frac{|q|^4|\phi|^2}{y^{4\Bbbk+2}}y^{d-3}\,\mathrm{d}y\lesssim \int_{1\le y\le B_0} \frac{|q|^4}{y^{4\Bbbk+5-d}}\,\mathrm{d}y+\int_{y>B_0} \frac{|q|^4}{y^{4\Bbbk+5-d}}\,\mathrm{d}y. 
\end{equation*}
By \textnormal{(\romannumeral4)} of Lemma \ref{coercivity-determined esti on q}, Definition \ref{bootstrap assump} and recall that $d=4\hbar+4\delta+2\gamma+2,$ we have
\begin{align}
	\int_{1\le y\le B_0} \frac{|q|^4}{y^{4\Bbbk+5-d}}\,\mathrm{d}y&\lesssim \left\|\frac{y^{d-4}|q|^2}{y^{2(2\Bbbk-1)}}\right\|_{L^{\infty}(y\ge 1)}\left\|\frac{y^{d-4}|q|^2}{y^{2(2l+2\hbar+3)}}\right\|_{L^{\infty}(y\ge 1)}\int_{1\le y\le B_0} y^{4l-4\delta-2\gamma+5}\,\mathrm{d}y\notag\\ &\lesssim \mathscr{E}_{2\Bbbk}\mathscr{E}_{2(l+\hbar+2)}B_0^{4l-4\delta-2\gamma+6}\lesssim b_1^{(1+\gamma-K\eta)+2L+2(1-\delta)(1+\eta)}.\label{esti of A000 1}
\end{align}
Similarly, 
\begin{align}
	\int_{y>B_0}\frac{|q|^4}{y^{4\Bbbk+5-d}}\,\mathrm{d}y&\lesssim \left\|\frac{y^{d-4}|q|^2}{y^{2(2\Bbbk-2l-1)}}\right\|_{L^{\infty}(y\ge1)}\left\|\frac{y^{d-4}|q|^2}{y^{2(2l+2\hbar+1)}}\right\|_{L^{\infty}(y\ge1)}\int_{y>B_0}y^{-4\delta-2\gamma+1}\,\mathrm{d}y\notag\\ &\lesssim \mathscr{E}_{2(\Bbbk-l)}\mathscr{E}_{2(l+\hbar+1)}B_0^{-4\delta-2\gamma+2}\lesssim b_1^{2L+2(1-\delta)(1+\eta)+(1+\gamma)-2(K+1-\delta)\eta}.\label{esti of A000 2}
\end{align}\par 
When $k\ge 1$ and $i=k.$ By Leibniz rule, 
\begin{equation*}
	|\partial_y^n Z|^2\lesssim \sum\limits_{j=0}^n\frac{|\partial_y^j q|^2}{y^{2+2n-2j}},\,\,\,\text{for all}\,\,\,n\in \mathbb{N}.\Longrightarrow A_{k,k,m}\lesssim \sum\limits_{j=0}^m\sum\limits_{l=0}^{k-m}\int_{y\ge 1}\frac{|\partial_y^j q|^2|\partial_y^l q|^2}{y^{4\Bbbk-2j-2l+2}}. 
\end{equation*}
Direct computation implies
\begin{align}
	B_{j,l}:&=\int_{y\ge 1}\frac{|\partial_y^j q|^2|\partial_y^l q|^2}{y^{4\Bbbk-2j-2l+2}}y^{d-3}\,\mathrm{d}y\notag \\&=\int_{1\le y\le B_0}\frac{(y^{d-4}|\partial_y^j q|^2)(y^{d-4}|\partial_y^l q|^2)}{y^{4\Bbbk-2j-2l+4\hbar+6}}y^{7-4\delta-2\gamma}\,\mathrm{d}y+\int_{y>B_0}\frac{(y^{d-4}|\partial_y^j q|^2)(y^{d-4}|\partial_y^l q|^2)}{y^{4\Bbbk-2j-2l+4\hbar}}y^{-(4\delta+2\gamma-1)}\,\mathrm{d}y\notag \\ &\lesssim \left\|\frac{(y^{d-4}|\partial_y^j q|^2)(y^{d-4}|\partial_y^l q|^2)}{y^{4\Bbbk-2j-2l+4\hbar+6}}\right\|_{L^{\infty}(y\ge 1)}b_1^{2\delta+\gamma-4}+\left\|\frac{(y^{d-4}|\partial_y^j q|^2)(y^{d-4}|\partial_y^l q|^2)}{y^{4\Bbbk-2j-2l+4\hbar}}\right\|_{L^{\infty}(y\ge 1)}b_1^{2\delta+\gamma-1}\notag \\ &=\left\|\frac{(y^{d-4}|\partial_y^j q|^2)(y^{d-4}|\partial_y^l q|^2)}{y^{2J_1+2J_2-2j-2l}}\right\|_{L^{\infty}(y\ge 1)}b_1^{2\delta+\gamma-4}+\left\|\frac{(y^{d-4}|\partial_y^j q|^2)(y^{d-4}|\partial_y^l q|^2)}{y^{2J_3+2J_4-2j-2l}}\right\|_{L^{\infty}(y\ge 1)}b_1^{2\delta+\gamma-1}\notag \\ &=:B_{j,l,J_1,J_2}b_1^{2\delta+\gamma-4}+B_{j,l,J_3,J_4}b_1^{2\delta+\gamma-1}, \label{revised combine}
\end{align}
where $J_1+J_2=2\Bbbk+2\hbar+3,$ $J_3+J_4=2\Bbbk+2\hbar,$ and those $J$ are to be determined. For the first right hand side term in (\ref{revised combine}), when $l\ge 3,$ we choose $J_1=2\Bbbk-2l+3,$ $J_2=2\hbar+2l.$ Then by \textnormal{(\romannumeral4)} of Lemma \ref{coercivity-determined esti on q} and Definition \ref{bootstrap assump}, we see that
\begin{align*}
	B_{j,l,J_1,J_2}&\lesssim \left\|\frac{y^{d-4}|\partial_y^j q|^2}{y^{2J_1-2j}}\right\|_{L^{\infty}(y\ge 1)}\left\|\frac{y^{d-4}|\partial_y^l q|^2}{y^{2J_2-2l}}\right\|_{L^{\infty}(y\ge 1)}\\ &\lesssim \mathscr{E}_{J_1+1}\sqrt{\mathscr{E}_{J_2}\mathscr{E}_{J_2+2}} \\ &\lesssim b_1^{2L-l+3(1-\delta)+4-\frac{3}{2}K\eta+\frac{l}{2l-\gamma}(2l-2\delta-\gamma)}\lesssim b_1^{2L+4(1-\delta)+3-\frac{\gamma}{2}-\frac{3}{2}K\eta},
\end{align*}
where in the last inequality we used the fact that $\frac{l}{2l-\gamma}>\frac{1}{2}.$ When $l=1$ or $2,$ we choose $J_1=2\Bbbk-1,$ $J_2=2\hbar+4,$ then similarly one gets for $l=2,$ 
\begin{equation*}
		B_{j,l,J_1,J_2}\lesssim \mathscr{E}_{J_1+1}\sqrt{\mathscr{E}_{J_2}\mathscr{E}_{J_2+2}}
	\lesssim b_1^{2L+2(1-\delta)(1+\eta)+5-2\delta-\frac{\gamma}{2}-\frac{1}{2}K\eta},
\end{equation*} 
and for $l=1,$ there holds 
\begin{equation*}
		B_{j,l,J_1,J_2}\lesssim \mathscr{E}_{J_1+1}\sqrt{\mathscr{E}_{J_2}\mathscr{E}_{J_2+2}}
	\lesssim b_1^{2L+2(1-\delta)(1+\eta)+5-2\delta-K\eta}.
\end{equation*}
Hence \begin{equation}\label{revised ingre 1}
	B_{j,l,J_1,J_2}b_1^{2\delta+\gamma-4}\lesssim \left\{\begin{aligned}
		&b_1^{2L+1+2(1-\delta)(1+\eta)+\frac{\gamma}{2}-\left[\frac{3}{2}K+2(1-\delta)\right]\eta},\,\,\,\text{for}\,\,\,l\ge 3,\\
		&b_1^{2L+1+2(1-\delta)(1+\eta)+\frac{\gamma}{2}-\frac{1}{2}K\eta},\,\,\,\text{for}\,\,\,l=2,\\
		&b_1^{2L+1+2(1-\delta)(1+\eta)+\gamma-K\eta},\,\,\,\text{for}\,\,\,l=1.
	\end{aligned}\right. 
\end{equation}
For the second right side hand term in (\ref{revised combine}), we similarly have
\begin{equation*}
	B_{j,l,J_3,J_4}\lesssim \sqrt{\mathscr{E}_{J_3}\mathscr{E}_{J_3+2}\mathscr{E}_{J_4}\mathscr{E}_{J_4+2}}\lesssim b_1^{2L+4(1-\delta)-\frac{\gamma}{2}-\frac{3}{2}K\eta},
\end{equation*}
where we choose $J_3=2\Bbbk-2l,$ $J_4=2\hbar+2l.$ Thus
\begin{equation}\label{revised ingre 2}
	B_{j,l,J_3,J_4}b_1^{2\delta+\gamma-1}\lesssim b_1^{2L+1+2(1-\delta)(1+\eta)+\frac{\gamma}{2}-\left[\frac{3}{2}K+2(1-\delta)\right]\eta}.
\end{equation} 
Therefore, injecting the estimates (\ref{revised ingre 1}) and (\ref{revised ingre 2}) into (\ref{revised combine}), we derive the estimate of $B_{j,l}$ as   
\begin{equation}\label{esti of B_{j,l}}
	B_{j,l}\lesssim 
	b_1^{2L+1+2(1-\delta)(1+\eta)+\frac{\gamma}{2}-\left[\frac{3}{2}K+2(1-\delta)\right]\eta},
\end{equation}
which holds for any positive integer $l.$\par 
When $k\ge 1$ and $i\le k-1.$ Again by Leibniz rule, we further write
\begin{equation*}
	A_{k,i,m}\lesssim \sum\limits_{j=0}^m\sum\limits_{l=0}^{i-m}\int_{y\ge 1}\frac{|\partial_y^j q|^2|\partial_y^l q|^2||\partial_y^{k-i} \phi|^2}{y^{4\Bbbk-2j-2l+2-2(k-i)}}.
\end{equation*}
We shall need pointwise estimates of $\partial_y^n \phi,$ for $n\in \mathbb{N}_{+}.$ Note that by degrees of $T_k$ and $S_k$, we derive
\begin{align*}
	|\partial_y^n \widetilde{Q}_b|&=\left|\partial_y^n\left(Q+\sum\limits_{k=1}^L\chi_{B_1}b_kT_k+\sum\limits_{k=2}^{L+2}\chi_{B_1}S_k\right)\right|\\ &\lesssim \frac{1}{y^{\gamma+n}}+\sum\limits_{k=1}^L\frac{b_1^ky^{2k}}{y^{\gamma+n}}1_{y\le 2B_1}\lesssim \frac{b_1^{-\eta(L+1)}}{y^{\gamma+n}}.
\end{align*}
By \textnormal{(\romannumeral4)} of Lemma \ref{coercivity-determined esti on q} and Definition \ref{bootstrap assump}, we have
for $1\le y\le B_0,$ \begin{equation*}
	|\partial_y^n q|^2\lesssim y^{2(2\Bbbk-1-n)}\left|\frac{\partial_y^n q}{y^{2\Bbbk-1-n}}\right|^2\lesssim y^{2(2\Bbbk-1-n)-(d-4)}\mathscr{E}_{2\Bbbk}\lesssim b_1^{\eta+\gamma+2(1-\delta)\eta},
\end{equation*}
for $y\ge B_0,$ \begin{align*}
	|\partial_y^n q|^2&\lesssim y^{2(2\hbar+2l+1-n)}\left|\frac{\partial_y^n q}{y^{2\hbar+2l+1-n}}\right|^2\\&\lesssim y^{2(2\hbar+2l+1-n)-(d-4)}\mathscr{E}_{2\hbar+2l+2}\lesssim y^{4l+4(1-\delta)}b_1^{n+\gamma+2l+2(1-\delta)-K\eta}.
\end{align*}
Thus \begin{equation*}
	|\partial_y^n \phi|^2\lesssim |\partial_y^n \widetilde{Q}_b|^2+|\partial_y^n q|^2\lesssim \left\{\begin{aligned}
		&\frac{b_1^{-2(L+1)\eta}}{y^{2\gamma+2n}},\,\,\,\text{when}\,\,\,1\le y\le B_0.\\ &b_1^{-C_{L,K}\eta+\gamma+n}(b_1y^2)^{2l+2(1-\delta)},\,\,\,\text{when}\,\,\,y\ge B_0.
	\end{aligned}\right.
\end{equation*}
Then similar to the proof of (\ref{esti of B_{j,l}}), we compute
\begin{align}
	A_{k,i,m}&\lesssim b_1^{-C_{L,K}\eta}\sum_{j=0}^m\sum_{l=0}^{i-m}\Bigg(\int_{1\le y\le B_0}\frac{|\partial_y^j q|^2|\partial_y^l q|^2}{y^{4\Bbbk-2j-2l+2+2\gamma}}y^{d-3}\,\mathrm{d}y\notag\\&\quad\,\,+b_1^{\gamma+\alpha}\int_{y>B_0}\frac{|\partial_y^j q|^2|\partial_y^l q|^2}{y^{4\Bbbk-2j-2l+2-2\alpha}}y^{d-3}\,\mathrm{d}y\Bigg)\notag\\&\lesssim b_1^{2L+1+2(1-\delta)(1+\eta)+\frac{\gamma}{2}-C_{K,L,\delta}\eta},\label{estimate of Akim}
\end{align}
where $\alpha:=k-i+2l+2(1-\delta).$\par 
In view of (\ref{esti of A000 1})-(\ref{estimate of Akim}), we conclude the proof of (\ref{most thought esti in energy estimate}), then by (\ref{esti of ALBbbk-1Nq2 integral for y less than 1}) and (\ref{most thought esti in energy estimate}), we complete the proof of (\ref{in energy id esti of Nq term 1}).\par
It remains to estimate the integral $\frac{1}{\lambda^{4\Bbbk-d+4}}\int \mathscr{L}^{\Bbbk}q\mathscr{L}^{\Bbbk}(\partial_s \xi_L)$ in (\ref{energy identity first sorted}). Let us further write
\begin{align}
	\frac{1}{\lambda^{4\Bbbk-d+4}}\int \mathscr{L}^{\Bbbk}q\mathscr{L}^{\Bbbk}(\partial_s \xi_L)&=\frac{\mathrm{d}}{\mathrm{d}s}\frac{1}{\lambda^{4\Bbbk-d+4}}\left(\int \mathscr{L}^{\Bbbk}q\mathscr{L}^{\Bbbk}\xi_L-\frac{1}{2}\int |\mathscr{L}^{\Bbbk}\xi_L|^2\right)\notag\\ &\quad\,\,+\frac{(4\Bbbk-d+4)}{\lambda^{4\Bbbk-d+4}}\frac{\lambda_s}{\lambda}\left(\int \mathscr{L}^{\Bbbk}q\mathscr{L}^{\Bbbk}\xi_L-\frac{1}{2}\int |\mathscr{L}^{\Bbbk}\xi_L|^2\right)\notag\\ &\quad\,\,-\frac{1}{\lambda^{4\Bbbk-d+4}}\int \mathscr{L}^{\Bbbk}(\partial_s q-\partial_s \xi_L) \mathscr{L}^{\Bbbk}\xi_L.\label{expreesion of the oscilation integral}
\end{align}
Recall the proof of (\ref{in energy id esti of widetildeMod term}), we deduce that \begin{equation}\label{esti of L2 norm of LBbbkxiL}
	\int |\mathscr{L}^{\Bbbk}\xi_L|^2\lesssim b_1^{2(1-\delta)\eta}\mathscr{E}_{2\Bbbk}.
\end{equation}
Then by H\"older, we get \begin{equation*}
	\left|\int \mathscr{L}^{\Bbbk}q\mathscr{L}^{\Bbbk}\xi_L\right|\lesssim b_1^{(1-\delta)\eta}\mathscr{E}_{2\Bbbk}.
\end{equation*}
Thus \begin{equation}\label{first line on RHS of the expreesion of oscilation integral}
	\frac{\mathrm{d}}{\mathrm{d}s}\frac{1}{\lambda^{4\Bbbk-d+4}}\left(\int \mathscr{L}^{\Bbbk}q\mathscr{L}^{\Bbbk}\xi_L-\frac{1}{2}\int |\mathscr{L}^{\Bbbk}\xi_L|^2\right)=\frac{\mathrm{d}}{\mathrm{d}t}\left\{\frac{\mathscr{E}_{2\Bbbk}}{\lambda^{4\Bbbk-d+2}}O(b_1^{(1-\delta)\eta})\right\},
\end{equation}
and \begin{equation}\label{second line on RHS of the expreesion of oscilation integral}
\left|\frac{(4\Bbbk-d+4)}{\lambda^{4\Bbbk-d+4}}\frac{\lambda_s}{\lambda}\left(\int \mathscr{L}^{\Bbbk}q\mathscr{L}^{\Bbbk}\xi_L-\frac{1}{2}\int |\mathscr{L}^{\Bbbk}\xi_L|^2\right)\right|\lesssim \frac{b_1^{1+(1-\delta)\eta}\mathscr{E}_{2\Bbbk}}{\lambda^{4\Bbbk-d+4}}.
\end{equation}\par 
On the third line on the RHS of (\ref{expreesion of the oscilation integral}). By (\ref{eq: eq of q}), we further write
\begin{align}
	&\quad\,\,\int \mathscr{L}^{\Bbbk}(\partial_s q-\partial_s \xi_L) \mathscr{L}^{\Bbbk}\xi_L\notag\\ &=\frac{\lambda_s}{\lambda}\int \Lambda q\mathscr{L}^{2\Bbbk}\xi_L-\int \mathscr{L}^{\Bbbk}q\mathscr{L}^{\Bbbk+1}\xi_L+\int \mathscr{L}^{\Bbbk}(-\widetilde{\Psi}_b-\widetilde{Mod}+\mathcal{H}(q)-\mathcal{N}(q))\mathscr{L}^{\Bbbk}\xi_L.\label{expression of the third line on RHS of oscilation integral}
\end{align}
By H\"older, (\ref{coercivity used in modulation estimates}) and (\ref{estimate of the fraction part}), we have 
\begin{align}
	&\quad\,\,\left|\frac{\lambda_s}{\lambda}\int \Lambda q\mathscr{L}^{2\Bbbk}\xi_L\right|\notag\\&\lesssim b_1\left(\frac{|\partial_y q|^2}{1+y^{4\Bbbk-2}}\right)^{\frac{1}{2}}\left(\int y^2(1+y^{4\Bbbk-2})\left|\mathscr{L}^{2\Bbbk}\left(\frac{\langle \mathscr{L}^L q,\chi_{B_0}\Lambda Q\rangle}{\langle \chi_{B_0}\Lambda Q,\Lambda Q\rangle}(1-\chi_{B_1})T_L\right)\right|^2\right)^{\frac{1}{2}} \lesssim b_1^{1+\eta(1-\delta)}\mathscr{E}_{2\Bbbk}.\label{1 esti of expression of the third line on RHS of oscilation integral}
\end{align}
Again by (\ref{estimate of the fraction part}), 
\begin{equation*}
	\int |\mathscr{L}^{\Bbbk+1}\xi_L|^2\lesssim \left|\frac{\langle \mathscr{L}^L q,\chi_{B_0}\Lambda Q\rangle}{\langle \chi_{B_0}\Lambda Q,\Lambda Q\rangle}\right|^2\int \left|\mathscr{L}^{\Bbbk+1}\left((1-\chi_{B_1})T_L\right)\right|^2\lesssim b_1^{2+2(2-\delta)\eta}\mathscr{E}_{2\Bbbk},
\end{equation*}
combined with H\"older, we see that 
\begin{equation}\label{2 esti of expression of the third line on RHS of oscilation integral}
	\left|\int \mathscr{L}^{\Bbbk}q\mathscr{L}^{\Bbbk+1}\xi_L\right|\lesssim b_1^{1+(2-\delta)\eta}\mathscr{E}_{2\Bbbk}.
\end{equation}
By H\"older, (\ref{in energy id esti of widetildePsib term}), (\ref{in energy id esti of widetildeMod term}) and (\ref{esti of L2 norm of LBbbkxiL}), we get \begin{align}
	\left|\int \mathscr{L}^{\Bbbk}(-\widetilde{\Psi}_b-\widetilde{Mod})\mathscr{L}^{\Bbbk}\xi_L \right|&\lesssim \left(\int |\mathscr{L}^{\Bbbk}(\widetilde{\Psi}_b+\widetilde{Mod})|^2\right)^{\frac{1}{2}}\left(\int |\mathscr{L}^{\Bbbk}\xi_L|^2\right)^{\frac{1}{2}}\notag\\ &\lesssim b_1^{\eta(1-\delta)+L+1+(1-\delta)(1+\eta)}\sqrt{\mathscr{E}_{2\Bbbk}}+b_1^{1+(1-\delta)\eta}\mathscr{E}_{2\Bbbk}.\label{3 esti of expression of the third line on RHS of oscilation integral}
\end{align}
Similarly, by H\"older, (\ref{in energy id esti of Hq term}), (\ref{in energy id esti of Nq term 2}) and (\ref{esti of L2 norm of LBbbkxiL}), we obtain
\begin{align}
	\left|\int \mathscr{L}^{\Bbbk}(\mathcal{H}(q)-\mathcal{N}(q))\mathscr{L}^{\Bbbk}\xi_L\right|&\lesssim \left(\int \frac{|\mathscr{L}^{\Bbbk-1}(\mathcal{H}(q)-\mathcal{N}(q))|^2}{1+y^4}\right)^{\frac{1}{2}}\left(\int (1+y^4)|\mathscr{L}^{\Bbbk+1}\xi_L|^2\right)^{\frac{1}{2}}\notag\\ &\lesssim b_1^{1+\eta(1-\delta)}\mathscr{E}_{2\Bbbk}+b_1^{\eta(1-\delta)+L+1+(1-\delta)(1+\eta)}\sqrt{\mathscr{E}_{2\Bbbk}}.\label{4 esti of expression of the third line on RHS of oscilation integral}
\end{align}
Substituting (\ref{1 esti of expression of the third line on RHS of oscilation integral})-(\ref{4 esti of expression of the third line on RHS of oscilation integral}) into (\ref{expression of the third line on RHS of oscilation integral}), and then taking (\ref{first line on RHS of the expreesion of oscilation integral})-(\ref{expression of the third line on RHS of oscilation integral}) into (\ref{expreesion of the oscilation integral}), we derive
\begin{align}
	\frac{1}{\lambda^{4\Bbbk-d+4}}\int \mathscr{L}^{\Bbbk}q\mathscr{L}^{\Bbbk}(\partial_s\xi_L)&=\frac{\mathrm{d}}{\mathrm{d}t}\left\{\frac{\mathscr{E}_{2\Bbbk}}{\lambda^{4\Bbbk-d+2}}O(b_1^{\eta(1-\delta)})\right\}\notag\\ &\quad\,\,+\frac{b_1}{\lambda^{4\Bbbk-d+4}}O(b_1^{\eta(1-\delta)}\mathscr{E}_{2\Bbbk}+b_1^{\eta(1-\delta)}b_1^{L+(1-\delta)(1+\eta)}\sqrt{\mathscr{E}_{2\Bbbk}}).\label{estimate of oscilation integral}
\end{align}\par 
Now substituting (\ref{boootstrap assump on E2Bbbk}) in the first term in the second line of (\ref{estimate of oscilation integral}), then collecting estimates (\ref{in energy id esti of widetildePsib term})-(\ref{in energy id esti of widetildeMod term}), (\ref{in energy id esti of Hq term})-(\ref{in energy id esti of Nq term 2}) and (\ref{estimate of oscilation integral}) into (\ref{energy identity first sorted}), we get the desired estimate (\ref{monotonicity for E2Bbbk}), this concludes the proof.
\end{pf}

\vspace{\baselineskip}
\section{Improved bootstrap, transverse crossing and conclusion}
Next we prove improved bootstrap estimates as follows.
\begin{prop}\label{improved bootstrap}
	Given initial data as in Definition \ref{def of initial data}, assuming for some large universal constant $K$ there is $s_0(K)\gg 1$ such that $(b(s),q(s))\in \mathcal{S}_K(s)$ on $s\in [s_0,s_1]$ for some $s_1\ge s_0.$ Then for all $s\in [s_0,s_1],$ we have
\begin{align}
|\mathcal{V}_1(s)|&\le s^{-\frac{\eta}{2}(1-\delta)},\label{improved esti for V1}\\
|b_k(s)|&\lesssim s^{-(k+\eta(1-\delta))},\,\,\,\text{for}\,\,\,l+1\le k\le L,\label{improved esti for bk from l+1 to L}\\
\mathscr{E}_{2m}&\le \left\{\begin{aligned}
&\frac{1}{2}Ks^{-\frac{l(4m-d+2)}{2l-\gamma}},\,\,\, \hbar+2\le m\le l+\hbar,\\
&\frac{1}{2}Ks^{-[2(m-\hbar-1)+2(1-\delta)-K\eta]},\,\,\, l+\hbar+1\le m\le \Bbbk-1,
\end{aligned}\right. \label{improved esti for E2m lower}\\
\mathscr{E}_{2\Bbbk}&\le \frac{1}{2}Ks^{-[2L+2(1-\delta)(1+\eta)]}.\label{improved esti for E2Bbbk}
\end{align}
\end{prop}
\begin{pf}
	\normalfont
	In order to make use of Proposition \ref{energy estimates}, let us firstly estimate $\lambda$ in the variable $s.$ By Lemma \ref{linearization of bk from 1 to l} and Definition \ref{bootstrap assump}, we have
	\begin{equation*}
		b_1(s)=\frac{l}{2l-\gamma}\frac{1}{s}+\frac{\mathcal{U}_1}{s}=\frac{l}{2l-\gamma}\frac{1}{s}+O\left(\frac{1}{s^{1+\frac{\eta(1-\delta)}{2}}}\right),
	\end{equation*} 
combined with (\ref{modu esti for lambda and bk from 1 to L-1}), it follows that
\begin{equation*}
	-\frac{\lambda_s}{\lambda}=b_1+O(b_1^{L+1+(1-\delta)(1+\eta)})=\frac{l}{2l-\gamma}\frac{1}{s}+O\left(\frac{1}{s^{1+\frac{\eta(1-\delta)}{2}}}\right).
\end{equation*}
Or say \begin{equation}\label{dynamic estimate of lambda}
	\left|\partial_s\ln\left(s^{\frac{l}{2l-\gamma}}\lambda\right)\right|\lesssim \frac{1}{s^{1+\frac{\eta(1-\delta)}{2}}}. 
\end{equation}
Integrating (\ref{dynamic estimate of lambda}) from $s_0$ to $s$ yields
\begin{equation*}
	e^{-s_0^{-\frac{\eta(1-\delta)}{2}}}\lesssim \frac{\lambda(s)s^{\frac{l}{2l-\gamma}}}{s_0^{\frac{l}{2l-\gamma}}}\lesssim e^{s_0^{-\frac{\eta(1-\delta)}{2}}},
\end{equation*}
hence\begin{equation}\label{estimate of lambda in s}
	\lambda(s)\simeq \left(\frac{s_0}{s}\right)^{\frac{l}{2l-\gamma}}.
\end{equation}\par 
Improved bound for $\mathscr{E}_{2\Bbbk}:$ Integrating (\ref{monotonicity for E2Bbbk}) from $s_0$ to $s,$ we get
\begin{align*}
	\frac{\mathscr{E}_{2\Bbbk}(s)}{\lambda(s)^{4\Bbbk-d+2}}\left(1+O(b_1(s)^{\eta(1-\delta)})\right)&\lesssim \mathscr{E}_{2\Bbbk}(s_0)\left(1+O(b_1(s_0)^{\eta(1-\delta)})\right)\\ &\quad\,\,+\int_{s_0}^s \frac{b_1(\tau)}{\lambda(\tau)^{4\Bbbk-d+4}}\Bigg(b_1(\tau)^{L+(1-\delta)(1+\eta)}\sqrt{\mathscr{E}_{2\Bbbk}(\tau)}\\&\quad\,\,+\frac{\mathscr{E}_{2\Bbbk}(\tau)}{M^{2\delta}}+b_1(\tau)^{2L+2(1-\delta)(1+\eta)}\Bigg)\,\mathrm{d}\tau,
\end{align*}
where we used the assumption $\lambda(s_0)=1.$ Then applying (\ref{estimate of lambda in s}), $b_1(s)\simeq \frac{c_1}{s},$ Definition \ref{def of initial data} and Definition \ref{bootstrap assump}, and for convenience discarding $\frac{1}{\lambda^2}$ in the integral (note that $\frac{1}{\lambda(\tau)^2}\lesssim \left(\frac{s_1}{s_0}\right)^{\frac{l}{2l-\gamma}}\lesssim 1$), we get
\begin{align*}
	\mathscr{E}_{2\Bbbk}(s)&\lesssim s_0^{-\frac{l}{2l-\gamma}[6L-4(1-\delta)+2\gamma]}s^{-\frac{l}{2l-\gamma}[4L+4(1-\delta)-2\gamma]}\\&\quad\,\,+\left(\sqrt{K}+\frac{K}{M^{2\delta}}+1\right)s^{-\frac{l}{2l-\gamma}(4\Bbbk-d+2)}\int_{s_0}^s \tau^{-(2L+1+2(1-\delta)(1+\eta))+\frac{l}{2l-\gamma}(4\Bbbk-d+2)}\,\mathrm{d}\tau\\&\lesssim s_0^{-(3L+\gamma)}s^{-(2L-\gamma)}+\left(\sqrt{K}+\frac{K}{M^{2\delta}}+1\right)s^{-(2L+2(1-\delta)(1+\eta))} \le \frac{K}{2}s^{-(2L+(1-\delta)(1+\eta))},
\end{align*}
where in the second inequality we used the fact that 
\begin{align*}
	&\quad\,\,\frac{l}{2l-\gamma}(4\Bbbk-d+2)-2L-2(1-\delta)(1+\eta)\\&=2\left(\frac{2l}{2l-\gamma}-1\right)L+2\left(\frac{2l}{2l-\gamma}-(1+\eta)\right)(1-\delta)-\frac{2l}{2l-\gamma}\gamma>0.
\end{align*}\par
Improved bound for $\mathscr{E}_{2m}$ with $\hbar+2\le m\le l+\hbar$ : Similarly, one integrates (\ref{monotoncity for E2m}) from $s_0$ to $s,$ then applying (\ref{estimate of lambda in s}), $b_1(s)\simeq \frac{c_1}{s},$ Definition \ref{def of initial data} and Definition \ref{bootstrap assump}, and discarding $\frac{1}{\lambda^2}$ in the integral, after some direct computations which we shall omit, one gets $\mathscr{E}_{2m}(s)\le \frac{K}{2}s^{-\frac{l}{2l-\gamma}(4m-d+2)}.$ For the same reason, we omit the proof of improved bound for $\mathscr{E}_{2m}$ with $l+\hbar+1\le m\le \Bbbk-1.$\par
Improved bound for $b_k$ with $l+1\le k\le L$ : We aim to prove (\ref{improved esti for bk from l+1 to L}) by induction on $k.$ When $k=L,$ denoting \begin{equation*}
	\widetilde{b}_L:=b_L+\frac{\langle \mathscr{L}^L q,\chi_{B_0}\Lambda Q\rangle}{\langle \chi_{B_0}\Lambda Q,\Lambda Q\rangle}.
\end{equation*}
Note that by (\ref{estimate of the fraction part}) and (\ref{boootstrap assump on E2Bbbk}), we have \begin{equation}\label{fraction part detailed}
	\left|\frac{\langle \mathscr{L}^L q,\chi_{B_0}\Lambda Q\rangle}{\langle \chi_{B_0}\Lambda Q,\Lambda Q\rangle}\right|\lesssim b_1^{L+\eta(1-\delta)},
\end{equation}
Hence $|\widetilde{b}_L|\lesssim b_1^L.$ Straight calculation gives
\begin{equation*}
	\frac{\mathrm{d}}{\mathrm{d}s}\frac{\widetilde{b}_L(s)}{\lambda(s)^{2L-\gamma}}=\frac{1}{\lambda^{2L-\gamma}}\left[(\widetilde{b}_L)_s+(2L-\gamma)b_1\widetilde{b}_L-(2L-\gamma)\left(\frac{\lambda_s}{\lambda}+b_1\right)\widetilde{b}_L\right].
\end{equation*}
Note that by (\ref{improved estimate for b_L}), (\ref{boootstrap assump on E2Bbbk}) and (\ref{fraction part detailed}), we see that \begin{align*}
	|(\widetilde{b}_L)_s+(2L-\gamma)b_1\widetilde{b}_L|&=\left|(b_L)_s+\partial_s \frac{\langle \mathscr{L}^L q,\chi_{B_0}\Lambda Q\rangle}{\langle \chi_{B_0}\Lambda Q,\Lambda Q\rangle}+(2L-\gamma)b_1b_L+(2L-\gamma)b_1\frac{\langle \mathscr{L}^L q,\chi_{B_0}\Lambda Q\rangle}{\langle \chi_{B_0}\Lambda Q,\Lambda Q\rangle}\right|\\ &\lesssim \frac{1}{B_0^{2\delta}}\left(c(M)\sqrt{\mathscr{E}_{2\Bbbk}}+b_1^{L+1+(1-\delta)-C\eta}\right)+b_1^{1+L+\eta(1-\delta)}\lesssim b_1^{L+1+\eta(1-\delta)}.
\end{align*}
Then by (\ref{modu esti for lambda and bk from 1 to L-1}), we get 
\begin{equation}\label{dynamic for widetildebL}
		\left|\frac{\mathrm{d}}{\mathrm{d}s}\frac{\widetilde{b}_L(s)}{\lambda(s)^{2L-\gamma}}\right|\lesssim \frac{b_1^{L+1+\eta(1-\delta)}}{\lambda^{2L-\gamma}}.
\end{equation}
Then integrating (\ref{dynamic for widetildebL}) from $s_0$ to $s,$ applying (\ref{estimate of lambda in s}), Definition \ref{def of initial data}, $b_1(s)\simeq \frac{c_1}{s}$ and the fact that \begin{equation*}
-(L+1+\eta(1-\delta))+(2L-\gamma)\frac{l}{2l-\gamma}=\frac{l\gamma}{2l-\gamma}\left(L-1-\frac{1}{2}\left(1-\frac{1}{l}\right)\gamma\right)-\eta(1-\delta)-1>-1,
\end{equation*}
we see that 
\begin{align*}
	|\widetilde{b}_L(s)|&\lesssim \lambda(s)^{2L-\gamma}\widetilde{b}_L(s_0)+\lambda(s)^{2L-\gamma}\int_{s_0}^{s}\frac{b_1(\tau)^{L+1+\eta(1-\delta)}}{\lambda(\tau)^{2L-\gamma}}\,\mathrm{d}\tau\\
		&\lesssim s^{-(L-\frac{\gamma}{2})}\cdot s_0^{-\left[\frac{3}{2}L+\frac{1}{2}\gamma-(1-\delta)\right]}+s^{-\frac{l}{2l-\gamma}(2L-\gamma)}\int_{s_0}^{s} \tau^{-(L+1+\eta(1-\delta))+\frac{l}{2l-\gamma}(2L-\gamma)}\,\mathrm{d}\tau\\
		&\lesssim s^{-L-\eta(1-\delta)}\cdot s_1^{\frac{1}{2}\gamma+\eta(1-\delta)}s_0^{-\left[\frac{3}{2}L+\frac{1}{2}\gamma-(1-\delta)\right]}+s^{-L-\eta(1-\delta)}\\
		&\lesssim s^{-L-\eta(1-\delta)},
	\end{align*}
	where in the last inequality, one may assume $s_1$ is a fixed multiplier of $s_0.$ Thus $|b_L|\lesssim |\widetilde{b}_L|+\left|\frac{\langle \mathscr{L}^L q,\chi_{B_0}\Lambda Q\rangle}{\langle \chi_{B_0}\Lambda Q,\Lambda Q\rangle}\right|\lesssim s^{-L-\eta(1-\delta)}.$ Assuming (\ref{improved esti for bk from l+1 to L}) holds for $k+1,$ we aim to show it holds true for $k.$ By induction hypothesis and (\ref{modu esti for lambda and bk from 1 to L-1}), we estimate
\begin{align}
	\frac{\mathrm{d}}{\mathrm{d}s}\frac{b_k(s)}{\lambda^{2k-\gamma}}&=\frac{1}{\lambda^{2k-\gamma}}\left((b_k)_s+(2k-\gamma)b_1b_k-b_{k+1}-(2k-\gamma)\left(\frac{\lambda_s}{\lambda}+b_1\right)b_k+b_{k+1}\right)\notag\\ &\lesssim \frac{1}{\lambda^{2k-\gamma}}(b_1^{L+1+(1-\delta)(1+\eta)}+b_1^{k+L+1+(1-\delta)(1+\eta)}+b_1^{k+1+\eta(1-\delta)})\lesssim \frac{b_1^{k+1+\eta(1-\delta)}}{\lambda^{2k-\gamma}}.\label{dynamic of b_k for l+1lekleL-1}
\end{align}
Similarly, integrating (\ref{dynamic of b_k for l+1lekleL-1}) from $s_0$ to $s$, then applying (\ref{estimate of lambda in s}), Definition \ref{def of initial data} and the fact that \begin{equation*}
	-(k+1+\eta(1-\delta))+\frac{l}{2l-\gamma}(2k-\gamma)=\frac{\gamma}{2l-\gamma}(k-l)-\eta(1-\delta)-1>-1,
\end{equation*}
we obtain
\begin{align*}
	|b_k(s)|&\lesssim \lambda(s)^{2k-\gamma}b_k(s_0)+\lambda(s)^{2k-\gamma}\int_{s_0}^s \frac{b_1(\tau)^{k+1+\eta(1-\delta)}}{\lambda(\tau)^{2k-\gamma}}\,\mathrm{d}\tau\\
	&\lesssim s^{-(k-\frac{\gamma}{2})}s_0^{-(5k-2\gamma)}+s^{-\frac{l}{2l-\gamma}(2k-\gamma)}\int_{s_0}^s \tau^{-(k+1+\eta(1-\delta))+\frac{l}{2l-\gamma}(2k-\gamma)}\,\mathrm{d}\tau\\
	&\lesssim s^{-k-\eta(1-\delta)}\cdot s_1^{\frac{1}{2}\gamma+\eta(1-\delta)}s_0^{-(5k-2\gamma)}+s^{-k-\eta(1-\delta)}\\
	&\lesssim s^{-k-\eta(1-\delta)},
\end{align*}
where again in the last inequality, one may assume $s_1$ is a fixed multiplier of $s_0,$ this consludes the proof of (\ref{improved esti for bk from l+1 to L}). \par 
Improved bound for $\mathcal{V}_1(s)$ : By direct computation, for any $1\le k\le l,$ 
\begin{equation*}
	s(\mathcal{V}_k)_s=\sum\limits_{j=1}^{l-1}(P_l)_{k,j}[s(\mathcal{U}_j)_s-(A_l\mathcal{U})_j]+(P_l)_{k,l}[s(\mathcal{U}_l)_s-(A_l\mathcal{U})_l]+(D_l\mathcal{V})_k.
\end{equation*}
Note that by Lemma \ref{linearization of bk from 1 to l}, Proposition \ref{modulation estimates}, Definition \ref{bootstrap assump} and (\ref{improved esti for bk from l+1 to L}), we get
\begin{align*}
	|s(\mathcal{U}_j)_s-(A_l\mathcal{U})_j|&\lesssim s^{j+i}[(b_j)_s+(2j-\gamma)b_1b_j-b_{j+1}]+O(|\mathcal{U}|^2)\\&\lesssim s^{-\eta(1-\delta)},\\ |s(\mathcal{U}_l)_s-(A_l\mathcal{U})_l|&\lesssim s^{l+1}[(b_l)_s+(2l-\gamma)b_1b_l-b_{l+1}]+s^{l+1}b_{l+1}+O(|\mathcal{U}|^2)\\ &\lesssim s^{-\eta(1-\delta)}.
\end{align*}
It immediately follows that
\begin{equation}\label{dynamic for mathcalV}
	s\mathcal{V}_s=D_l\mathcal{V}+O(s^{-\eta(1-\delta)}).
\end{equation}
In particular,\begin{equation}\label{dynamic for mathcalV1}
	|(s\mathcal{V}_1)_s|\lesssim s^{-\eta(1-\delta)}.
\end{equation}
Integrating (\ref{dynamic for mathcalV1}) from $s_0$ to $s$, then applying Definition \ref{def of initial data}, we see that \begin{equation*}
	\left(\frac{s_0}{s}\right)^{1-\frac{\eta(1-\delta)}{2}}-Cs^{-\frac{\eta(1-\delta)}{2}}\le s^{\frac{\eta(1-\delta)}{2}}\mathcal{V}_1(s)\le \left(\frac{s_0}{s}\right)^{1-\frac{\eta(1-\delta)}{2}}+Cs^{-\frac{\eta(1-\delta)}{2}}
\end{equation*}
for $s_0<s\le s_1.$ This concludes the proof of (\ref{improved esti for V1}).
\end{pf}

Next we describe the reduction to a finite-dimensional problem and transverse crossing property.
\begin{prop}\label{reduction to finite dim and transverse crossing}
	There exists $K_1\ge 1$ such that for any $K\ge K_1$, there exists $s_{0,1}(K)>1$ such that for all $s_0\ge s_{0,1}(K)$ the following holds. Given initial data at $s=s_0$ as in Definition \ref{def of initial data}, if $(b(s),q(s))\in \mathcal{S}_K(s)$ for all $s\in [s_0,s_1]$ with $(b(s),q(s))\in \partial \mathcal{S}_K(s_1)$ for some $s_1\ge s_0,$ then\\
	\textnormal{(\romannumeral1)} Reduction to a finite-dimensional problem: \begin{equation}\label{reduction to a finite dim pb}
		s_1^{\frac{\eta}{2}(1-\delta)}(\mathcal{V}_2(s_1),\cdots,\mathcal{V}_l(s_1))\in \partial \mathcal{B}_{l-1}(0,1).
	\end{equation}
    \textnormal{(\romannumeral2)} Transverse crossing:
    \begin{equation}\label{transverse crossing}
    	\frac{\mathrm{d}}{\mathrm{d}s}\bigg\vert_{s=s_1} \sum\limits_{i=2}^l |s^{\frac{\eta}{2}(1-\delta)}\mathcal{V}_i(s)|^2>0.
    \end{equation}
\end{prop} 
\begin{pf}
	\normalfont
\textnormal{(\romannumeral1)} is a direct consequence of Proposition \ref{improved bootstrap}. Let us prove \textnormal{(\romannumeral2)}. Note that by (\ref{dynamic for mathcalV}), 
\begin{equation*}
	s(\mathcal{V}_i)_s=\frac{i\gamma}{2l-\gamma}\mathcal{V}_i+O(s^{-\eta(1-\delta)}),\,\,\,\text{for}\,\,\,2\le i\le l.
\end{equation*}	
Combined with Definition \ref{bootstrap assump}, we have
\begin{align*}
	\frac{\mathrm{d}}{\mathrm{d}s}\Bigg|_{s=s_1}\sum\limits_{i=2}^l\left|s^{\frac{\eta(1-\delta)}{2}}\mathcal{V}_i(s)\right|^2&=2s^{\eta(1-\delta)-1}\sum\limits_{i=2}^l\left(\frac{\eta(1-\delta)}{2}\mathcal{V}_i^2+s(\mathcal{V}_i)_s\mathcal{V}_i\right)\\ &=2s^{\eta(1-\delta)-1}\left\{\sum\limits_{i=2}^l\left(\frac{i\gamma}{2l-\gamma}+\frac{\eta(1-\delta)}{2}\right)\mathcal{V}_i^2+O(s^{-\frac{3\eta(1-\delta)}{2}})\right\}\\ &\ge \frac{C}{s}\sum\limits_{i=2}^l\left|s^{\frac{\eta(1-\delta)}{2}}\mathcal{V}_i(s)\right|^2\Bigg|_{s=s_1}+O(s^{-\frac{\eta(1-\delta)}{2}-1})\gtrsim \frac{1}{s}>0.
\end{align*}
\end{pf}

Then we are in place to show the existence of solutions to (\ref{eq: renormalized flow}) that are trapped in $\mathcal{S}_K(s)$ for any large $s.$
\begin{prop}\label{existence od sol trapped for large rescaled time}
	There exists $K_2\ge 1$ such that for any $K\ge K_2,$ there exists $s_{0,2}(K)>1$ such that for all $s_0\ge s_{0,2},$ there exists initial data satisfying Definition \ref{def of initial data} such that $(b(s),q(s))\in \mathcal{S}_K(s)$ for all $s\ge s_0.$
\end{prop}
\begin{pf}
	\normalfont
Let us argue by contradiction. If not, defining
\begin{equation*}
	s_{*}:=\sup \{s \mid s\ge s_0\,\,\, \text{such that}\,\,\, (b(s),q(s))\in \mathcal{S}_K(s)\}.
\end{equation*}
as the exit time, then we have for any initial data satisfying Definition \ref{def of initial data}, $s_{*}<+\infty.$ Defining the map \begin{align*}
	\Xi\,\colon \,\mathcal{B}_{l-1}(0,1)&\longrightarrow \partial \mathcal{B}_{l-1}(0,1)\\ s_0^{\frac{\eta(1-\delta)}{2}}(\mathcal{V}_2(s_0),\cdots,\mathcal{V}_l(s_0))&\longmapsto  s_{*}^{\frac{\eta(1-\delta)}{2}}(\mathcal{V}_2(s_{*}),\cdots,\mathcal{V}_l(s_{*})).
\end{align*}
By \textnormal{(\romannumeral1)} of Propostion \ref{reduction to finite dim and transverse crossing}, $\Xi$ is well defined. By \textnormal{(\romannumeral2)} of Propostion \ref{reduction to finite dim and transverse crossing}, $\Xi$ restricted on ${\partial \mathcal{B}_{l-1}(0,1)}$ is the identity map. Then $\Xi$ is a continuous map and an identity on the boundry of a ball which is not possible in view of the Brouwer's fixed point theorem, this concludes the proof.
\end{pf}

Finally we finish the proof of Theorem \ref{main thm}.\par
Expression of $\lambda$ in original time variable. Recall by the proof of (\ref{estimate of lambda in s}), we get
\begin{equation*}
	-\lambda_t\lambda=c(u_0)\lambda^{\frac{2l-\gamma}{l}}(1+o(1)),
\end{equation*}
or say
\begin{equation*}
	\partial_t\left(\lambda^{\frac{\gamma}{l}}\right)=-c(u_0)(1+o(1)).
\end{equation*}
Then integrating from $t$ to $T,$ we see that
\begin{equation*}
	\lambda(t)=c(u_0)(1+o(1))(T-t)^{\frac{l}{\gamma}}.
\end{equation*}\par 
On the smallness of Sobolev norms of $q.$ By \textnormal{(\romannumeral3)} of Lemma \ref{coercivity-determined esti on q} and Definition \ref{bootstrap assump}, we have
\begin{equation*}
	\int |\partial_y^{2m}q|^2\lesssim \mathscr{E}_{2m}\rightarrow 0,\,\,\,\text{as}\,\,\,s\rightarrow \infty, \,\,\,\text{for}\,\,\,\hbar+2\le m\le \Bbbk.
\end{equation*}
This concludes the whole proof of Theorem \ref{main thm}. 
\vspace{\baselineskip}

\section*{\centerline{Acknowledgement}}
The author was supported by NSFC Grant of China No. 12271497, No. 12341102 and the National Key Research and Development Program of China No. 2020YFA0713100.
\vspace{\baselineskip}
\centering

\end{document}